\pdfoutput=1
\RequirePackage{ifpdf}
\ifpdf 
\documentclass[pdftex]{sigma}
\else
\documentclass{sigma}
\fi

\newsavebox{\largestimage}

\usepackage{bbm,mathtools, blkarray, enumitem, tikz-cd}
\usepackage[scr=boondox]{mathalpha}
\usepackage{subcaption}
\definecolor{dark-red}{rgb}{0.75,0.2,0.2}
\definecolor{dark-blue}{rgb}{0.25,0.15,0.55}
\definecolor{medium-blue}{rgb}{0.10,0.10,0.85}
\definecolor{ForestGreen}{rgb}{.133, .545, .133}

\numberwithin{equation}{section}

\newtheorem{Theorem}{Theorem}[section]
\newtheorem*{Theorem*}{Theorem}
\newtheorem{Corollary}[Theorem]{Corollary}
\newtheorem{Lemma}[Theorem]{Lemma}
\newtheorem{Proposition}[Theorem]{Proposition}
\newtheorem{Conjecture}[Theorem]{Conjecture}

\theoremstyle{definition}
\newtheorem{Definition}[Theorem]{Definition}

\newtheorem{Example}[Theorem]{Example}
\newtheorem{Remark}[Theorem]{Remark}
\newtheorem{Question}[Theorem]{Question}
\newtheorem{Notation}[Theorem]{Notation}

\usetikzlibrary{positioning}
\usetikzlibrary{calc}
\usetikzlibrary{decorations.markings}
\usetikzlibrary{arrows.meta}
\tikzset{> /.tip = {Stealth[round,length=7pt]}}
\tikzstyle over=[preaction={draw,line width=5pt,white}]
\usetikzlibrary{calc}
\usetikzlibrary{decorations.markings}

\newcommand{\qbin}[3]{\begin{bmatrix} \vspace{-0.7mm} \raise 1mm \hbox{$#1$}\vspace{-0.7mm}\\ #2\vspace{-0.8mm}\end{bmatrix}_{#3}}

\newcommand{\bracks}[1]{\langle#1\rangle}
\newcommand{\floor}[1]{\lfloor #1 \rfloor}
\newcommand{\mvec}[1]{{\smash{\text{\boldmath{$#1$}}}}}
\newcommand{\s}{\mvec{s}}
\newcommand{\x}{\mvec{\sigma^{(100)}}}
\newcommand{\y}{\mvec{\sigma^{(001)}}}
\newcommand{\Vt}{V(\mvec{t})}

\newcommand{\Vs}{V(\mvec{s})}
\newcommand{\VH}{V^H(\mvec{\lambda})}

\newcommand{\grf}[1]{
 \mathchoice{\raisebox{-1pt}{\scalebox{1.18}{$\mathsf{#1}$}}}
 {\raisebox{-1pt}{\scalebox{1.18}{$\mathsf{#1}$}}}
 {\raisebox{-0.15pt}{\scalebox{0.8}{$\mathsf{#1}$}}}
 {\raisebox{-0.15pt}{\scalebox{0.8}{$\mathsf{#1}$}}}}
\newcommand{\sroots}{\grf{\Delta}}
\newcommand{\roots}{\grf{\Phi}}
\newcommand{\psroots}{{\sroots^+}}
\newcommand{\proots}{{\roots^+}}
\newcommand{\maps}{\grf{\Psi}}

\newcommand{\Apoly}{\Delta_{\mathcal{A}}}
\newcommand{\R}{\mathscr{R}}
\renewcommand{\H}{\X_{12}}
\newcommand{\X}{\mathscr{X}}
\newcommand{\tr}{\operatorname{tr}}
\newcommand{\trace}{\operatorname{trace}}
\newcommand{\mdim}{\mathsf{mdim}}
\newcommand{\Hom}{\operatorname{Hom}}
\newcommand{\End}{\operatorname{End}}
\newcommand{\Uwmod}{\mathcal{C}}
\newcommand{\UHwmod}{\mathcal{C}^H}
\newcommand{\braid}{\boldsymbol{c}^H}
\newcommand{\redbraid}{\boldsymbol{c}}
\newcommand{\sbraid}{{\overline{\boldsymbol{c}}}^H}
\renewcommand{\t}{\mvec{t}}
\renewcommand{\P}{\mathbb{T}}
\newcommand{\E}{\Upsilon}
\newcommand{\qRchar}{R}
\newcommand{\Rt}{\qRchar^\bullet}
\newcommand{\lcoev}{\stackrel{\longleftarrow}{\operatorname{coev}}}
\newcommand{\lev}{\stackrel{\longleftarrow}{\operatorname{ev}}}
\newcommand{\rev}{\stackrel{\longrightarrow}{\operatorname{ev}}}
\newcommand{\rcoev}{\stackrel{\longrightarrow}{\operatorname{coev}}}

\usepackage[labelfont=bf,labelsep=period]{caption}

\begin{document}

\allowdisplaybreaks

\newcommand{\arXivNumber}{2008.06983}

\renewcommand{\PaperNumber}{025}

\FirstPageHeading

\ShortArticleName{A Non-Abelian Generalization of the Alexander Polynomial from Quantum $\mathfrak{sl}_3$}

\ArticleName{A Non-Abelian Generalization \\ of the Alexander Polynomial from Quantum $\boldsymbol{\mathfrak{sl}_3}$}

\Author{Matthew HARPER}

\AuthorNameForHeading{M.~Harper}

\Address{Department of Mathematics, Michigan State University, East Lansing, MI 48824, USA}
\Email{\mail{mrhmath@proton.me}, \mail{harpe111@msu.edu}}
\URLaddress{\url{https://mrhmath.github.io/}}

\ArticleDates{Received July 28, 2025, in final form February 18, 2026; Published online March 17, 2026}

\Abstract{One construction of the Alexander polynomial is as a quantum invariant associated with representations of restricted quantum $\mathfrak{sl}_2$ at a fourth root of unity. We generalize this construction to define a link invariant $\Delta_{\mathfrak{g}}$ for any semisimple Lie algebra $\mathfrak{g}$ of rank $n$, taking values in $n$-variable Laurent polynomials. Focusing on the case $\mathfrak{g}=\mathfrak{sl}_3$, we establish a direct relation between $\Delta_{\mathfrak{sl}_3}$ and the Alexander polynomial. We show that certain parameter evaluations of $\Delta_{\mathfrak{sl}_3}$ recover the Alexander polynomial on knots, despite the $R$-matrix not satisfying the Alexander--Conway skein relation at these points. We tabulate $\Delta_{\mathfrak{sl}_3}$ for all knots up to seven crossings and various other examples, including the Kinoshita--Terasaka knot and Conway knot mutant pair which are distinguished by this invariant.}


\Classification{57K16; 17B37}

\section{Introduction}

 \subsection{An overview of quantum group invariants} One of the goals of quantum topology is to construct combinatorial and algorithmically computable invariants of knots and 3-manifolds with significant implications for low-dimensional topology. Given a representation of a quantum group, the Reshetikhin--Turaev construction produces an invariant of links \cite{RT}. The most well-known of these invariants is the Jones polynomial, obtained from the fundamental representation of $U_q(\mathfrak{sl}_2)$ \cite{Jones}. Other representations of~$U_q(\mathfrak{sl}_2)$ define the so-called colored-Jones polynomials which are related to the Jones polynomial of cablings of knots and higher-dimensional representations. Higher rank analogs of the Jones polynomial are computed from representations of the quantum groups $U_q(\mathfrak{g})$, where~$\mathfrak{g}$ is a simple Lie algebra. These type-$\mathfrak{g}$ invariants include specializations of the {HOMFLY} and Kauffman polynomials, and the Kuperberg polynomial \cite{HOMFLY, KauffmanIsotopyInvariant, KuperbergG2}.

 The Alexander polynomial, an invariant from classical topology, also arises as a quantum invariant from a family of representations of \smash{$U_{\sqrt{-1}}(\mathfrak{sl}_2)$} with a slight modification to the construction \cite{Murakami, MurakamiStateModel, Ohtsuki}. If $\omega$ is a primitive $n$-th root of unity, then $U_\omega(\mathfrak{sl}_2)$ admits a family of $n/\!\gcd(n,2)$-dimensional representations $V(t)$, with arbitrary nonzero highest weight~$t$. The associated invariants are called the ADO invariants \cite{ADO} and include the Alexander polynomial in the case $\omega=\sqrt{-1}$.

 In the present paper, we initiate the study of higher rank Lie type analogs of the Alexander polynomial associated with representations of $U_{\sqrt{-1}}(\mathfrak{g})$ which have arbitrary nonzero highest weights and are denoted here by $\Delta_{\mathfrak{g}}$. We focus on the case $\mathfrak{g}=\mathfrak{sl}_3$, which is the simplest generalization in terms of algebraic complexity and appears to have the most direct classical topological relevance. In contrast to the invariants described in the first paragraph, these polynomials belong to the set of $n$-variable Laurent polynomials, where $n$ is the rank of $\mathfrak{g}$. More generally, one may consider invariants $\Delta_{\mathfrak{g},\omega}$ associated to representations of $U_\omega(\mathfrak{g})$ at roots of unity. The well-definedness of these invariants has been shown for roots of unity $\omega$ with odd order at least three in \cite{GPtrace}, but the invariants themselves have not been computed explicitly and additional properties of these (non-super) invariants beyond rank one are not known. We~summarize this invariantology in Table~\ref{tab:invariantology} below.

 One other family of invariants worth mentioning here are associated to quantum supergroups $\mathfrak{gl}(m|n)$ (or $\mathfrak{sl}(m|n)$) at generic $q$. The Alexander polynomial appears among these invariants, derived from representations of $\mathfrak{gl}(1|1)$ \cite{KauffmanSaleur,Sartori,Viro}. In higher rank, the Links--Gould invariants are polynomials in two variables \cite{DeWitinfinite,GPsl21,LinksGould} which admit specializations to a product of Alexander polynomials or the Alexander polynomial in the variable $t^2$ \cite{DeWitIshiiLinks,KohliPatureau}. The Links--Gould invariants are known to improve on the genus bound determined by the Alexander polynomial \cite{Kohli-Tahar, NvdV}.

 \begin{table}[!ht]
 \centering\centering\renewcommand{\arraystretch}{1.2}
 \begin{tabular}{c|c}
 $\substack{\mbox{generic $q$}\\\mbox{(polynomials in $q$)}}$ &$\substack{\mbox{$\omega$ is a root of unity}\\\mbox{(polynomials in $t_1,\dots, t_n$)}}$\\
 \hline
 $\big(U_q(\mathfrak{sl}_2),V_2\big)$ Jones polynomial
 &
 $\big(U_{\sqrt{-1}}(\mathfrak{sl}_2),V(t)\big)$ Alexander polynomial
 \\[1ex]
 $\big(U_q(\mathfrak{sl}_2),V_m\big)$ colored-Jones polynomial
 &
 $\big(U_{\omega}(\mathfrak{sl}_2),V(t)\big)$ ADO invariants
 \\[1ex]
 $\big(U_q(\mathfrak{g}),V_n\big)$ type-$\mathfrak{g}$ polynomial
 &
 $\big(U_{\omega}(\mathfrak{g}),V(t_1,\dots, t_n)\big)$ higher Alexander/ADO
 \end{tabular}
 \caption{Some link polynomials from non-super quantum groups.}
 \label{tab:invariantology}
 \end{table}

 Despite both being derived from non-semisimple categories, invariants from quantum groups at roots of unity are qualitatively different from the Links--Gould invariants in that they can have more than two variables. Moreover, while cabling of knots for the invariants at generic $q$ produces ``colored'' invariants, associated to higher-dimensional representations, tensor products of the representations in the root of unity case are ``self-similar'' and do not provide any significant refinement.

 The root of unity link invariants and the Links--Gould polynomials do share another feature which contrasts them against the type-$\mathfrak{g}$ polynomials derived from semisimple representation categories. These non-semisimple invariants have quantum dimension zero, implying that the naive RT invariant assigns the value of zero to any closed tangle. To compute meaningful invariants from these ``negligible'' objects, we use the modified trace construction formalized by Geer, Patureau-Mirand, and Turaev \cite{GPT}.

 The introduction of these link polynomials $\Delta_\mathfrak{g}$ leads to exciting questions about which properties they share with and refine from the Alexander polynomial, their topological implications, and relations to other invariants.

 \subsection{Main results} We consider the restricted quantum group $\overline{U}_\zeta(\mathfrak{g})$ associated to a simple Lie algebra $\mathfrak{g}$ of rank $n$ at a primitive fourth root of unity $\zeta$. This quantum group is the quotient of $U_\zeta(\mathfrak{g})$ by the Hopf ideal generated by the square of all root generators.

 Let $\proots$ ($\psroots$) denote a choice of positive (simple) roots for the root system of $\mathfrak{g}$. Each character $\t$ on the Cartan subalgebra, which we identify with $(t_1,\dots,t_n)\in (\mathbb{C}^\times)^n\cong \text{Map}(\psroots,\mathbb{C}^\times)$, determines a Verma module $\Vt$ of dimension $2^{|\proots|}$ over $\overline{U}_\zeta(\mathfrak{g})$ \cite{KDC90}. To such a quantum group with a family of representations, we denote the associated (modified) Turaev $R$-matrix invariant \cite{GPT,RT, Turaev88} by $\Delta_{\mathfrak{g}}$, which assigns a Laurent polynomial in $\mathbb{C}\bigl[ t_1^{\pm1},\dots, t_n^{\pm1}\bigr]$ to every link $\mathcal{L}$.
	
	These invariants are not to be confused with the multivariable Alexander polynomial. The number of variables in $\Delta_{\mathfrak{g}}$ depends on the rank of~$\mathfrak{g}$ and not on the number of components of~$\mathcal{L}$. If $\mathcal{L}$ has $m$ components, one can consider a ``multi-colored'' version of~$\Delta_{\mathfrak{g}}$ in which each component of $\mathcal{L}$ is assigned a representation with a different highest weight, but we do not investigate this generalization in detail here.

 We give particular attention to $\Delta_{\mathfrak{sl}_3}$ which is a two-variable Laurent polynomial invariant of links.

\begin{Theorem} \label{Thm:list}
 The invariant $\Delta_{\mathfrak{sl}_3}$ has the following properties:
 \begin{enumerate}[label=\textup{(\arabic*)}]\itemsep=0pt
 \item it dominates the Alexander polynomial on knots $($Theorem~{\rm \ref{thm:plugin})},
 \item \label{thm:sym} for all links $\mathcal{L}$ $($Section~{\rm \ref{sec:sym})},
 \[
 \Delta_{\mathfrak{sl}_3}(\mathcal{L})(t_1,t_2)=\Delta_{\mathfrak{sl}_3}(\mathcal{L})(t_2,t_1) =\Delta_{\mathfrak{sl}_3}(\mathcal{L})\bigl(t_1^{-1},t_2^{-1}\bigr)\in \mathbb{Z}\bigl[t_1^{\pm2}, t_2^{\pm2}\bigr],
 \]

 \item \label{Thm:nonabelian} it can detect mutation and knots with zero Alexander module, and is therefore non-abelian in the sense of~{\rm \cite{Cochran}}
 $($Figures {\rm \ref{fig:mutants}} and {\rm \ref{fig:doubles})},
 \item there is a $9$-term skein relation for $\Delta_{\mathfrak{sl}_3}$ \label{Thm:skein}
 $($Proposition~{\rm \ref{prop:skein})}.
 \end{enumerate}
\end{Theorem}

 The dominance of $\Delta_{\mathfrak{sl}_3}$ over the Alexander polynomial $\Apoly$ is implied by the following theorem, which shows that $\Delta_{\mathfrak{sl}_3}$ is a generalization of the classical knot invariant.

\begin{Theorem}[Theorem~\ref{thm:plugin}]\label{thm:introplugin}
 Let $\mathcal{K}$ be any knot. Then
		\begin{equation*}
		\Delta_{\mathfrak{sl}_3}(\mathcal{K})(t,\pm 1)=\Delta_{\mathfrak{sl}_3}(\mathcal{K})(\pm 1,t)=\Delta_{\mathfrak{sl}_3}(\mathcal{K})\bigl(t,\pm \sqrt{-1}/t\bigr)=\Apoly(\mathcal{K})\bigl(t^4\bigr).
		\end{equation*}
		Moreover, if for all knots $\mathcal{K}$ the equality $\Delta_{\mathfrak{sl}_3}(\mathcal{K})(t,s)=\Apoly(\mathcal{K})\bigl(t^4\bigr)$ holds, then either $s^2=1$ or $s^2=-t^{-2}$.
\end{Theorem}

	This equality of invariants is not obvious. The rank one relation $\Delta_{\mathfrak{sl}_2}(\mathcal{L})(t)=\Apoly(\mathcal{L})\bigl(t^2\bigr)$ for any link~$\mathcal{L}$ is straightforward to prove from the minimal polynomial of the $\mathfrak{sl}_2$ $R$-matrix because it satisfies the Alexander--Conway skein relation \cite{Murakami,MurakamiStateModel,Ohtsuki}. In contrast, the $R$-matrix evaluated at $t_2=1$ (for example) in the $\mathfrak{sl}_3$ case does not satisfy this skein relation, but nevertheless yields the Alexander polynomial. Consequently, the tangle invariant obtained from the evaluated $\mathfrak{sl}_3$ $R$-matrix is different from the Alexander ($\mathfrak{sl}_2$) tangle invariant. Given that $\Delta_{\mathfrak{sl}_3}$ and the Alexander polynomial (as a quantum invariant) is computed here using a partial trace, Theorem~\ref{thm:introplugin} is the statement that the two tangle invariants agree on long-knots.

 The parameter evaluations of Theorem~\ref{thm:introplugin} are natural from a representation-theoretic point of view. To each $\alpha\in\proots$ we associate a curve in $\X_\alpha\subset (\mathbb{C}^\times)^2$, see Figure~\ref{fig:Xa}. A point $\t=(t_1,t_2)$ on exactly one such curve determines an evaluation of $\Delta_{\mathfrak{sl}_3}$ to the Alexander polynomial as presented in Theorem~\ref{thm:introplugin}. Let $\R_\alpha$ be the set of points in $\X_\alpha$ which are disjoint from some other $\X_\beta$. Then $\R_\alpha$ parameterizes the highest weights~$\t$ such that~$\Vt$ is reducible with a~four-dimensional (irreducible) head $W_\alpha(\t)$. If $\t\in\R_\alpha$, then the knot invariants derived from~$\Vt$ and $W_\alpha(\t)$ are the same. Theorem~\ref{thm:introplugin} is now proven as a consequence of the following.

 \begin{Theorem} [Theorem~\ref{thm:smallAlexander}]\label{thm:introsmallAlexander}
 Fix a positive root $\alpha$. Assume $\t\in\R_\alpha$ so that it is of the form $(\sigma,t)$, $(t,\sigma)$, or $\bigl(t,\sigma\sqrt{-1}/ t\bigr)$ where $\sigma^2=1$ and $t\in\mathbb{C}^\times$. The $R$-matrix invariant of a link colored by $W_\alpha(\t)$ is equal to the Alexander--Conway polynomial evaluated at~$t^4$.
 \end{Theorem}

 \begin{figure}[!ht]
 \centering
 \begin{tikzpicture}[scale=1.1]
 \draw[gray!50, thin, step=0.5] (-1.5,-1.5) grid (1.5,1.5);
 \draw[very thick,->] (-1.5,0) -- (1.7,0) node[right] {$t_1^2$};
 \draw[very thick,->] (0,-1.5) -- (0,1.7) node[right] {$t_2^2$};

 \draw[medium-blue, very thick] (.5,-1.5) -- (.5,1.5);

 \draw[dark-red, very thick] (-1.5,.5) -- (1.5,.5);

 \draw[very thick, domain=0.165:1.5, smooth, variable=\x, ForestGreen] plot ({\x}, {-.25/\x});

 \draw[very thick, domain=0.165:1.5, smooth, variable=\x, ForestGreen] plot (-{\x}, {.25/\x});

 \draw (.5,1.5) node[ right] {
 {\color{medium-blue}\footnotesize$\X_1$}};

 \draw (1.5,.5) node[ above] {
 {\color{dark-red}\footnotesize$\X_2$}};

 \draw (1.2,-.5) node[] {
 {\color{ForestGreen}\footnotesize$\X_{12}$}};
 \draw (-.5,1.2) node[] {
 {\color{ForestGreen}\footnotesize$\X_{12}$}};

 \foreach \x in {.5} \draw (\x,0.05) -- (\x,-0.05);
 \foreach \y in {.5} \draw (-0.05,\y) -- (0.05,\y);

 \draw (.45,0) node[below]
 {\small$1$};

 \draw (0,.45) node[left ]
 {\small$1$};
\end{tikzpicture}
 \caption{Sketch of the curves $\X_\alpha\subset (\mathbb{C}^\times)^2$:
$\X_1=\bigl\{(t_1,t_2)\mid t_1^2=1\bigr\}$, $\X_2=\bigl\{(t_1,t_2)\mid t_2^2=1\bigr\}$, $\X_{12}=\bigl\{(t_1,t_2)\mid (t_1t_2)^2=-1\bigr\}$.
 Each point on a unique $\X_\alpha$ determines an evaluation to the Alexander polynomial and is a highest weight of $\Vt$ with irreducible subrepresentation $W_\alpha(\t)$.}
 \label{fig:Xa}
 \end{figure}

 \subsection{Tabulation of the invariant} We include a tabulation of $\Delta_{\mathfrak{sl}_3}$ on all prime knots up to seven crossings in Figure~\ref{fig:sl3leq7} as well as several other knots in Figure~\ref{fig:sl3geq8}. Most notable among them is the Conway knot $\mathsf{11_{n34}}$ and the Kinoshita--Terasaka knot $\mathsf{11_{n42}}$ which are a mutant pair and are distinguished by $\Delta_{\mathfrak{sl}_3}$. The values of $\Delta_{\mathfrak{sl}_3}$ on $\mathsf{11_{n34}}$ and $\mathsf{11_{n42}}$ are determined from its coefficients in Figure~\ref{fig:mutants}. The constant term of the polynomial is boxed and establishes the coordinate $(0,0)$. The coefficient of $t_1^{2a}t_2^{2b}$ is in position $(a,b)$ and by Theorem~\ref{Thm:list}\,\ref{thm:sym} is equal to those in positions $(b,a)$, $(-a,-b)$, and $(-b,-a)$, which are not shown.

 In addition to $\mathsf{11_{n34}}$ and $\mathsf{11_{n42}}$, untwisted Whitehead doubles of knots have trivial Alexander module and Alexander polynomial equal to 1 \cite{Rolfsen}. Recall that the Alexander polynomial is an abelian knot invariant in the sense of \cite{Cochran} in that it is determined by the first two terms of the derived series of the knot group, whereas the Jones polynomial is non-abelian. Abelian invariants are limited in their ability to distinguish knots with a trivial Alexander module, such as untwisted Whitehead doubles, from the unknot. We find that~$\Delta_{\mathfrak{sl}_3}$ is nontrivial on the Whitehead double of the trefoil $\text{Wh}^0(\mathsf{3_{1}})$, see Figure~\ref{fig:doubles}. Thus proving Theorem~\ref{Thm:list}\,\ref{Thm:nonabelian}.

	{	\begin{figure}[!ht]\centering
			\begin{subfigure}{.35\linewidth}
				\centering
				\includegraphics[]{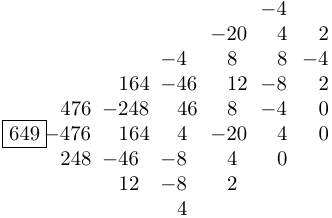}
				\caption{$\mathsf{11_{n34}}$}
			\end{subfigure}\qquad
			\begin{subfigure}{.35\linewidth}
				\centering
				\vspace{8.4mm}
				\includegraphics[]{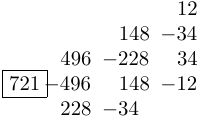}
				\vspace{8.2mm}
				\caption{$\mathsf{11_{n42}}$}
			\end{subfigure}
			\caption{The values of $\Delta_{\mathfrak{sl}_3}$ on the mutant pair $\mathsf{11_{n34}}$ and $\mathsf{11_{n42}}$.
			}\label{fig:mutants}
	\end{figure}}
	{	\begin{figure}[!ht]
			\begin{subfigure}{\linewidth}
				\centering
				\includegraphics[]{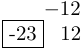}
				\caption{$\text{Wh}^0(\mathsf{3_{1}})$}
			\end{subfigure}
			\caption{The value of $\Delta_{\mathfrak{sl}_3}$ on the untwisted Whitehead double of $\mathsf{3_{1}}$.
			}\label{fig:doubles}
	\end{figure}}
 A limitation in our computation of $\Delta_{\mathfrak{sl}_3}$ is that we compute it from braid presentations of knots. A knot with braid index $k$ requires the multiplication of $8^k \times 8^k$ sparse symbolic matrices. Under a simple implementation, it took about 24 hours to compute each of $\Delta_{\mathfrak{sl}_3}(\mathsf{11_{n34}})$ and $\Delta_{\mathfrak{sl}_3}(\mathsf{11_{n42}})$ in \textsc{Maple 2018.0} on The Ohio State University's unity high performance computing cluster using the $k=4$ presentations on~\cite{TKA}. All invariants in this paper can be computed using \textsc{SymPy~1.14.0} with the domainmatrix module in a few hours. Perhaps these computations could be made more efficient by implementing the methods of \cite{BNV:PG,BNV:API}.

	The $\mathfrak{sl}_3$ invariant admits a nine-term skein relation via the minimal polynomial of the $R$-matrix represented in $\Vt\otimes \Vt$, see Theorem~\ref{Thm:list}\,\ref{Thm:skein} and Proposition~\ref{prop:skein}. Using a recursion determined by the square of the $R$-matrix, we compute an explicit formula for $\Delta_{\mathfrak{sl}_3}$ on $(2n+1,2)$ torus knots.

\begin{Theorem}
 [Theorem~\ref{thm:torusknots}]\label{thm:introtorusknots}
 The value of $\Delta_{\mathfrak{sl}_3}$ on a $(2n+1,2)$ torus knot is given by
 \begin{align*}
		\frac{
			\bigl(t_1-t_1^{-1}\bigr)\bigl(t_1^{4n+2}+t_1^{-(4n+2)}\bigr)} {\bigl(t_2+t_2^{-1}\bigr)\bigl(t_1^2+t_1^{-2}\bigr)\bigl(t_1t_2-t_1^{-1}t_2^{-1}\bigr)}&+
		\frac{	\bigl(t_2-t_2^{-1}\bigr) \bigl(t_2^{4n+2}+t_2^{-(4n+2)}\bigr)}{\bigl(t_1+t_1^{-1}\bigr)\bigl(t_2^2+t_2^{-2}\bigr)\bigl(t_1t_2-t_1^{-1}t_2^{-1}\bigr)}\\
&+
	\frac{\bigl(t_1t_2+t_1^{-1}t_2^{-1}\bigr) \bigl(t_1^{4n+2}t_2^{4n+2}+t_1^{-(4n+2)}t_2^{-(4n+2)}\bigr)} {\bigl(t_1^2t_2^2+t_1^{-2}t_2^{-2}\bigr)\bigl(t_1+t_1^{-1}\bigr)\bigl(t_2+t_2^{-1}\bigr)}.
		\end{align*}
\end{Theorem}

	\subsection{Relation to other invariants}

Non-semisimple quantum invariants from the quantum supergroups $\mathfrak{gl}(m|n)$ (or $\mathfrak{sl}(m|n)$) are also known to generalize the Alexander polynomial. The representations used in the construction of the Links--Gould invariants $LG^{m,n}\in\mathbb{Z}[t,q]$ also have an arbitrary highest weight taking the role of the polynomial variable, however it is not necessary that~$q$ be a root of unity to define such a~representation. The relation between the Links--Gould invariants and the Alexander polynomial
\begin{align*}
	LG^{m,n}(\mathcal{L})\bigl(t, {\rm e}^{{\rm i}\pi/m}\bigr)=\bigl(\Apoly(\mathcal{L})\bigl(t^{2m}\bigr)\bigr)^n
\end{align*}
given by specializing $q$ to be a $2m$-th root of unity is conjectured to hold for all~$m$ and~$n$, and has been proven when either $m$ or $n$ equals~1 \cite{DeWitIshiiLinks,KohliPatureau}. Compare this with Conjecture~\ref{conj} below. We also note that there does not appear to be an evaluation of~$\Delta_{\mathfrak{sl}_3}$ which equals a~higher power of the Alexander polynomial. This can be checked by solving for $t_2$ in the system $\Delta_{\mathfrak{sl}_3}(\mathcal{K})(t_1,t_2)=\Apoly(\mathcal{K})(t_1^m)^n$ for $\mathcal{K}\in \{\mathsf{3_1}, \mathsf{4_1}\}$. The system is further simplified by assuming $t_1=2$, for example, and it has been verified that there is no solution for positive integers $m,n\leq 20$ except $(m,n)=(4,1)$.

 The low rank invariants $V_1$ and $\Lambda_{-1}$ of knots constructed in~\cite{GK} were conjectured to coincide with RT polynomials $LG$ and $\Delta_{\mathfrak{sl}_3}$ of links, and this was proven affirmatively in~\cite{GHKST}. Specifically, $\Lambda_{-1}$ extends to a link invariant valued in $\mathbb{Z}\bigl[t_1^{\pm1},t_2^{\pm1}\bigr]$ and satisfies the relation
 \[
 \Lambda_{-1}(\mathcal{L})\bigl(t_1^{-2}, t_2^{-2}\bigr) = \Delta_{\mathfrak{sl}_3}(\mathcal{L})(t_1,t_2).
 \]
 The $R$-matrix used in the construction of the invariants of Garoufalidis--Kashaev is natural from the topological perspective in that its matrix entries are valued in $\mathbb{Z}[t^\pm,s^\pm]$ and the link invariant is defined in canonical (non-squared) variables. The results of \cite{GK} and \cite{GHKST} imply Theorem~\ref{Thm:list}\,\ref{thm:sym} as well as the symmetry $\Delta_{\mathfrak{sl}_3}(t_1,t_2)=\Delta_{\mathfrak{sl}_3}\bigl(\frac{\zeta}{t_1t_2},t_2\bigr)$. We include a self-contained proof of the first equality of Theorem~\ref{Thm:list}\,\ref{thm:sym} in Proposition~\ref{cor:symm} from the perspective of Dynkin diagram automorphisms.

 \subsection{Further questions}
	Here we propose additional conjectures regarding the properties of the invariants $\Delta_{\mathfrak{g}}$ and the representations studied in this paper.
	Following \cite{Piccirillo}, since $\Delta_{\mathfrak{sl}_3}$ distinguishes $\mathsf{11_{n34}}$ and $\mathsf{11_{n42}}$, it~is natural to pose the following question.
 \begin{Question}
 Does $\Delta_{\mathfrak{sl}_3}$ contain information on sliceness, such as a generalized Fox--Milnor condition?
 \end{Question}
 Nevertheless, we suspect $\Delta_{\mathfrak{sl}_3}$ is related to other geometrically constructed invariants that are sensitive to knots with trivial Alexander modules. Knot Floer homology, for example, is nontrivial on the Whitehead double of $\mathsf{4_1}$ \cite{Hedden}. Another example is the set of twisted Alexander polynomials for a particular matrix group. The set of twisted invariants derived from all parabolic ${\rm SL}_2(\mathbb{F}_7)$ representations, up to conjugacy, of the knot groups of~$\mathsf{11_{n34}}$ and~$\mathsf{11_{n42}}$ are enough to distinguish the pair of mutant knots from each other and the unknot~\cite{Wada}.
	
 Higher Alexander modules \cite{Cochran}, which use terms further in the derived series of the knot group, improve the Alexander polynomial genus bound and can detect mutation on knots with nontrivial Alexander polynomial \cite{Horn}. However, these modules are trivial on knots with Alexander polynomial 1, such as $\mathsf{11_{n34}}$ and $\mathsf{11_{n42}}$.

 Theorem~\ref{thm:introplugin} (Theorem~\ref{thm:plugin}) may be stated in terms of $\Delta_{\mathfrak{sl}_3}$ and $\Delta_{\mathfrak{sl}_2}$, since for any link $\mathcal{L}$ $\Delta_{\mathfrak{sl}_2}(\mathcal{L})(t)=\Apoly(\mathcal{L})\bigl(t^2\bigr)$. This motivates the conjecture that the set of invariants $\Delta_\mathfrak{g}$ indexed by $\mathfrak{g}$ are partially ordered according to dominance, and in this ordering $\Delta_{\mathfrak{g}'}\leq\Delta_\mathfrak{g}$ if and only if $\mathfrak{g}'\subseteq \mathfrak{g}$. Let $\mathfrak{g}$ be a simple Lie algebra of rank $n+1$. Suppose $\Vt$ is a representation of $\overline{U}_\zeta(\mathfrak{g})$ with highest weight $\t\in (\mathbb{C}^\times)^{n+1}$ and dimension $2^{n+1}$. We denote by $v_{{\rm lowest}}^{\t}$ a lowest weight vector of $\Vt$.
	\begin{Conjecture}\label{conj}
		Choose $\t\in (\mathbb{C}^\times)^{n+1}$ such that for exactly one $\alpha\in\proots$ and all $\alpha_i\in\psroots$, $E_iE_\alpha v_{{\rm lowest}}^{\t}=0$ in $\Vt$. Suppose $\alpha_j$ has a nonzero component in $\alpha$. Let $\t'$ be obtained by deleting the $j$-th entry from $\t$. Then for any knot $\mathcal{K}$, $\Delta_{\mathfrak{sl}_{n+1}}(\mathcal{K})(\t)=\Delta_{\mathfrak{sl}_{n}}(\mathcal{K})\bigl((\t')^{2(n-1)}\bigr)$.
	\end{Conjecture}
 It is also natural to investigate Conjecture \ref{conj} on quantum groups in other Lie types and at other roots of unity, and how it extends to the case of small roots of unity for non-simply laced types as studied in \cite{Lentner}.

	In Theorem~\ref{thm:smallAlexander}, we prove that for each positive root $\alpha$ the family of four-dimensional $\overline{U}_\zeta(\mathfrak{sl}_3)$ representations $W_\alpha(\t)$ determines the Alexander polynomial of singly-colored links. For each family of representations, we claim that the relations for the Conway potential function, given in \cite{Jiang}, are also satisfied.
	\begin{Conjecture}
		The multi-variable invariant of links with components colored by the palette
		$\{W_\alpha(\t)\mid\t\in\R_\alpha\}$ for each $\alpha\in\proots$ is the Conway potential function.
	\end{Conjecture}
	
 There is a natural identification between the Burau representation and the braid representations from $R$-matrices acting on quantum $\mathfrak{sl}_2$ representations~$V(t)^{\otimes n}$~\cite{Ohtsuki}. The ADO invariants have appeared as traces of certain homological representations in \cite{AnghelHomological,Ito,MartelWilletts}. One may construct braid representations from $\Vt$ by restricting to certain weight spaces, but there does not appear to be a simple interpretation as a homological representation.
\begin{Question}
 Is there a higher rank, multivariable analog of the Burau representation which recovers $\Delta_{\mathfrak{sl}_3}$ as a determinant? What is the geometric interpretation of such a representation?
\end{Question}

	 It is also shown in \cite{BCGP} that the Reidemeister torsion is recovered from TQFTs based on the $\mathfrak{sl}_2$ representations $V(t)$. We expect that applying their TQFT to higher rank quantum groups at a~fourth root of unity generalizes Reidemeister torsion and implies a Turaev surgery formula~\cite{TuraevReidemeister} in terms of $\Delta_{\mathfrak{sl}_3}$. Such a formula is likely to appear in the relation between the CGP invariant~\cite{costantino2014} and the $\widehat{Z}$-invariant~\cite{gukov2017,gukov2020} for $\mathfrak{sl}_3$, extending the results of the invariants in rank one at certain roots of unity~\cite{costantino2023,CHRY,ferrari2024}.

	\subsection{Structure of paper} In Section~\ref{sec:qgrp}, we recall the restricted quantum group $\overline{U}_\zeta(\mathfrak{sl}_3)$ and show directly that it is a quotient of the Kac--De Concini quantum group by a Hopf ideal. We study its representations $\Vt$ and its composition series for certain nondegenerate parameters in Section~\ref{sec:reps}. We make use of the tensor product decompositions of Theorem~\ref{thm:directsum} to characterize the $R$-matrix action on these submodules and quotients~$W_\alpha(\t)$ of $\Vt$.
	
	We recall the unrolled restricted quantum group in Section~\ref{sec:unrolled}, which admits a braiding on its category of weight representations. The pivotal structure and $R$-matrix are normalized so that they do not depend on the $H_i$-weights $\mvec{\lambda}$ of the unrolled restricted quantum group representations. Thus, the $R$-matrix acts on $\Vt\otimes\Vt$ and we express it in terms of the direct sum basis from \cite{HarperVerma}.	We give an overview on computing invariants and the modified trace in Section~\ref{sec:invariants}. Here we discuss ambidexterity of $\Vt$ and well-definedness of the unframed link invariant, then prove that the four-dimensional representations $W_\alpha(\t)$ yield the Alexander polynomial in the variable $t^4$ for any link $\mathcal{L}$.
	
	Section~\ref{sec:properties} is concerned with the same properties of $\Delta_{\mathfrak{sl}_3}$ from Theorem~\ref{Thm:list}. We prove Theorem~\ref{thm:plugin}, describe the $\Delta_{\mathfrak{sl}_3}$ skein relation, and a method to compute $\Delta_{\mathfrak{sl}_3}$ for families of torus knots. The invariant $\Delta_{\mathfrak{sl}_3}$ is tabulated on prime knots up to seven crossings along with several other examples in Section~\ref{sec:values}. We also make several observations regarding these polynomials and their presentation.

 Proofs of Proposition~\ref{prop:quasiR} and Lemma~\ref{lem:intertwineraction}, which involve longer computations, are given in Appendices \ref{sec:quasiR} and \ref{sec:intertwineraction}.

\section{Restricted quantum \texorpdfstring{$\boldsymbol{\mathfrak{sl}_3}$}{sl\_3}}\label{sec:qgrp}
	We recall the restricted quantum group $\overline{U}_\zeta(\mathfrak{sl}_3)$, which is a quotient of the Kac--De Concini--Procesi ``unrestricted specialization.'' For convenience of the reader, we only present $\mathfrak{g}=\mathfrak{sl}_3$ at a fourth root of unity, and the main results of this paper will only be stated for this case. In future work we consider other roots of unity and Lie types.

 \begin{Notation}
 Throughout this paper, $\zeta$ is a fixed primitive fourth root of unity.
 \end{Notation}

	Let $\proots$ ($\psroots$) be a set of positive (simple) roots for the $\mathsf{A}_2$ root system. Let $A=\left(\begin{smallmatrix}
	 2&-1\\-1&2
	\end{smallmatrix}\right)$ denote the Cartan matrix with corresponding bilinear pairing $\langle\alpha_i,\alpha_j \rangle=A_{ij}$. Fix the word presentation $w=w_1w_2w_1$ for the longest element $s_\bullet$ in the Weyl group. This presentation determines an ordering $\prec_w$ on $\proots$
\[\alpha_1\prec_w\alpha_1+\alpha_2\prec_w\alpha_2.\]

 For $m,n\in\mathbb{N}_0$, quantum numbers, factorials, and binomials are denoted
\begin{align*}
& [n]=\dfrac{\zeta^n-\zeta^{-n}}{\zeta-\zeta^{-1}},
\qquad
 [n]!=\prod_{j=1}^n [j],
\qquad
 \mbox{and}
\qquad
 \qbin{m+n}{n}{}
 =\dfrac{[m+n]!}{[m]![n]!},
\end{align*}
and take values in $\mathbb Z[\zeta]\,$. We also use the notation
\[
\floor{x}=\dfrac{x-x^{-1}}{\zeta-\zeta^{-1}}.
\]
	
	 The following is the Kac--De~Concini quantum group for $\mathfrak{sl}_3$, also known as the \textit{unrestricted specialization} of the quantum group at a root of unity. This algebra was first studied for simple~$\mathfrak{g}$ primarily at odd roots of unity in a series of papers \cite{KDC90,KDC92,KDCP}.
	 \begin{Definition}
	 	 Let $U_\zeta(\mathfrak{sl}_3)$ be the algebra over $\mathbb{Q}(\zeta)$ generated by $E_i$, $F_i$, and $K^{\pm1}_i$ for ${1\leq i\leq 2}$ subject to the relations:
	 	\begin{align}
	 		&K_iK_i^{-1}=1, \qquad K_iK_j=K_jK_i,\notag
 \\
	 		&K_iE_j=\zeta^{A_{ij}}E_jK_i, \qquad K_iF_j=\zeta^{-A_{ij}}F_jK_i,\label{eqn:qgroup2}
\\
	& 		[E_i,F_j]=\delta_{ij}\frac{K_i-K_i^{-1}}{\zeta-\zeta^{-1}},\nonumber
\\
\label{eq:SerreE}
	 		&\sum_{r+s=1-A_{ij}}(-1)^s\qbin{1-A_{ij}}{s}{}E_i^rE_jE_i^s=0 \qquad \text{for }i\neq j,
	 		\\\label{eq:SerreF}
	 		&\sum_{r+s=1-A_{ij}}(-1)^s\qbin{1-A_{ij}}{s}{} F_i^rF_jF_i^s=0 \qquad \text{for }i\neq j.
	 	\end{align}
 We write $U$ to denote ${U}_\zeta(\mathfrak{sl}_3)$.
	 \end{Definition}
 The Hopf algebra structure on $U$ is defined on generators by
\begin{alignat}{4}
	&\Delta(E_i)=E_i\otimes K_i+1\otimes E_i, \qquad && S(E_i)=-E_iK_i^{-1},\qquad && \epsilon(E_i)=0,\notag&\\
	&\Delta(F_i)=F_i\otimes 1+ K_i^{-1}\otimes F_i, \qquad && S(F_i)=-K_iF_i,\qquad && \epsilon(F_i)=0,\notag&\\
	&\Delta(K_i)=K_i\otimes K_i,\qquad && S(K_i)=K_i^{-1},\qquad && \epsilon(K_i)=1 .&\label{eq:HopfE}
\end{alignat}

 Let $\mho\colon U\to U^{\rm op}$ be the anti-involution on $U$ determined from
 \begin{align}\label{eq:mho}
 & \mho(E_i)=E_i, \qquad \mho(F_i)=F_i,
 \qquad \mho(K_i)=K_i^{-1}.
\end{align}
In \cite{Lusztig}, Lusztig defines a set of automorphisms indexed by $1\leq i\leq n$ on quantum groups given by
\begin{align*}
& T_i(E_i)=-F_iK_i,\qquad T_i(F_i)=-K_i^{-1}E_i,
 \qquad T_i(K_j)=K_jK_i^{-A_{ij}},
\end{align*}
and for $i\neq j$,
\begin{align*}
& T_i(E_j)=\sum_{r+s=-A_{ij}}\frac {(-1)^r\zeta^{-s}} {[r]![s]!}E_i^rE_jE_i^{s},
 \qquad T_i(F_j)=\sum_{r+s=-A_{ij}}\frac {(-1)^r\zeta^{s}} {[s]![r]!}F_i^{s}F_jF_i^{r}.
\end{align*}
These actions together with our chosen presentation $w$ of $s_\bullet$ determine expressions for non-simple root vectors \begin{align*}
	& E_{12}=T_1(E_2)=-(E_1E_2+\zeta E_2E_1) \qquad \text{and}\qquad F_{12}=T_1(F_2)=-(F_2F_1-\zeta F_1F_2).
\end{align*}

	\begin{Definition}
 Define the \emph{restricted quantum group} $\overline{U}_\zeta(\mathfrak{sl}_3)=U_\zeta(\mathfrak{sl}_3)/\bigl\langle E_\alpha^{2},F_\alpha^{2}\mid \alpha\in\proots\bigr\rangle$.
 We write $\overline{U}$ to denote $\overline{U}_\zeta(\mathfrak{sl}_3)$.
	\end{Definition}

The Serre relations in \eqref{eq:SerreE} and \eqref{eq:SerreF} vanish in $\overline{U}_\zeta(\mathfrak{sl}_3)$ since $[2]_\zeta=0$. The relation $E_{12}^2=0$ in $\overline{U}$ is equivalent to $E_1E_2E_1E_2=E_2E_1E_2E_1$. The latter equality holds for any choice of presentation of $E_{12}$ then imposing $E_{12}^2=0$.

\begin{Remark}
For simple $\mathfrak{g}$ of rank greater than 1, different presentations of $s_\bullet$, the longest element of the Weyl group, determine different expressions for $E_\alpha$ and $F_\alpha$. \emph{A priori} an ideal generated by a set of some powers of $E_\alpha^2$ and $F_\alpha^2$ may determine different quotients of $\overline{U}_\zeta(\mathfrak{g})$ depending on the presentation of $s_\bullet$. In other words, it is not obvious that $\overline{U}_\zeta(\mathfrak{g})$ is well defined.

 In \cite{HK}, we consider all roots of unity $\zeta$ (not necessarily of order four) where $U_\zeta(\mathfrak{g})$ is defined. We characterize which ideals generated by subsets of $E_\alpha^{\ell_\alpha}$ are Hopf ideals in $U_\zeta(\mathfrak{g})$, here $\ell_\alpha$ is the order of the root of unity $\zeta^{2d_\alpha}$ and $d_\alpha$ is a constant related to the length of $\alpha$ in the root system of $\mathfrak{g}$. We prove that the restricted quantum group $\overline{U}_\zeta(\mathfrak{g})$, defined as the quotient by the ideal generated by all $E_\alpha^{\ell_\alpha}$ and $F_\alpha^{\ell_\alpha}$, is a well-defined Hopf algebra, independent of the presentation of $s_\bullet$, for any simple $\mathfrak{g}$.
\end{Remark}

	\begin{Proposition}
		The Hopf algebra structure on $\overline{U}_\zeta(\mathfrak{sl}_3)$ is inherited from $U_\zeta(\mathfrak{sl}_3)$.
	\end{Proposition}
	\begin{proof}
		We verify that the two-sided ideal $J$ generated by $\bigl\{E_\alpha^2 \mid \alpha\in \proots \bigr\}$ is a Hopf ideal, the proof is analogous for $\bigl\{F_\alpha^2 \mid \alpha\in \proots \bigr\}$. It is enough to show that $\Delta(J)\subseteq J\otimes U+U\otimes J$ and $S(J)\subseteq J$. These relations are readily verified on the generators $E_1^2$ and $E_2^2$ from~\eqref{eq:HopfE}. We now consider $E_{12}^2$,
		\begin{align*}
			& E_{12}^2=(E_1E_2+\zeta E_2E_1)^2=(E_1E_2)^2+\zeta E_1E_2^2E_1+\zeta E_2E_1^2E_2-(E_2E_1)^2.
		\end{align*}
		It is enough to show $\Delta(E_1E_2)^2-\Delta(E_2E_1)^2\in J\otimes U+U\otimes J$, as the other terms clearly belong to $J$. We have
		\begin{gather*}
			\Delta(E_1E_2)^2 =(E_1E_2\otimes K_1K_2+E_1\otimes K_1E_2+E_2\otimes E_1K_2+1\otimes E_1E_2)^2,\\
			\Delta(E_1E_2)^2 +J\otimes U+U\otimes J=(E_1E_2)^2\otimes (K_1K_2)^2+E_1E_2E_1\otimes E_2K_1^2K_2\\
\hphantom{\Delta(E_1E_2)^2}{}
+\zeta E_1E_2\otimes E_2E_1K_1K_2+E_1\otimes E_2E_1E_2K_1+E_2E_1E_2\otimes E_1K_1K_2^2\\
\hphantom{\Delta(E_1E_2)^2}{} +\zeta E_2E_1\otimes E_1E_2K_1K_2+
			E_2\otimes E_1E_2E_1K_2+1\otimes (E_1E_2)^2+
			J\otimes U+U\otimes J.
		\end{gather*}
		The computation for $\Delta(E_2E_1)^2$ is identical to the above except the indices are switched. Thus, $\Delta(E_2E_1)^2-\Delta(E_1E_2)^2\in J\otimes U+U\otimes J$.
		
		To verify the antipode relation, we will again show the computation for the $E_{12}^2$ case. Since $S(E_1E_2)=-\zeta E_2E_1K_1^{-1}K_2^{-1}$, it follows that $S(E_{12})^2=\bigl(-\zeta E_2E_1K_1^{-1}K_2^{-1}+ E_1E_2K_1^{-1}K_2^{-1}\bigr)^2$ and it
 is easily seen to belong to $J$.
	\end{proof}

 \begin{Notation}\label{nota:psi}
 Let $\P$ denote the multiplicative characters on the Cartan subalgebra $\bigl\langle\hspace{-0.5pt} K_1^{\pm1}, K_2^{\pm1}\hspace{-0.5pt}\bigr\rangle$. There is a natural identification $\P\cong(\mathbb{C}^\times)^{2}$ by mapping $\t\in\P$ to its values on the pair $(K_1,K_2)$. There is a group structure on $\P$ under entrywise multiplication with identity $\mvec{1}=(1,1)$.

 Let $\maps$ denote the space of maps $\{0,1\}^{\proots}$. For $\psi\in\maps$, let \begin{align*}
 E^\psi=\prod_{\alpha\in \proots} E_\alpha^{\psi(\alpha)}
\qquad
 \mbox{and}
\qquad
 F^\psi=\prod_{\alpha\in \proots} F_\alpha^{\psi(\alpha)},
\end{align*}
where the product is ordered according to $\prec_w$. Write ${\psi}^\vee=\sum_{\alpha\in \proots}\psi(\alpha)\alpha\in\mathbb{Z}^\psroots$ and $\mvec{\sigma_\psi}(\cdot)=\zeta^{-\langle \psi^\vee, \,\cdot\, \rangle }\in \P$ so that
\begin{align*}
 K_iE^\psi=\mvec{\sigma_\psi}(\alpha_i)^{-1} E^\psi K_i
\qquad
 \mbox{and}
\qquad
 K_iF^\psi=\mvec{\sigma_\psi}(\alpha_i) F^\psi K_i.
\end{align*}
For $\psi,\psi'\in\maps$, the identity $\mvec{\sigma_\psi}\mvec{\sigma_{\psi'}}=\mvec{\sigma_{\psi+\psi'}}$ holds by linearity of the pairing. If $\psi^\vee=\alpha_i$, then we may also write $\mvec{\sigma_i}=\mvec{\sigma_\psi}$.
\end{Notation}

\section{Representations of \texorpdfstring{$\boldsymbol{\overline{U}_\zeta(\mathfrak{sl}_3)}$}{U\_z(sl\_3)}}\label{sec:reps}
	Here we recall the representation $\Vt$ as a Verma module over $\overline{U}$. We then characterize the structure of $\Vt$ when it has a four-dimensional irreducible subrepresentation. The Jordan--H\"older series in these cases are implied by the exact sequences given in Propositions~\ref{prop:Xexact} and~\ref{prop:Hexact}. Theorems~\ref{thm:Vdecomp}, \ref{thm:directsum}, and~\ref{thm:mixed} state the tensor product decompositions for these representations. The category of representations of~$\overline{U}$ is studied further in~\cite{HarperVerma}.

	\subsection{Induced representations}
 Let $(V,\rho)$ be a finite-dimensional $\overline{U}$-module, where $\rho\colon \overline{U}\to\End_{\mathbb{C}}(V)$ is called the \emph{action} of~$\overline{U}$ on~$V$.
 We say that $(V,\rho)$, or simply $V$ when understood from context, is a \emph{weight module} if $K_1$ and $K_2$ act \emph{semisimply} on~$V$, i.e., there exists a basis of $V$ for which the actions of $K_1$ and $K_2$ are simultaneously diagonalizable. Such a basis is called a \emph{weight basis} and it is generated by so-called \emph{weight vectors}. Let $\Uwmod$ be the category of $\overline{U}$-weight modules and their $\overline{U}$-linear maps.

 Let $B=\bracks{E_\alpha,K_i^\pm\mid \alpha\in \proots,\, 1\leq i\leq 2}$ be the Borel subalgebra of $\overline{U}$. Each character $\t\in\P$ extends to a character $\gamma_{\t}\colon B\rightarrow \mathbb{C}$ by
\begin{align*}
	\gamma_{\t}(K_i)=t_i \qquad \text{and} \qquad \gamma_{\t}(E_i)=0.
\end{align*}

	\begin{Definition}
	 Let $V_{\t}=\bigl\langle v_0^{\t}\bigr\rangle$ be the \hbox{one-dimensional} left $B$-module determined by $\gamma_{\t}$, i.e., for each $b\in B$, $b v_0^{\t}=\gamma_{\t}(b)v_0^{\t}$.
		Define $V(\t)$ to be the induced module
		\begin{align*}
		\Vt\coloneq{\text{{Ind}} _{B}^{\overline{U}}(V_{\t})=\overline{U}\otimes_{B}V_{\t}}.
		\end{align*}
	\end{Definition}

 These representations are naturally defined for any restricted quantum group and are referred to as \emph{diagonal modules} in the Kac--De~Concini/unrestricted quantum group setting \cite{KDC90}.

From the PBW basis \cite{Lusztig}, we have that $\Vt\cong U^-$ as vector spaces and
\begin{gather*}
\{1, F_1, F_2,F_1F_2, F_{12}, F_1F_{12}, F_{12}F_2, F_1F_{12}F_2\}\nonumber\\
\qquad =
		\bigl\{ F^{(000)},F^{(100)},F^{(001)},F^{(101)},F^{(010)},F^{(110)},F^{(011)},F^{(111)}\bigr\} 
\end{gather*}
	is an ordered basis of $U^-$. This basis determines \textit{the standard PBW basis} for $\Vt$ by tensoring with $v_0^{\t}$.
	
\begin{figure}[!ht]
		\centering
		\includegraphics[scale=1.35]{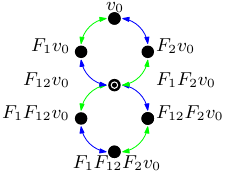}
		\caption{The action of $\overline{U}$ on the weight spaces of $\Vt$.}\label{fig:action}
	\end{figure}

We give the actions of $E_1$ and $E_2$ on the standard PBW basis in \hbox{Table \ref{table:actions}} below. We also provide a graphical description of the action of $\overline{U}$ on $\Vt$ in terms of weight spaces labeled by the standard PBW basis in Figure~\ref{fig:action}.

 \begin{Notation}
 Figure~\ref{fig:action} is read as follows. Each solid vertex indicates a one-dimensional weight space of $\Vt$, and the ``dotted'' vertex indicates the two-dimensional weight space spanned by~$F^{(101)}v_0$ and $F^{(010)}v_0$. An upward pointing edge is drawn between vertices if the action of either~$E_1$ or~$E_2$ is nonzero between the associated weight spaces. Downward edges are used to indicate nonzero matrix elements of~$F_1$ and~$F_2$. Green (blue) colored edges indicate actions of~$E_1$ and~$F_1$ ($E_2$~and~$F_2$). For atypical values of~$\t$, see Definition~\ref{def:typical} below, arrows are deleted from the graph because matrix elements of~$E_1$ and~$E_2$ vanish.
 \end{Notation}

	\begin{table}[!ht]	
\centering
\caption{Actions of $E_1$ and $E_2$ on $\Vt$ expressed in the standard PBW basis. The remaining actions are zero for all $\t\in\P$.}\label{table:actions} \vspace{1mm}

$\begin{aligned}
		&E_1F^{(100)}v_0^{\t}=\floor{t_1}F^{(000)}v_0^{\t},\qquad &&E_2F^{(001)}v_0^{\t}=\floor{t_2}F^{(000)}v_0^{\t},&\\
		&E_1F^{(101)}v_0^{\t}=\floor{\zeta t_1}F^{(001)}v_0^{\t},\qquad && E_2F^{(101)}v_0^{\t}=\floor{t_2}F^{(100)}v_0^{\t},&\\
		&E_1F^{(010)}v_0^{\t}=\zeta t_1F^{(001)}v_0^{\t},\qquad && E_2F^{(010)}v_0^{\t}=-t_2^{-1}F^{(100)}v_0^{\t},&\\
		&E_1F^{(110)}v_0^{\t}=\zeta t_1F^{(101)}v_0^{\t}-\floor{\zeta t_1}F^{(010)}v_0^{\t},\qquad
		&&E_2F^{(011)}v_0^{\t}= t_2^{-1}F^{(101)}v_0^{\t}+ \floor{t_2}F^{(010)}v_0^{\t},&\\
		&E_1F^{(111)}v_0^{\t}=\floor{t_1}F^{(011)}v_0^{\t},\qquad &&E_2F^{(111)}v_0^{\t}= \floor{t_2}F^{(110)}v_0^{\t}.&
		\end{aligned}$
\end{table}

\begin{Remark}\label{rem:dual}
 It is straightforward to verify that the dual of the weight module $\Vt$ is another weight module $V(\mvec{-t^{-1}})$.
\end{Remark}
	
	We now state the genericity condition on the highest weight of $\Vt$. Let 	
\begin{align*}
	\X_1=\bigl\{\t\in\P \mid t_1^2=1\bigr\}, \qquad
	\X_2=\bigl\{\t\in\P \mid t_2^2=1\bigr\}, \qquad
	\H=\bigl\{\t\in\P \mid (t_1t_2)^{2}=-1\bigr\},
	\end{align*}
	then set $\R$ to be the union of $\X_1$, $\X_2$, and $\H$. Expressing $E^{(111)}F^{(111)}v_0^{\t}$ in the standard PBW basis proves the following.

 	\begin{Proposition}[\cite{HarperVerma}]
		The representation $\Vt$ of $\overline{U}$ is irreducible if and only if $\t\notin\R$.
	\end{Proposition}

 \begin{Definition}\label{def:typical}
 We say that $\t\in\P$ is \emph{typical} if $\t\notin\R$ and is \emph{atypical} otherwise.
 \end{Definition}

 Partition $\R$ into disjoint subsets indexed by nonempty subsets $I\!\subseteq \!\proots$ with $\R_{I}=\bigl(\bigcap_{\alpha\in I} \X_\alpha\bigr)\setminus \bigl(\bigcup_{\alpha\notin I}\X_\alpha\bigr)$. Note that $R_{\proots}=\varnothing$. We denote by $\R_\varnothing$ the set $\P\setminus\R$ of typical $\t$. If $\t$ belongs to~$\R_1$,~$\R_2$, or~$\R_{12}$, then $\Vt$ has a unique simple submodule. The socle, i.e., the sum of irreducible submodules, of~$\Vt$ is an irreducible subrepresentation of dimension four and the head, i.e., the quotient by the intersection of all maximal submodules, is four-dimensional and has highest weight $\t$.

	\subsection{Representations \texorpdfstring{$\boldsymbol{W_i(\t)}$}{W\_i(t)}}\label{sec:Wi} We first consider the ``simple'' degeneracies $\t\in\X_i$ for $i\in\{1,2\}$. Use $B_i$ to denote the subalgebra of $\overline{U}$ generated by $B$ and $F_i$.
	\begin{Definition}
		Suppose $\t\in\X_i$. Let $\gamma_{\t}^{W_i}$ be the extension of the character $\gamma_{\t}$ on $B$ to $B_i$ with $\gamma_{\t}^{W_i}(F_i)=0$. Set $W_{i,\t}=\bigl\langle w_0^{i,\t}\bigr\rangle$ to be the one-dimensional $B_i$-module determined by $\gamma_{\t}^{W_i}$ and define
		the four-dimensional module \begin{align*}
		W_i(\t)={\text{{Ind}} _{B_i}^{\overline{U}}(W_{i,\t})=\overline{U}\otimes_{B_i}W_{i,\t}}.
		\end{align*}
	\end{Definition}
	\begin{Remark}
		The representation $W_i(\t)$ is defined if and only if $\t\in \X_i$. Indeed
 \begin{align*}
			0=[E_i,F_i]w_0^{i,\t}=\floor{K_i}w_0^{i,\t}=\floor{t_i}w_0^{i,\t}
		\end{align*} if and only if $t_i^2=1$.
	\end{Remark}

 Recall $\mvec{\sigma_\psi}$ from Convention \ref{nota:psi} which implies $\t \cdot \mvec{\sigma_\psi}$ is the weight of $F^\psi v_0\in V(\t)$.
	
	\begin{Proposition}\label{prop:Xexact}
		For $\t\in\X_i$, we have the exact sequence
		\begin{align*}
		0\rightarrow W_i(\t \cdot \mvec{\sigma_i})\rightarrow &V(\t)\rightarrow W_i(\t)\rightarrow0.
		\end{align*}
	\end{Proposition}
	
\begin{figure}[!ht]
		\centering
 \begin{subfigure}[b]{.3\linewidth}
			\centering
			\includegraphics[scale=1.35]{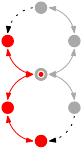}
			\caption{ $\t\in\R_1$}\label{fig:Y}
		\end{subfigure}\quad
		\begin{subfigure}[b]{.3\linewidth}
			\centering
			\includegraphics[scale=1.35]{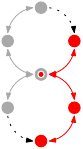}
			\caption{ $\t\in\R_2$}\label{fig:X}
		\end{subfigure}
		\caption{Reducible $\Vt$ with subrepresentation $W_i(\t \cdot \mvec{\sigma_i})$.}\label{fig:XY}
	\end{figure}

	As a subrepresentation of $\Vt$, $W_i(\t \cdot \mvec{\sigma_i})$ has a basis given by \begin{align*}
		\bigl \langle F_iv_0^{\t}, F_jF_iv_0^{\t}, F_iF_jF_iv_0^{\t}, F_jF_iF_jF_iv_0^{\t}\bigr\rangle
	\end{align*} where $\{i,j\}=\{1,2\}$.
This subrepresentation is indicated by the red points in Figure~\ref{fig:XY}. The quotient representation is colored gray and the action of $F_i$ which vanishes under the identification is indicated by a dotted arrow. Moreover, assuming $\t\in\R_i$ is equivalent to assuming both $W_i(\t\cdot \mvec{\sigma_i})$ and its quotient in $\Vt$ are irreducible.
	
	\subsection{Representations \texorpdfstring{$\boldsymbol{W_{12}(\t)}$}{W\_\{12\}(t)}}\label{sec:W12} Motivating the $\t\in \H$ case, we consider a quotient of $\Vt$ such that there is a linear dependence between the vectors $F_{12}v_0^{\t}$ and $F_{21}v_0^{\t}$, where $F_{21}=-(F_1F_2-\zeta F_2F_1)$.
 \begin{Proposition}
 Suppose $\t\in\X_{12}$. There is a quotient of $\Vt$ in which there is a linear dependence between the nonzero vectors $F_{12}v_0^{\t}$ and $F_{21}v_0^{\t}$, and this quotient is unique up to isomorphism.
 \end{Proposition}
 \begin{proof}
 We consider a quotient of $\Vt$ as a vector space by the subspace $\bigl\langle F_{12}v_0^\t- x F_{21}v_0^\t\bigr\rangle$ for some nonzero $x\in \mathbb{Q}(\zeta, t_1,t_2)$. We show that there is a unique value of~$x$ which makes this vector space into a representation. It is sufficient to consider the image of this subspace under~$E_1$ and~$E_2$. Solving for~$x$ in each equation of the system
 \begin{align*}
		E_1\bigl(F_{12}v_0^\t- x F_{21}v_0^\t\bigl)=0,\qquad E_2\bigl(F_{12}v_0^\t- x F_{21}v_0^\t\bigr)=0
	\end{align*}
 shows that $x=-\zeta t_1^2$ and $x=\zeta t_2^2$. Since $\t\in \X_{12}$, $x$ has a well-defined value and the uniqueness of the solution implies uniqueness of the quotient.
 \end{proof}

 Expanding $F_{12}v_0^{\t}-xF_{21}v_0^{\t}=0$, using the value of $x$ from the above proof, implies \[\zeta t_1^{-1}\bigl(t_1-t_1^{-1}\bigr)F_1F_2v_0^{\t}-t_2\bigl(t_2-t_2^{-1}\bigr)F_2F_1v_0^{\t}=0.\] Hence, we set \begin{align*}
 	B_{12}=\bigl \langle B, \zeta F_1F_2\floor{K_1}K_1^{-1}-F_2F_1\floor{K_2}K_2\bigr\rangle
 \end{align*} and let $\gamma_{\t}^{W_{12}}$ be the character on $B_{12}$ which is an extension of $\gamma_{\t}$ on $B$ and is zero otherwise.

	\begin{Definition}
		Let $\t\in\X_{12}$ and let $W_{12,\t}=\bigl \langle w_0^{12,\t}\bigr\rangle$ be the one-dimensional $B_{12}$-module determined by $\gamma_{\t}^{W_{12}}$. We define the four-dimensional module $W_{12}(\t)$ by induction
		\begin{align*}
		W_{12}(\t)=\text{Ind}_{B_{12}}^{\overline{U}}(W_{12,\t})=\overline{U}\otimes_{B_{12}}W_{12,\t}.
		\end{align*}
	\end{Definition}

\begin{Remark}
	To define $W_{12}(\t)$, we require $\t\in\H$ so that \[E_i\cdot\bigl(\zeta F_1F_2\floor{K_1}K_1^{-1}-F_2F_1\floor{K_2}K_2\bigr)w_0^{12,\t}=0.\] The dependence between $F_1F_2w_0^{12,\t}$ and $F_2F_1w_0^{12,\t}$ implies that $W_{12}(\t)$ is four-dimensional.
\end{Remark}

 For $\t\in\X_{12}$, there is an inclusion of $W_{12}(\t\cdot \mvec{\sigma_{(010)}})$ into $V(\t)$ which is determined by mapping~$w_0^{12,\t}$ to $\zeta t_1^{-1}\bigl(t_1-t_1^{-1}\bigr)F_1F_2v_0^{\t}-t_2\bigl(t_2-t_2^{-1}\bigr)F_2F_1v_0^{\t}$. Quotienting out this submodule returns us to the situation considered at the beginning of this subsection.

	\begin{Proposition}\label{prop:Hexact}
		If $\t\in\H$, we have the following exact sequence:
		\begin{align*}
		0\rightarrow W_{12}(\t\cdot \mvec{\sigma_{(010)}})\rightarrow V(\t)\rightarrow W_{12}(\t)\rightarrow 0.
		\end{align*}
	\end{Proposition}

\begin{figure}[!ht]
		\centering
		\begin{subfigure}[b]{.3\linewidth}\centering
			\includegraphics[scale=1.35]{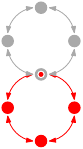}
			\caption{ $\t\in\R_{12}$}
		\end{subfigure}
		\caption{Reducible $\Vt$ with subrepresentation $W_{12}(\t\cdot \mvec{\sigma_{(010)}})$.}	\label{fig:W}
	\end{figure}

	In Figure~\ref{fig:W}, we assume $\t\in\R_{12}$ so that both $W_{12}(\t)$ and $W_{12}(\t\cdot \mvec{\sigma_{(010)}})$ are irreducible. Again, the subrepresentation is colored red and the resulting quotient is gray. Unlike Figure~\ref{fig:XY}, the trivialized actions of $F_1$ and $F_2$ are not indicated by dotted arrows because the lowest weight of $W_{12}(\t)$ is the same as the highest weight of $W_{12}(\t\cdot \mvec{\sigma_{(010)}})$, and both $F_1$ and $F_2$ act nontrivially on this weight space in the subrepresentation.
	
\subsection{Tensor product decompositions} We state three theorems on the tensor decompositions of $\Vt$ and $W_\alpha(\t)$. Each decomposition is given with a set of explicit highest weight vectors.
 \begin{Theorem}[\cite{HarperVerma}]\label{thm:Vdecomp}
	Assume that $\t$, $\s$, and $\t\s\cdot \mvec{\sigma_\psi}$ are typical for all $\psi\in \maps$. Then there is an isomorphism
	\begin{align*}
		\Vt\otimes\Vs\cong\bigoplus_{\psi\in\maps}V(\t\s\cdot \mvec{\sigma_\psi}).
	\end{align*}
\end{Theorem}
A highest weight vector for the summand $V(\t\s\cdot \mvec{\sigma_\psi})$ is $\Delta\bigl(E^{(111)}F^{(111)}\bigr)\bigl(v_0^{\t}\otimes F^\psi v_0^{\s}\bigr)$.
	
	\begin{Theorem}\label{thm:directsum}
		The following isomorphisms hold for any $\t,\s\in\P$ which satisfy the indicated additional constraints:
		\begin{align*}
& W_i(\t)\otimes W_i(\s)
\cong
		W_i(\t\s)\oplus W_i(\t\s\cdot \mvec{\sigma_{(010)}} \mvec{\sigma_{j}})\oplus V(\t\s\cdot \mvec{\sigma_j})
\qquad \mbox{for $\t,\s,\t\s\in\R_i$,}\\
&		W_{12}(\t)\otimes W_{12}(\s)
\cong W_{12}(\t\s\cdot \mvec{\sigma_i})\oplus W_{12}(\t\s\cdot \mvec{\sigma_j})\oplus V(\t\s)\qquad \mbox{for $\t,\s,\t\s\cdot \mvec{\sigma_i},\t\s\cdot \mvec{\sigma_j}\in\R_{12}$,}
		\end{align*}
 where $\{i,j\}=\{1,2\}$.
	\end{Theorem}
 The assumptions on the parameters are such that all representations are well defined, and all tensorands and summands are irreducible.
	\begin{proof}
		To establish the first isomorphism, we consider a module homomorphism \begin{align*}
			f\colon \ V(\t\s)\oplus V(\t\s\cdot \mvec{\sigma_{(010)}} \mvec{\sigma_{j}})\oplus V(\t\s\cdot \mvec{\sigma_{j}})\rightarrow W_i(\t)\otimes W_i(\s)
		\end{align*} completely determined by the image of a highest weight vector in each summand. We choose the respective images of the highest weight vectors under $f$ to be
		\begin{align*}
		w_0^{i,\t}\otimes w_0^{i,\s}, \qquad \Delta(E_jE_iE_j)\bigl(F_jF_iF_jw_0^{i,\t}\otimes F_jF_iF_jw_0^{i,\s}\bigr),\qquad \Delta(E_j)\bigl(F_jw_0^{i,\t}\otimes F_jw_0^{i,\s}\bigr).
		\end{align*}
		These three vectors have the correct weights and are clearly annihilated by $E_1$ and $E_2$. By the assumption on parameters, $\t\s\cdot \mvec{\sigma_{(010)}} \mvec{\sigma_{j}},\t\s\cdot \mvec{\sigma_{j}}\notin\R_j\cup\R_{12}$. Therefore, these vectors are nonzero and have distinct weights. By Proposition~\ref{prop:Xexact}, $V(\t\s)$ and $V(\t\s\cdot \mvec{\sigma_{(010)}} \mvec{\sigma_{j}})$ have heads $W_i(\t\s)$ and $W_i(\t\s\cdot \mvec{\sigma_{(010)}} \mvec{\sigma_{j}})$. The choice of parameters also imply that $V(\t\s\cdot \mvec{\sigma_{j}})$ is irreducible. Thus, the head of each of $V(\t\s)$, $V(\t\s\cdot \mvec{\sigma_{(010)}} \mvec{\sigma_{j}})$, and $V(\t\s\cdot \mvec{\sigma_{j}})$ is mapped to a distinct nonzero subspace under $f$. The socles of $V(\t\s)$ and $V(\t\s\cdot \mvec{\sigma_{(010)}} \mvec{\sigma_{j}})$ are irreducible and have highest weights $\t\s\cdot \mvec{\sigma_{i}}$ and $\mvec{-ts}$. Therefore, they must belong to $\ker{f}$. Quotienting this kernel yields the desired isomorphism.
		
 In the $W_{12}(\t)\otimes W_{12}(\s)$ case, the respective generating vectors are
		\begin{align*}
		\Delta(E_1)\bigl(F_1w_0^{12,\t}\otimes F_1w_0^{12,\s}\bigr),\qquad \Delta(E_2)\bigl(F_2w_0^{12,\t}\otimes F_2w_0^{12,\s}\bigr),\qquad \text{and} \qquad w_0^{12,\t}\otimes w_0^{12,\s}.\!\! \tag*{\qed}
		\end{align*}
 \renewcommand{\qed}{}
	\end{proof}
	
	Although we will not use them in this paper, we include the data of mixed tensor products for completeness. As in Theorem~\ref{thm:directsum}, parameters are assumed so that all tensorands and direct summands are well defined and irreducible.
	\begin{Theorem}\label{thm:mixed}
		The following isomorphisms hold for any $\t,\s\in\P$ which satisfy the indicated additional constraints:
		\begin{align*}
&		W_i(\t)\otimes W_j(\s)\cong V(\mvec{ts})\oplus V(\mvec{ts}\cdot \mvec{\sigma_{(010)}})\qquad \mbox{for $\t\in\R_i, \s\in\R_j, \t\s\cdot\mvec{\sigma_{(010)}}\in \R_\varnothing$,}\\
&		W_i(\t)\otimes W_{12}(\s)\cong V(\mvec{ts})\oplus V(\mvec{ts}\cdot \mvec{\sigma_{j}})
\qquad \mbox{for $\t\in\R_i, \s\in\R_{12}, \t\s\cdot\mvec{\sigma_{j}}\in \R_\varnothing$,}
\\
& V(\t)\otimes W_i(\s)\cong V(\t\s) \oplus V(\t\s\cdot \mvec{\sigma_{j}})\oplus V(\t\s\cdot \mvec{\sigma_{(010)}})\oplus V(\mvec{ts}\cdot \mvec{\sigma_{(010)}}\mvec{\sigma_j})
 \\ & \qquad
 \mbox{for $\t,\t\s\cdot\mvec{\sigma_{j}},\t\s\cdot\mvec{\sigma_{(010)}}\in\R_\varnothing, \s\in\R_{i}$,}
 \\
&		V(\t)\otimes W_{12}(\s)\cong V(\t\s)\oplus V(\t\s\cdot \mvec{\sigma_{i}})\oplus V(\t\s\cdot \mvec{\sigma_{j}})\oplus V(\t\s\cdot \mvec{\sigma_{(010)}})
 \\ & \qquad
 \mbox{for $\t,\t\s,\t\s\cdot\mvec{\sigma_{(010)}}\in\R_\varnothing, \s\in\R_{12}$,}
		\end{align*}
 where $\{i,j\}=\{1,2\}$.
	\end{Theorem}
	\begin{proof}
		Using the same argument as above, we only provide highest weight vectors which generate an irreducible representation under the action of $F_1$ and $F_2$. We then check the weights of these generating vectors, which indicate the isomorphism class of the resulting representation:
		\begin{gather*}
			W_i(\t)\otimes W_j(\s)\colon \ w_0^{i,\t}\otimes w_0^{j,\s},\ \Delta(E^{(111)})\bigl(F_jF_iF_jw_0^{i,\t}\otimes F_iF_jF_iw_0^{j,\s}\bigr),\\
W_i(\t)\otimes W_{12}(\s)\colon \ w_0^{i,\t}\otimes w_0^{12,\s},\ \Delta(E_j)\bigl(F_jw_0^{i,\t}\otimes F_jw_0^{12,\s}\bigr), \\
V(\t)\otimes W_i(\s)\colon \ v_0^{\t}\otimes w_0^{i,\s},\ \Delta(E_j)\bigl(F_jv_0^{\t}\otimes F_jw_0^{i,\s}\bigr),\ \Delta(E^{(111)})\bigl(F_iF_jF_iv_0^{\t}\otimes F_jF_iF_jw_0^{i,\s}\bigr),
\\
\hphantom{V(\t)\otimes W_i(\s)\colon}{} \ \Delta(E^{(111)})\bigl(F^{(111)}v_0^{\t}\otimes F_jF_iF_jw_0^{i,\s}\bigr),\\
V(\t)\otimes W_{12}(\s) \colon \
			v_0^{\t}\otimes w_0^{12,\s},\ \Delta(E_1)\bigl(F_1v_0^{\t}\otimes F_1w_0^{12,\s}\bigr),\ \Delta(E_2) \bigl(F_2v_0^{\t}\otimes F_2w_0^{12,\s}\bigr), \\
\hphantom{V(\t)\otimes W_{12}(\s) \colon}{} \
 \Delta(E^{(111)})\bigl(F_1F_2v_0^{\t}\otimes F^{(111)}w_0^{12,\s}\bigr).\tag*{\qed}
\end{gather*}\renewcommand{\qed}{}
	\end{proof}

 \section{Unrolled restricted quantum \texorpdfstring{$\boldsymbol{\mathfrak{sl}_3}$}{sl\_3} and braiding}\label{sec:unrolled}
\subsection{The unrolled quantum group}	We recall the unrolled restricted quantum group in Definition~\ref{defn:unrolled}. According to~\cite{GPqgps}, at odd roots of unity, the category of weight representations of an unrolled quantum group admits a~braiding $\boldsymbol{c}$. We show directly that there is a braiding for~$\mathfrak{g}=\mathfrak{sl}_3$ at a primitive fourth root of unity. We then provide a renormalization of the braiding that removes the dependence on the $H_i$-weights $\mvec{\lambda}$ up to exponentiation and thus descends to an operator on~$\Uwmod$. We end this section with its renormalized action on the tensor decomposition of~$\Vt\otimes\Vt$.
	
	\begin{Definition}\label{defn:unrolled}
 The \emph{unrolled restricted quantum group} \smash{$\overline{U}^H_\zeta(\mathfrak{sl}_3)$} is the algebra $\overline{U}_\zeta(\mathfrak{sl}_3)[H_1,H_2]$ modulo the relations for $i,j\in\{1,2\}$:
		\begin{align}\label{eq:Hcomm}
		H_iK_j^{\pm}=K_j^{\pm}H_i, \qquad H_iE_j-E_jH_i=A_{ij}E_j, \qquad H_iF_j-F_jH_i=-A_{ij}F_j.
		\end{align}
 We will use $\overline{U}^H$ as a shorthand for $\overline{U}^H_\zeta(\mathfrak{sl}_3)$.
	\end{Definition}

 A \smash{$\overline{U}^H$}-module $V$ is a \emph{weight module} if it is a direct sum of $(H_1,H_2)$-weight spaces and $H_iv=\lambda_i v$ implies $K_iv=\zeta^{\lambda_i} v$. {There exist representations where the equality $K_i=\zeta^{H_i}$ holds because the commutation relations with $E_j$ and $F_j$ in~\eqref{eqn:qgroup2} are an exponentiation of those~in~\eqref{eq:Hcomm}. } Let $\UHwmod$ denote the category of \smash{$\overline{U}^H$}-weight modules. There is a functor $\mathsf{F}_H\colon \UHwmod\to \Uwmod$ which forgets the actions of~$H_1$ and~$H_2$ and is the identity on morphisms.

\begin{Definition}
Fix a character $\t\in\P$. Choose $\mvec{\lambda}=(\lambda_1,\lambda_2)\in\mathbb{C}^2$ such that $\zeta^{\mvec{\lambda}}=\t$, by which we mean $\zeta^{\lambda_i}={\rm e}^{2\pi \sqrt{-1}\lambda_i/4}=t_i$ for each $i\in\{1,2\}$. We define $V^H(\mvec{\lambda})$ to be the unique \smash{$\overline{U}^H$}-module satisfying $\mathsf{F}_H\bigl(V^H(\mvec{\lambda})\bigr)=\Vt$ and $H_iv_0^{\mvec{\lambda}}=\lambda_iv_0^{\mvec{\lambda}}$ for each $i$, where $v_0^{\mvec{\lambda}}\in V^H(\mvec{\lambda})$ is a highest weight vector.
	\end{Definition}

 We say that $\mvec{\lambda}$ is \emph{typical} (for $V^H(\mvec{\lambda})$) if $V^H(\mvec{\lambda})$ is irreducible, or equivalently, if $\t$ is typical (for $\Vt$).

	\subsection{The \texorpdfstring{$\boldsymbol{R}$}{R}-matrix} A formula for the $R$-matrix as an operator on representations of unrolled (restricted) quantum groups at odd roots of unity is given in~\cite{GPqgps}. The formula naturally extends to even roots of unity as stated in~\cite{CR,Rupert}. We give a direct and self-contained computation that the expression in~\eqref{eqn:checkR} satisfies the quasi $R$-matrix relations for $\mathfrak{g}=\mathfrak{sl}_3$ at a fourth root of unity. Thus, defining a braiding in~\eqref{eqn:R}. We refer to~\cite{EGNO15} for the definition of a braiding.
	
For each pair of representations $(V,\rho),(W,\rho')\in\UHwmod$, we define an automorphism $\E_{V,W}$ as follows. Let $v\in V$ and $w\in W$ be weight vectors such that $H_iv=\lambda_i v$ and $H_jw=\mu_j w$, then
	\begin{align*}
	\E_{V,W}(v\otimes w)=
		\zeta^{{\sum_{ij}(A^{-1})_{ij}\lambda_i\mu_j}}(v\otimes w)=
 \zeta^{\frac{2}{3}(\lambda_1 \mu_1+\lambda_2\mu_2)+\frac{1}{3}(\lambda_1\mu_2+\lambda_2\mu_1)}(v\otimes w).
	\end{align*}
Thus, $\E_{V,W}$ can be thought of as the formal expression \smash{$\zeta^{{\sum_{ij}(A^{-1})_{ij}H_i\otimes H_j}}$}, which one may formalize in a topological completion of \smash{$\overline{U}^H$} but it will not be necessary in our treatment here. Let $\Psi_\zeta$ be the automorphism of \smash{$\overline{U}^H\otimes \overline{U}^H$} defined so that for all \smash{$x,y\in\overline{U}^H$} of weights $\alpha$ and $\beta$, respectively,
		\begin{align*}
		\Psi_{\zeta}(x\otimes y)=\zeta^{-\bracks{\alpha, \beta}}\bigl(xK_\beta^{-1}\otimes yK_\alpha^{-1}\bigr). 
		\end{align*}
	Following the computations given in \cite[Proposition~10.1.19]{CP} and \cite[Lemma~40]{GPqgps}, $\E_{V,W}$ implements $\Psi_\zeta$ on tensor products of weight representations in the sense that for all \smash{$x,y\in\overline{U}^H$} the following relation holds:
	\begin{align*}
		(\rho\otimes\rho')(\Psi_\zeta(x\otimes y))=\E_{V,W}^{-1}\circ (\rho(x)\otimes\rho'(y))\circ \E_{V,W}.
	\end{align*}

\begin{Definition}
	An invertible element \smash{$R\in \overline{U}^H\otimes \overline{U}^H$} is called a \emph{quasi-$R$-matrix} if it satisfies the following relations:
	\begin{align*}
		(\Psi_\zeta)_{23}({R}_{13}){R}_{23}=(\Delta\otimes 1)({R}),\qquad	
		(\Psi_\zeta)_{12}({R}_{13}){R}_{12}=(1\otimes \Delta)({R}),
	\end{align*}
 and ${R}\Delta(x)=\Psi_\zeta(\Delta^{\rm op}(x)){R}$ for all $x\in \overline{U}^H$.
\end{Definition}

For each $\alpha\in\proots$, define the \emph{elementary} quasi-$R$-matrix
\begin{align*}
 {\Rt_\alpha}=1\otimes 1+ \bigl(\zeta-\zeta^{-1}\bigr)E_{\alpha}\otimes F_{\alpha}\in \overline{U}^H\otimes \overline{U}^H
\end{align*}
with inverse $({\Rt_\alpha})^{-1}=1\otimes 1- \bigl(\zeta-\zeta^{-1}\bigr)E_{\alpha}\otimes F_{\alpha}.$
Set
	\begin{align}\label{eqn:checkR}
		{\Rt}=\prod_{\alpha\in\proots}{\Rt_\alpha},
	\end{align}
 with the ordered product multiplying on the right for larger $\alpha$ with respect to $\prec_w$. Indeed ${\Rt}$ is invertible. We prove the following in Appendix \ref{sec:quasiR}.

\begin{Proposition}\label{prop:quasiR}
	The element ${\Rt}$ is a quasi-$R$-matrix.
\end{Proposition}

For $(V,\rho),(W,\rho')\in\UHwmod$, define \begin{align}\label{eqn:R}
	 	\braid_{V,W}=P_{V,W}\circ \E_{V,W}\circ (\rho\otimes \rho')({\Rt})\in\Hom_{\overline{U}^H}(V\otimes W, W\otimes V),
	 \end{align}
 where $P_{V,W}\colon V\otimes W\to W\otimes V$ is the tensor swap $v\otimes w\mapsto w\otimes v$ for all $v\in V$ and $w\in W$.

\begin{Proposition}\label{prop:braid}
The morphism $\braid_{V,W}$ is a braiding on $\UHwmod$.
\end{Proposition}
\begin{proof}
	 Fix representations $(V,\rho),(W,\rho'),(U,\rho'')\in\UHwmod$. Since $P_{V,W}$, $\E_{V,W}$, and ${\Rt}$ are invertible, $\braid_{V,W}$ is an isomorphism. Routine computations prove that $\boldsymbol{c}_{V,W}$ is an intertwiner and satisfies the hexagon (triangle) identities: \begin{align*}
	 	\bigl(\braid_{ V, U}\otimes \mathrm{id}_W\bigr)\circ \bigl(\mathrm{id}_V\otimes \braid_{W, U}\bigr)=\braid_{V\otimes W, U},\qquad
	 	\bigl(\mathrm{id}_W\otimes \braid_{ V, U}\bigr)\circ\bigl(\braid_{V, W}\otimes \mathrm{id}_U\bigr)=\braid_{V,W\otimes U}.\!\!\!\tag*{\qed}
	 \end{align*}
 \renewcommand{\qed}{}
\end{proof}

\subsection{Duality morphisms}\label{ss:RT} A {pivot} on $\overline{U}^H$ is implemented by $K_{2\rho}^{1-r}=K_1^{-2}K_2^{-2}$, as in~\cite{GPqgps} for $r=2$ and where $2\rho$ is the sum of positive roots. We take the natural isomorphism $\varphi_V\colon V^{**}\to V$ to be the pivotal structure on the category of weight representations, which canonically identifies ${\rm eval}_v\in V^{**}$ with $v\in V$ and multiplies by $h_V=K_1^{-2}K_2^{-2}$. Given any basis $(e_i)$ of $V$ and corresponding dual basis~$(e_i^*)$, the left and right duality structures on $V$ are defined as
	\begin{alignat*}{3}
	& \lev_V(e_i^*\otimes e_j)=e_i^*(e_j),\qquad && \rev_V(e_i\otimes e_j^*)= e_j^*(h_V\cdot e_i),&\\
	& \lcoev_V(1)=\sum_i e_i\otimes e_i^*, \qquad && \rcoev_V(1)= \sum_i e_i^*\otimes \bigl(h_V^{-1}\cdot e_i\bigr),&
	\end{alignat*}
	and do not depend on the choice of basis.

 \begin{Definition}
		Fix an intertwiner $f\in\End_{\overline{U}^H}(V^{\otimes n})$. The \emph{right} or \emph{$n$-th partial quantum trace of $f$} is the intertwiner on $V^{\otimes n-1}$ given by
		\begin{align*}
 \tr_R(f)
 =
 \bigl(\mathrm{id}_{V^{\otimes n-1}}\otimes \rev_V\bigr)\circ (f\otimes \mathrm{id}_{V*})\circ \bigl(\mathrm{id}_{V^{\otimes n-1}}\otimes \lcoev_V\bigr).
		\end{align*}
		The \emph{left} or \emph{first partial quantum trace of $f$} is defined similarly,
		\begin{align*}
		\tr_L(f)=\bigl(\lev_V\otimes \mathrm{id}_{V^{\otimes n-1}}\bigr)\circ (\mathrm{id}_{V*}\otimes f)\circ \bigl(\rcoev_V\otimes \mathrm{id}_{V^{\otimes n-1}}\bigr) .
		\end{align*}
	\end{Definition}

 For a vector space $V$, let $\trace\colon \End_{\mathbb{C}}(V)\to\mathbb{C}$ denote the ordinary (vector space) trace on endomorphisms of~$V$. When $V\in \UHwmod$ and we restrict to endomorphisms belonging to \smash{$\End_{\overline{U}^H}(V)$}, recall that $\trace$ does not commute with the action of \smash{$\overline{U}^H$} and similarly in the category of $\overline{U}$-weight modules. In other words, $\trace$~is not an intertwiner which is in contrast to~$\tr_R$ and~$\tr_L$. The notation $\trace_i$ indicates the partial trace over the $i$-th tensor factor of an endomorphism of $V^{\otimes n}$.

 Let \smash{$f\in\End_{\overline{U}^H}(V^{\otimes n})$}, now observe that $\tr_R$ and $\tr_L$ may be expressed in terms of the partial trace:
 \begin{align}\label{eq:trR}
 \tr_R(f)
 =
 \trace_n((\mathrm{id}_{V^{\otimes n-1}}\otimes h_V)\circ f), \qquad \tr_L(f)=\trace_1\bigl(\bigl(h_V^{-1}\otimes \mathrm{id}_{V^{\otimes n-1}}\bigr)\circ f\bigr).
 \end{align}
 If, in addition, $V\in \UHwmod$ is irreducible, then by Schur's lemma $(\tr_R)^{n-1}(f)=a\cdot \mathrm{id}_V$ for some $a\in \mathbb{C}$ and we have
 \begin{align}\label{eq:trRprops}
 \trace\bigl(\tr_R^{n-1}(f)\bigr)=a\trace(\mathrm{id}_V)=a\dim(V).
 \end{align}
 The presentation of $\tr_R$ in \eqref{eq:trR} and the dimension factor appearing in \eqref{eq:trRprops} will be useful for the linear algebraic presentation of link invariants given in Proposition~\ref{prop:invntpres}.

 It is important that we do not consider the full quantum trace of intertwiners on irreducible $\VH$, which simply evaluate to zero. Observe, once again by Schur's lemma,
 \begin{align}\label{eq:tr0} \tr_R^n(f)=\tr_R(a\cdot\mathrm{id}_{\VH}) = a \trace(h_{\VH})=0,
 \end{align}
 where the last equality is computed from the following expression for $h_{\VH}=K_1^{-2}K_2^{-2}$ given in the standard PBW basis:
 \[
 h_{\VH}=t_1^{-2}t_2^{-2}\cdot\mathrm{diag}(1,-1,-1,1,1,-1,-1,1).
 \]

\subsection{Ribbon normalization} Define the family of maps $\theta_V^H=\tr_R\bigl(\braid_{V,V}\bigr)$ for $V\in \UHwmod$.

 \begin{Lemma}\label{lem:rib}
	If $\mvec{\lambda}$ is typical, then $\theta_{\VH}=\theta_{\mvec{\lambda}}\mathrm{id}_{\VH}$ where
 \[
 \theta_{\mvec{\lambda}}=\zeta^{-2(\lambda_1+\lambda_2)+{\sum_{ij}(A^{-1})_{ij}\lambda_i\lambda_j}}.
 \]
 Moreover, $\theta_V^H$ determines a ribbon structure on $\UHwmod$.
 \end{Lemma}
\begin{proof}
	Write $V$ for $\VH$ with basis $\{v_k\}$ and highest weight vector $v_0^{\mvec{\lambda}}$. We compute the action of $\theta_V^H$ on $v_0^{\mvec{\lambda}}$. Observe that for every $k$, $\braid_{V,V}\bigl(v_0^{\mvec{\lambda}}\otimes v_k\bigr)=\E\bigl(v_k\otimes v_0^{\mvec{\lambda}}\bigr)$ since ${\Rt}$ acts as the identity on $v_0^{\mvec{\lambda}}\otimes v_k$. Then as $\rev_V\bigl(v_0^{\mvec{\lambda}}\otimes v_k^*\bigr)=\delta_{0k}\zeta^{-2(\lambda_1+\lambda_2)}$, we have
 \[
\theta_V(v_0^{\mvec{\lambda}})=\rev_V(\E(v_0^{\mvec{\lambda}}\otimes v_0^{\mvec{\lambda}} ))=\zeta^{-2(\lambda_1+\lambda_2)+{\sum_{ij}(A^{-1})_{ij}\lambda_i\lambda_j}}v_0^{\mvec{\lambda}}=
\zeta^{-2(\lambda_1+\lambda_2)+{\sum_{ij}(A^{-1})_{ij}\lambda_i\lambda_j}}v_0^{\mvec{\lambda}}.
 \]

It remains to prove that $\theta_V^H$ is a ribbon structure, which will follow from \cite[Theorem~9]{GPtrace}. Since $\UHwmod$ is generically semisimple in the sense of loc.\ cit., it is sufficient to prove that $(\theta_{V^H(\mvec{\lambda})})^*=\theta_{V^H(\mvec{\lambda})^*}$, or equivalently $\theta_{\mvec{\lambda}}=\theta_{\mvec{-(\lambda-2)}}$. Indeed,
\[
\theta_{\mvec{-(\lambda-2)}}
=
\zeta^{2(\lambda_1+\lambda_2-4)-\frac{2}{3}((\lambda_1-2)^2+(\lambda_2-2)^2+(\lambda_1-2)(\lambda_2-2))}
=
\zeta^{-2(\lambda_1+\lambda_2)+{\sum_{ij}(A^{-1})_{ij}\lambda_i\lambda_j}}
\theta_{\mvec{\lambda}}. \tag*{\qed}
\]
\renewcommand{\qed}{}
\end{proof}

\begin{Definition}\label{defn:nbraid}
 For any two $\overline{U}^H$-weight modules $(V,\rho)$ and $(W,\rho')$, define a natural transformation
\begin{align}\label{eq:R}
 	\sbraid_{V, W}=\bigl(\theta_{W}^{-1}\otimes \mathrm{id}_V\bigr)\circ\braid_{V, W},
 \end{align}
 which we call the \emph{normalized braiding}.
\end{Definition}

It is readily verified that $\sbraid_{V, W}$ satisfies the Yang--Baxter equation by naturality of $\braid_{V,W}$ and that $\braid_{V,W}$ satisfies the Yang--Baxter equation itself by Proposition~\ref{prop:braid}. If $V=W=\VH$, we denote $\sbraid_{V,W}$ and $\braid_{V,W}$ by \smash{$\sbraid_{(\mvec{\lambda},\mvec{\lambda})}$} and \smash{$\braid_{(\mvec{\lambda},\mvec{\lambda})}$}, respectively. Although $\braid_{V,W}$ is a formal braiding in~$\UHwmod$, \smash{$\sbraid_{V,W}$}~is not since one of the hexagon identities is not valid.

 \begin{Remark}
For typical $\mvec{\lambda}$, Lemma~\ref{lem:rib} implies that the normalized braiding has unit partial trace
\begin{align*}
	\tr_R\bigl(\sbraid_{(\mvec{\lambda},\mvec{\lambda})}\bigr)=\tr_R\bigl(\bigl(\theta_{\mvec{\lambda}}^{-1}\otimes \mathrm{id}_{\VH}\bigr)\circ\braid_{(\mvec{\lambda},\mvec{\lambda})}\bigr)=\theta_{\mvec{\lambda}} ^{-1}\tr_R\bigl(\braid_{(\mvec{\lambda},\mvec{\lambda})}\bigr)=\mathrm{id}_{\VH}.
\end{align*}
\end{Remark}

\begin{Proposition}\label{prop:lambda}
		Suppose $\mvec{\lambda}, \mvec{\lambda'}\in\mathbb{C}^2$ satisfy $\zeta^{\mvec{\lambda}}=\zeta^{\mvec{\lambda'}}=\mvec{t}$. Then \smash{$\sbraid_{(\mvec{\lambda},\mvec{\lambda})}$} and \smash{$\sbraid_{(\mvec{\lambda'},\mvec{\lambda'})}$} define the same operator \smash{$\sbraid_{\mvec{\lambda}}\in\End_{\overline{U}^H}\bigl(V^H(\mvec{\lambda})\otimes V^H(\mvec{\lambda})\bigr)$}. Therefore, the operator
 \[
 \redbraid_{\t}\coloneq \mathsf{F}_H\bigl(\sbraid_{\mvec{\lambda}}\bigr)\in\End_{\overline{U}}(\Vt\otimes\Vt)
 \]
 is well defined.
	\end{Proposition}
\begin{proof}
	Write $V=\VH$. We compute the action of \smash{$\sbraid_{(\mvec{\lambda},\mvec{\lambda})}$} directly.	We may assume that~$\mvec{\lambda}$ is typical so that Theorem~\ref{thm:Vdecomp} extends to $V\otimes V$. Therefore, \smash{$\sbraid_{(\mvec{\lambda},\mvec{\lambda})}$} acts by a constant on each multiplicity-one summand and as an amplified $2\times 2$ matrix on the set of multiplicity-two summands. To compute these values, we consider the action of \smash{$\sbraid_{(\mvec{\lambda},\mvec{\lambda})}$} on the highest weight vector of each summand. Since \smash{$\sbraid_{(\mvec{\lambda},\mvec{\lambda})}$} is an intertwiner,
	\begin{align*}
		& \sbraid_{(\mvec{\lambda},\mvec{\lambda})} \bigl(\Delta\bigl(E^{(111)}F^{(111)}\bigr)\bigl(v_0\otimes F^{\psi}v_0^{\mvec{\lambda}}\bigr)\bigr)=
		\Delta\bigl(E^{(111)}F^{(111)}\bigr)\cdot\sbraid_{(\mvec{\lambda},\mvec{\lambda})}\bigl(v_0^{\mvec{\lambda}}\otimes F^{\psi}v_0^{\mvec{\lambda}}\bigr)\\
		& \qquad=\Delta\bigl(E^{(111)}F^{(111)}\bigr)\cdot\bigl(\theta_{\mvec\lambda}^{-1}\otimes \mathrm{id}_{\VH}\bigr)\circ P_{V,V}\circ \zeta^{({\sum_{ij}(A^{-1})_{ij}H_i\otimes H_j})}\bigl(v_0^{\mvec{\lambda}}\otimes F^{\psi}v_0^{\mvec{\lambda}}\bigr)
		\\ & \qquad =
		 \Delta\bigl(E^{(111)}F^{(111)}\bigr)\cdot P_{V,V}\circ\bigl(\mathrm{id}_{\VH}\otimes \theta_{\mvec\lambda}^{-1}\bigr)\circ\zeta^{({\sum_{ij}(A^{-1})_{ij}\lambda_i\otimes H_j})}\bigl(v_0^{\mvec{\lambda}}\otimes F^{\psi}v_0^{\mvec{\lambda}}\bigr).
	\end{align*}
	For each $\psi\in\maps$, we compute the action of $\theta_{\mvec{\lambda}}^{-1}\zeta^{({\sum_{ij}(A^{-1})_{ij}\lambda_i H_j})}$ on $F^{\psi}v_0^{\mvec{\lambda}}$:
	\begin{align*}
		\theta_{\mvec{\lambda}}^{-1}\zeta^{({\sum_{ij}(A^{-1})_{ij}\lambda_i H_j})} F^{\psi}v_0^{\mvec{\lambda}}=(t_1t_2)^2\zeta^{({\sum_{ij}(A^{-1})_{ij}\lambda_i (H_j-\lambda_j)})} F^{\psi}v_0^{\mvec{\lambda}}.
	\end{align*}
	Observe that $\sum_{ij}\bigl(A^{-1}\bigr)_{ij}\lambda_i (H_j-\lambda_j)F_k v_0^{\mvec{\lambda}}=-\sum_{ij}(A^{-1})_{ij}A_{jk}\lambda_i F_k v_0^{\mvec{\lambda}}=-\lambda_kF_kv_0^{\mvec{\lambda}}$. Therefore,
	\begin{align*}
		\theta_{\mvec{\lambda}}^{-1}\zeta^{({\sum_{ij}(A^{-1})_{ij}\lambda_i H_j})}F^{\psi}v_0^{\mvec{\lambda}}=(t_1t_2)^2\zeta^{-{\sum_\alpha \psi(\alpha)\lambda_\alpha }}F^{\psi}v_0^{\mvec{\lambda}}=\Biggl(\prod_{\alpha\in\proots} t_\alpha^{1-\psi(\alpha)}\Biggr)F^{\psi}v_0^{\mvec{\lambda}}.
	\end{align*}
	It remains to compute
\smash{$\Delta\bigl(E^{(111)}F^{(111)}\bigr)\cdot P_{V,V}\circ \bigl(v_0^{\mvec{\lambda}}\otimes F^{\psi}v_0^{\mvec{\lambda}}\bigr)=\Delta\bigl(E^{(111)}F^{(111)}\bigr)\bigl(F^{\psi}v_0^{\mvec{\lambda}}\otimes v_0^{\mvec{\lambda}}\bigr)
$} in terms of \smash{$\Delta\bigl(E^{(111)}F^{(111)}\bigr) \bigl(v_0^{\mvec{\lambda}}\otimes F^{\psi}v_0^{\mvec{\lambda}}\bigr)$}. However, these expressions will be independent of~$\mvec{\lambda}$ since they do not involve any~$H_i$. A~computation for the action of \smash{$\sbraid_{(\mvec{\lambda'},\mvec{\lambda'})}$} is identical and also given entirely in terms of~$\t$. Thus, $\redbraid_{\t}$ is well defined in $\End_{\overline{U}}(\Vt\otimes \Vt)$.
\end{proof}
	
\begin{Remark}
	A similar computation shows that \smash{$\redbraid_{(\mvec{\mu},\mvec{\lambda})}\circ \redbraid_{(\mvec{\lambda},\mvec{\mu})}\in\End_{\overline{U}^H}\bigl(\VH \otimes V^H(\mvec{\mu})\bigr)$} can be expressed in terms of~$\zeta^{\mvec{\lambda}}$ and~$\zeta^{\mvec{\mu}}$. The above arguments, now starting from typical $\mvec{\lambda}$ and~$\mvec{\mu}$, produce a well-defined operator in \smash{$\End_{\overline{U}}\bigl(V\bigl(\zeta^{\mvec{\lambda}}\bigr)\otimes V(\zeta^{\mvec{\mu}})\bigr)$}.
\end{Remark}
	
	Given a sequence $(a_j)_{j=1}^k\in\{1,2\}^k$ of length $k$, set $F_{(a_j)}=F_{i_1}\cdots F_{i_k}$. Recall the anti-involution $\mho$ on $U$ from \eqref{eq:mho}, which descends to $\overline{U}$.
	
	\begin{Proposition}
		 For every sequence $(a_j)_{j=1}^k$, $\redbraid_{\t}\bigl(
		 \Delta(E^{(111)}F^{(111)}) \bigl(v_0^{\t}\otimes F_{(a_j)}v_0^{\t}\bigr)\bigr)$ equals
		\begin{align*}
			\left(\prod_{i=1}^k-t_{a_i}^{-2}\right)(t_1t_2)^2\zeta^{-\sum_{1\leq i<j\leq k} \bracks{\alpha_{a_i},\alpha_{a_j}}}\Delta\bigl(E^{(111)}F^{(111)}\bigr)\bigl(v_0^{\t}\otimes \mho(F_{(a_{j})})v_0^{\t}\bigr).
		\end{align*}
	\end{Proposition}
\begin{proof}
	Continuing from the proof of Proposition~\ref{prop:lambda}, the action of $\redbraid_{\t}$ in the direct sum decomposition is given by
	\begin{align*}
		\redbraid_{\t}\bigl(\Delta\bigl(E^{(111)}F^{(111)}\bigr)\bigl(v_0^{\t}\otimes F^{\psi}v_0^{\t}\bigr)\bigr)=\biggl(\prod_\alpha t_\alpha^{1-\psi(\alpha)}\biggr)\Delta\bigl(E^{(111)}F^{(111)}\bigr)\cdot P_{V,V}\circ \bigl(v_0^{\t}\otimes F^{\psi}v_0^{\t}\bigr).
	\end{align*}
	This extends to products of simple root vectors $F_{(a_j)}$,
	\begin{align*}
&	\redbraid_{\t}\bigl(\Delta\bigl(E^{(111)}F^{(111)}\bigr)\bigl(v_0^{\t}\otimes F_{(a_j)}v_0^{\t}\bigr)\bigr)\\
& \qquad
=	\left(\prod_{i=1}^kt_{a_i}^{-1}\right)(t_1t_2)^2\Delta\bigl(E^{(111)}F^{(111)}\bigr)\cdot P_{V,V}\circ \bigl(v_0^{\t}\otimes F_{(a_j)}v_0^{\t}\bigr).
	\end{align*}	
	To describe this action coherently, we must express each \smash{$\Delta\bigl(E^{(111)}F^{(111)}\bigr) \bigl(F_{(a_j)}v_0^{\t}\otimes v_0^{\t}\bigr)$} in terms of some \smash{$\Delta\bigl(E^{(111)}F^{(111)}\bigr) \bigl(v_0^{\t}\otimes F_{(b_j)}v_0^{\t}\bigr)$}.

	For any $i$, we have
\[
\Delta\bigl(E^{(111)}F^{(111)}\bigr)(F_i\otimes 1)=\Delta\bigl(E^{(111)}F^{(111)}\bigr)\bigl(-K^{-1}_i\otimes F_i\bigr).
\]
 This is a consequence of the identity \smash{$\Delta\bigl(F^{(111)}\bigr)\Delta(F_{i})=0$} which follows from
 \[
 \Delta\bigl(F^{(111)}\bigr)\Delta(F_{i})=-\Delta\bigl((F_1F_2)^2F_i\bigr)
 \]
 and that $(F_1F_2)^2=(F_2F_1)^2$.
	Thus,
	\begin{align*}
	&	\Delta\bigl(E^{(111)}F^{(111)}\bigr) \bigl(F_{(a_j)}v_0^{\t}\otimes v_0^{\t}\bigr)\\
&\qquad =
		-\Delta\bigl(E^{(111)}F^{(111)}\bigr)\bigl(K^{-1}_{{i_1}}F_{{i_2}}\cdots F_{{i_k}}v_0^{\t}\otimes F_{{i_1}}v_0^{\t}\bigr)
		\\
&\qquad =
		-t_{{i_1}}^{-1}\zeta^{-\sum_{j=2}^n\bracks{\alpha_{a_1},\alpha_{a_j}}}\Delta\bigl(E^{(111)}F^{(111)}\bigr) \bigl(F_{{i_2}}\cdots F_{{i_k}}v_0^{\t}\otimes F_{{i_1}}v_0^{\t}\bigr).
	\end{align*}
	Proceeding inductively, we find $\Delta\bigl(E^{(111)}F^{(111)}\bigr) \bigl(F_{(a_j)}v_0^{\t}\otimes v_0^{\t}\bigr)$ is equal to
	\begin{align*}
		\left(\prod_{i=1}^k-t_{a_i}^{-1}\right)\zeta^{-\sum_{1\leq i<j\leq k} \bracks{\alpha_{a_i},\alpha_{a_j}}} \Delta\bigl(E^{(111)}F^{(111)}\bigr)\bigl(v_0^{\t}\otimes F_{i_k}\cdots F_{i_1}v_0^{\t}\bigr).
	\end{align*}
	Together with the previous computation, this proves the proposition.
\end{proof}

	\begin{Corollary}\label{cor:Rdecomp}
		Suppose that $(\t,\t)\in\P^2$ satisfy the assumptions of Theorem~{\rm \ref{thm:Vdecomp}}. Under the tensor product decomposition of $\Vt\otimes\Vt$ given in Theorem~{\rm \ref{thm:Vdecomp}}, we have
		\begin{equation*}
		\begin{tikzcd}
		\Vt \otimes \Vt \arrow{r}{\redbraid_{\t}} \arrow[swap]{d}{\cong} & \Vt \otimes \Vt \arrow{d}{\cong} \\%
		\bigoplus_{\psi\in\maps}V\bigl(\mvec{\sigma^\psi}\t^2\bigr) \arrow{r}{ r\otimes \mathrm{id}_{8\times 8}}& \bigoplus_{\psi\in\maps}V\bigl(\mvec{\sigma^\psi}\t^2\bigr)
		\end{tikzcd}	
		\end{equation*}
		with $r$ given by
		\begin{align*} \operatorname{diag}\bigl(t_1^2t_2^2,-t_2^2,-t_1^2\bigr)\oplus
		\left[\begin{smallmatrix}
		0 &-\zeta\\
		-\zeta &0
		\end{smallmatrix}\right]
		\oplus \operatorname{diag}\bigl(-t_1^{-2},-t_2^{-2},t_1^{-2}t_2^{-2}\bigr)
		\end{align*} in the basis determined by the highest weight vectors \smash{$\Delta\bigl(E^{(111)}F^{(11 1)}\bigr)\bigl(v_0^{\t}\otimes fv_0^{\t}\bigr)$} for
		\begin{align*}
			f\in\{1,F_1,F_2,F_1F_2,F_2F_1,F_1F_2F_1,F_2F_1F_2,F_1F_2F_1F_2\}.
		\end{align*}
	\end{Corollary}
	
	\section{Link invariants from \texorpdfstring{$\boldsymbol{\overline{U}_\zeta(\mathfrak{sl}_3)}$}{U\_z(sl\_3)}}\label{sec:invariants}

	The goal of this section is to prove Theorem~\ref{thm:smallAlexander}. We begin with our conventions for the Reshetikhin--Turaev functor \cite{RT,Turaevbook}, then show that we obtain an unframed invariant of oriented 1-tangles (or long knots) from ambidextrous weight representations of \smash{$\overline{U}^H$}. In Section~\ref{ss:Alexander}, we show that the quantum invariant associated to an irreducible representation $W_\alpha(\t)$ is the Alexander--Conway polynomial in the variable $t^4$.
	
The Reshetikhin--Turaev functor $F$ assigns linear maps to tangles. For $V,W\in\UHwmod$, we use the conventions of Figure~\ref{fig:RT} to define the functor on elementary tangles. As noted above, these assignments also restrict to $\overline{U}$-modules.

\begin{figure}[htpb!]\centering
\begin{tikzpicture}[baseline=0, xscale=.4, yscale=.4]
\draw[very thick,->] (0,0) to (0,2);
\node at (1,0) {$V$};
\node at (1,2) {$V$};
\draw[-to] (1,.5) to (1,1.5);
\node[right] at (1,1) {$\mathrm{id}_V$};
\end{tikzpicture}
\qquad
\begin{tikzpicture}[baseline=0, xscale=.4, yscale=.4]
\draw[very thick,->] (7,0) to [] (5,2);
\draw[over,very thick,->] (5,0) to [] (7,2);
\node[below] at (9.5,.5) {$\phantom{W}V\otimes W\phantom{V}$};
\node[above] at (9.5,1.5) {$\phantom{V}W\otimes V\phantom{W}$};
\draw[-to] (9.5,.5) to (9.5,1.5);
\node[right] at (9.5,1) {$R_{VW}$};
\end{tikzpicture}
\qquad
\begin{tikzpicture}[baseline=0, xscale=.4, yscale=.4]
\draw[very thick,->] (4,0) to [out=90, in=180] (5,1.5);
\draw[very thick,] (5,1.5) to [out=0, in=90] (6,0);
\node[below] at (8,.5) {$\phantom{^*}V\otimes V^*$};
\node[above] at (8,1.5) {$\mathbb{C}$};
\draw[-to] (8,.5) to (8,1.5);
\node[right] at (8,1) {$\rev_V$};
\end{tikzpicture}
\qquad
\begin{tikzpicture}[baseline=0, xscale=.4, yscale=.4]
\draw[very thick,->] (10,0) to [out=90, in=0] (9,1.5);
\draw[very thick,] (9,1.5) to [out=180, in=90] (8,0);
\node[below] at (12.5,.5) {$V^*\otimes V\phantom{^*}$};
\node[above] at (12.5,1.5) {$\mathbb{C}$};
\draw[-to] (12.5,.5) to (12.5,1.5);
\node[right] at (12.5,1) {$\lev_V$};
\end{tikzpicture}
\\[2em]
\begin{tikzpicture}[baseline=0, xscale=.4, yscale=.4]
\draw[very thick,->] (0,2) to (0,0);
\node at (1.25,0) {$\phantom{^*}V^*$};
\node at (1.25,2) {$\phantom{^*}V^*$};
\draw[-to] (1.25,.5) to (1.25,1.5);
\node[right] at (1.25,1) {$\mathrm{id}_{V^*}$};
\end{tikzpicture}
\qquad
\begin{tikzpicture}[baseline=0, xscale=.4, yscale=.4]
\draw[very thick,->] (5,0) to (7,2);
\draw[over,very thick,->] (7,0) to (5,2);
\node[below] at (9.5,.5) {$\phantom{V}W\otimes V\phantom{W}$};
\node[above] at (9.5,1.5) {$\phantom{W}V\otimes W\phantom{V}$};
\draw[-to] (9.5,.5) to (9.5,1.5);
\node[right] at (9.5,1) {$R_{VW}^{-1}$};
\end{tikzpicture}
\qquad
\begin{tikzpicture}[baseline=0, xscale=.4, yscale=.4]
\draw[very thick,->] (6,2) to [out=-90, in=0] (5,.5);
\draw[very thick,] (5,.5) to [out=180, in=-90] (4,2);
\node[above] at (8,1.5) {$\phantom{^*}V\otimes V^*$};
\node[below] at (8,.5) {$\mathbb{C}$};
\draw[-to] (8,.5) to (8,1.5);
\node[right] at (8,1) {$\lcoev_V$};
\end{tikzpicture}
\qquad
\begin{tikzpicture}[baseline=0, xscale=.4, yscale=.4]
\draw[very thick,->] (4,2) to [out=-90, in=180] (5,.5);
\draw[very thick,] (5,.5) to [out=0, in=-90] (6,2);
\node[above] at (8,1.5) {$V^*\otimes V\phantom{^*}$};
\node[below] at (8,.5) {$\mathbb{C}$};
\draw[-to] (8,.5) to (8,1.5);
\node[right] at (8,1) {$\rcoev_V$};
\end{tikzpicture}
\caption{
 A graphical definition of the Reshetikhin--Turaev functor on oriented
 elementary tangle diagrams.}
\label{fig:RT}
\end{figure}

\subsection{Ambidextrous representations}\label{subsec:ambi} In this subsection, we recall the notion of an ambidextrous representation as described in \cite{GPT}. These representations produce well-defined nonzero quantum invariants of links via the modified trace construction.

Let \smash{$\End^{\mathrm{RT}}_{\overline{U}^H}(V^{\otimes n})$} be the set of all endomorphisms of $V\otimes V$ generated by the image of all $(n,n)$-tangles with upward oriented boundary points under the Reshetikhin--Turaev functor. Let \smash{$f=F(\mathcal{T})\in\End^{\mathrm{RT}}_{\overline{U}^H}\bigl(V^{\otimes 2}\bigr)$} be the image of a $(2,2)$-tangle $\mathcal{T}$. To obtain a meaningful quantum link invariant associated to the closure of $\mathcal{T}$, we will not consider $\trace\bigl(h_V^{\otimes 2} f\bigr)$ which evaluates to zero in the case $V=\VH$, shown in~\eqref{eq:tr0}, but rather the $(1,1)$-tangle invariants $\tr_L(f)$ and $\tr_R(f)$. If these two partial traces are equal, then we declare the common morphism to be an invariant of the closed link. We say that $V$ is \emph{ambidextrous} if and only if $\tr_L(f)= \tr_R(f)$ for any \smash{$f\in\End_{\overline{U}^H}\bigl(V^{\otimes 2}\bigr)$}. Also define $V$ to be \emph{weakly ambidextrous} if it is ambidextrous relative to every morphism in \smash{$\End^{\mathrm{RT}}_{\overline{U}^H}\bigl(V^{\otimes 2}\bigr)$}.

Since $\UHwmod$ is a ribbon Ab-category in the sense of~\cite{GPT}, following Lemma~\ref{lem:rib}, there is a well-defined invariant of closed ribbon graphs colored by ambidextrous representations. Consider an oriented framed link $\mathcal{L}$ whose components are colored by an ambidextrous irreducible representation~$V$. Cutting a component of $\mathcal{L}$ colored by $V$ yields a $(1,1)$-ribbon tangle $\mathcal{L}^{\rm cut}$ identified with an endomorphism $F_V(\mathcal{L}^{\rm cut})$ of $V$ via the Reshetikhin--Turaev functor using the conventions given in Figure~\ref{fig:RT} and maps defined in Section~\ref{ss:RT}. Since $V$ is irreducible, $F_V(\mathcal{L}^{\rm cut})$ is a~scalar multiple of the identity and we write $F_V(\mathcal{L}^{\rm cut})=\bracks{\mathcal{L}^{\rm cut}}_V\mathrm{id}_V$. Ambidexterity of~$V$ implies~$\bracks{\mathcal{L}^{\rm cut}}_V$ is independent of the cut point and is therefore an invariant of $\mathcal{L}$ as a framed link. We denote this invariant by $F'_V(\mathcal{L})$.

 The following theorem is a straightforward adaptation of \cite[Section~5.7]{GPqgps} from the odd root of unity case to the fourth root of unity case. There is a minor difference in that we take~$x$ and~$y$ to be proportional to $F^{(111)}$ and $E^{(111)}$, respectively, corresponding to taking powers equal to $\operatorname{ord}\bigl(q^2\bigr)-1$.
	
	\begin{Theorem}\label{thm:ambi}
		If $V^H(\mvec{\lambda})\in\UHwmod$ is irreducible, then it is ambidextrous.
	\end{Theorem}

Let $\mathsf{F}\colon \UHwmod\to\mathrm{Vec}$ be the fiber functor, forgetting the action of $\overline{U}^H$. Later, we will not distinguish this functor from the one which is defined on $\Uwmod$ and forgets the action of $\overline{U}$.
\begin{Lemma}\label{lem:FRT}
 Let $\mathcal{T}_{\mvec{\lambda}}$ and $\mathcal{T}_{\mvec{\mu}}$ be any tangles with components colored by the representations $V^H(\mvec{\lambda})$ and $V^H(\mvec{\mu})$, respectively. The following diagram commutes
 \begin{equation*}
 \begin{tikzcd}
 \mathcal{T}_{\mvec{\lambda}} \arrow{r}{{F}} \arrow[swap]{d}{\odot_{\mvec{\lambda},\mvec{\mu}}}& F(\mathcal{T}_{\mvec{\lambda}})\arrow{r}{\mathsf{F}} & \mathsf{F}\circ F(\mathcal{T}_{\mvec{\lambda}}) \arrow{d}{\vert_{\mvec{\lambda}\mapsto\mvec{\mu}}} \\%
 \mathcal{T}_{\mvec{\mu}} \arrow{r}{{F}} & F(\mathcal{T}_{\mvec{\mu}})\arrow{r}{ \mathsf{F}}& \mathsf{F}\circ F(\mathcal{T}_{\mvec{\mu}}),
 \end{tikzcd}	
 \end{equation*}
 where $\odot_{\mvec{\lambda},\mvec{\mu}}$ is the change of color map on tangles and $\vert_{\mvec{\lambda}\mapsto\mvec{\mu}}$ is the evaluation of $\mvec{\lambda}$ at $\mvec{\mu}$ on matrices.
\end{Lemma}
\begin{proof}
	 It suffices to consider elementary tangles as these functors respect composition. In the standard PBW basis, the matrices assigned to each of the tangles colored by either $\mvec{\lambda}$ or $\mvec{\mu}$ are related by evaluation.
\end{proof}

 \begin{Lemma}\label{lem:P1P2}
 Suppose $V^H(\mvec{\mu})\in\UHwmod$ is not necessarily irreducible. The following hold:
 \begin{enumerate}[label=\textup{(\arabic*)}]\itemsep=0pt
 \item \label{lem:P2} $V^H(\mvec{\mu})$ is weakly ambidextrous.
 \item \label{lem:P1} If \smash{$f\in \End^{\mathrm{RT}}_{\overline{U}^H}\bigl(V^H(\mvec{\mu})\bigr)$}, then $f=\langle f\rangle \cdot \mathrm{id}_{V^H(\mvec{\mu})}$.
 \item \label{lem:RTmu} The map $\mathcal{L}\mapsto F'_{V^H(\mvec{\mu})}(\mathcal{L})$ is an invariant of framed links.
 \end{enumerate}
 \end{Lemma}
 \begin{proof}
 To prove item \ref{lem:P2}, let $\mathcal{T}_{\mvec{\lambda}}$ be an arbitrary $(2,2)$-tangle with upward oriented boundary points and its components colored by an irreducible representation $V^H(\mvec{\lambda})$. We consider the left and right partial closures of $\mathcal{T}_{\mvec{\lambda}}$ which realize the left and right partial quantum traces of $F(\mathcal{T}_{\mvec{\lambda}})$ under the Reshetikhin--Turaev functor. Ambidextrerity of $V^H(\mvec{\lambda})$, given by Theorem~\ref{thm:ambi}, implies that $\tr_R(F(\mathcal{T}_{\mvec{\lambda}}))=\tr_L(F(\mathcal{T}_{\mvec{\lambda}}))$.

 Now consider the same partial closures of $\mathcal{T}_{\mvec{\lambda}}$ but now with components colored by $V^H(\mvec{\mu})$ and denote the colored tangle by $\mathcal{T}_{\mvec{\mu}}$. By Lemma~\ref{lem:FRT}, there are equalities
 \[ (\mathsf{F}\circ\tr_R(F(\mathcal{T}_{\mvec{\lambda}})))\vert_{\mvec{\lambda\mapsto\mvec{\mu}}}=\mathsf{F}\circ\tr_R(F(\mathcal{T}_{\mvec{\mu}}))
 \qquad \mbox{and} \qquad
 (\mathsf{F}\circ\tr_L(F(\mathcal{T}_{\mvec{\lambda}})))\vert_{\mvec{\lambda\mapsto\mvec{\mu}}}=\mathsf{F}\circ\tr_L(F(\mathcal{T}_{\mvec{\mu}})).
 \]
 Thus, $\mathsf{F}\circ\tr_R(F(\mathcal{T}_{\mvec{\mu}})) = \mathsf{F}\circ\tr_L(F(\mathcal{T}_{\mvec{\mu}}))$ and so $
 \tr_R(F(\mathcal{T}_{\mvec{\mu}})) = \tr_L(F(\mathcal{T}_{\mvec{\mu}}))$. Since $\mathcal{T}_{\mvec{\mu}}$ is an arbitrary $(2,2)$-tangle colored by $V^H(\mvec{\mu})$ and upward oriented boundary points, we have shown that $V^H(\mvec{\mu})$ is weakly ambidextrous.

 The proof of \ref{lem:P1} follows from irreducibility of $V^H(\mvec{\lambda})$. As $\mathsf{F}\circ \tr_R(F(\mathcal{T}_{\mvec{\lambda}}))$ is a scalar multiple of the identity matrix, the matrix $\mathsf{F}\circ \tr_R(F(\mathcal{T}_{\mvec{\mu}}))$ is also scalar. Moreover, such a scalar is assigned to $F'_{V^H(\mvec{\mu})}(\mathcal{L})$ for any framed link $\mathcal{L}$. By~weak ambidexterity, this is independent of the presentation of $\mathcal{L}^{{\rm cut}}$, the $(1,1)$-tangle representative of $\mathcal{L}$ and proves item~\ref{lem:RTmu}.
 \end{proof}

 \subsection{Unframed link invariants}\label{subsec:unframed}

For each framed link presented as a link diagram with blackboard framing, we produce a framed link with zero framing numbers by performing the transformation $z$ at every crossing as defined in Figure~\ref{fig:twist}. This transformation changes the framing but not the underlying link type by applying unframed Reidemeister 1 moves. We have positioned the twists so that they agree with our definitions of $\sbraid_{\mvec{\lambda}}$, $\redbraid_{\mvec{t}}$, and their inverses under the Reshetikhin--Turaev functor. This allows us to define an invariant of unframed links by composition with the framed link invariant~$F'$.
	
	\begin{figure}[!ht]\centering
 $
 z\left(~\begin{tikzpicture}[baseline=9,scale=.8]
 \draw[very thick, ->] (1,0)--(0,1);
 \draw[very thick, ->,over] (0,0)--(1,1);
 \end{tikzpicture}
 ~
 \right)=
 \begin{tikzpicture}[baseline=9,scale=.8]
 \coordinate (ct) at (-.25,1.15);
 \coordinate (cl) at (-.5,.9);
 \coordinate (cb) at (-.25,.65);
 \coordinate (cr) at (0,.8);
 \draw[very thick, ->] (cb) to[out=0,in=-135]
 (cr) to (.35,1.25);
 \draw[very thick, over] (1,0) to[in=-45, out=135] (0,1) to[out=135, in=-15] (ct) to[out=180,in=90] (cl) to[out=270,in=180] (cb)to[out=0,in=-135]
 (cr);
 \draw[very thick, ->,over] (0,0)--(1,1);
 \end{tikzpicture}
 \qquad \qquad
 z\left(~\begin{tikzpicture}[baseline=9,scale=.8]
 \draw[very thick, ->] (0,0)--(1,1);
 \draw[very thick, ->,over] (1,0)--(0,1);
 \end{tikzpicture}
 ~\right)=
 \begin{tikzpicture}[baseline=9,scale=.8]
 \coordinate (ct) at (-.25,.35);
 \coordinate (cl) at (-.5,.1);
 \coordinate (cb) at (-.25,-.15);
 \coordinate (cr) at (0,0);
 \draw[very thick] (.25,-.25) to (.25,0) to[out=135,in=0] (ct) to[out=180,in=90] (cl) to[out=-90,in=180] (cb);
 \draw[very thick,->,over] (cb) to[out=0,in=-135]
 (cr) to[out=45,in=-135] (1,1);
 \draw[very thick, ->,over] (1,-.25) --(1,0)--(0,1);
 \end{tikzpicture}$
		\caption{Transformation $z$ defined locally on signed crossings.}\label{fig:twist}
	\end{figure}

\begin{Lemma}\label{lem:unframedinvariants}
 Let $V$ be an irreducible and weakly ambidextrous representation. For any framed link $\mathcal{L}$, $F'_V(z(\mathcal{L}))$ is an invariant of $\mathcal{L}$ as an unframed link which evaluates to~$1$ on the unknot. We denote the invariant of unframed links $F'_V\circ z$ by $\Delta_V$.
\end{Lemma}

 There is a natural extension of ${z}$ to braids and tangles, however we will not distinguish these extensions from $z$ itself. Let $cl$ indicate the full closure of a braid or tangle diagram, as a topological operation, yielding a link. The operations $cl$ and $z$ commute.

 We define $\psi_n^H(b)$ to be the action of a braid $b\in B_n$ on $V^{\otimes n}\in\UHwmod$, where each braid group generator $\sigma_i$ acts by the normalized braiding \smash{$F_V(z(\sigma_i))=\bigl(\theta_V^{-1}\otimes \mathrm{id}_V\bigr)\braid_{V,V}=\sbraid_{V,V}$}, as defined in~\eqref{eq:R}, in tensor positions $i$ and $i+1$ of $V^{\otimes n}$. A simple verification proves $\psi_n^H$ is a braid group representation. Therefore, $F_V(z(b))=\psi_n^H(b)$ and we have the following proposition.

	 \begin{Proposition}\label{prop:invntpres}
	 	Let $V$ be an irreducible and weakly ambidextrous \smash{$\overline{U}^H$}-weight module. For each unframed link $\mathcal{L}$ with braid representative $b\in B_n$, we have
 \begin{align*}
	 		\Delta_V(\mathcal{L})=\frac{1}{\dim V}\trace \bigl(\bigl(\mathrm{id}_{V}\otimes h_V^{\otimes n-1}\bigr)\circ\psi_n^H(b)\bigr).
	 	\end{align*}
	 \end{Proposition}
\begin{proof}
	Since the closure of $b$ is a presentation of $\mathcal{L}$, $cl(z(b))$ is a presentation of $z(\mathcal{L})$. Its modified trace is given by \smash{$\frac{1}{\dim V}\trace \bigl(\bigl(\mathrm{id}_{V}\otimes h_V^{\otimes n-1}\bigr)\circ \psi_n^H(b)\bigr)$} which computes $F'_V(z(\mathcal{L}))$.
\end{proof}

In the case $V=V^H(\mvec{\lambda})$, Lemma~\ref{lem:P1P2} implies that $\Delta_V$ is an unframed link invariant for any $\mvec{\lambda}$ and is also computed by the expression in Proposition~\ref{prop:invntpres}.

By Proposition~\ref{prop:lambda}, any morphism $F\circ z(\mathcal{T}_\mvec{\lambda})$ depends only on $\t=\zeta^{\mvec{\lambda}}$. This defines a~morphism associated to tangles colored by~$\Vt$. For any zero framed tangle $\mathcal{T}=z(\mathcal{T})$, we have $F_{\Vt}(\mathcal{T})=
\mathsf{F}_H\circ F_{\VH}(\mathcal{T})$, where $\mathsf{F}_H$ forgets the actions of $H_1$ and $H_2$. To simplify notation, we may also say that $\mathcal{T}_{\t}$ is colored by $\Vt$ and write $F(\mathcal{T}_{\t})$ for $F_{\Vt}(\mathcal{T})$.

We now have an analog of Lemma~\ref{lem:FRT} for zero framed tangles $\mathcal{T}$ which are colored by $\Vt$. The proof follows from the identity $F_{\Vt}(\mathcal{T})=
\mathsf{F}_H\circ F_{\VH}(\mathcal{T})$. In the lemma below, $\mathsf{F}$ denotes the fiber functor on $\Uwmod$.

\begin{Lemma}\label{lem:FRTVt}
 Let $\mathcal{T}_{\t}$ and $\mathcal{T}_{\s}$ be any tangles in the image of $z$ with components colored by the representations $\Vt$ and $\Vs$, respectively. The following diagram commutes
 \begin{equation*}
		\begin{tikzcd}
		\mathcal{T}_{\t} \arrow{r}{{F}} \arrow[swap]{d}{\odot_{\t,\s}}& F(\mathcal{T}_{\t})\arrow{r}{\mathsf{F}} & \mathsf{F}\circ F(\mathcal{T}_{\t}) \arrow{d}{\vert_{\t\mapsto\s}} \\%
		\mathcal{T}_{\s} \arrow{r}{{F}} & F(\mathcal{T}_{\s})\arrow{r}{ \mathsf{F}}& \mathsf{F}\circ F(\mathcal{T}_{\s}),
		\end{tikzcd}	
		\end{equation*}
 where $\odot_{\t,\s}$ is the change of color map on tangles and $\vert_{\t\mapsto\s}$ is the evaluation of $\t$ at $\s$ on matrices.
 \end{Lemma}

Before stating an analog of Lemma~\ref{lem:P1P2}, we introduce the following notion. Consider the set of endomorphisms
\begin{gather}
 \End_{\overline{U}}^{\mathrm{zRT}}(V^{\otimes n})=\{F_V\circ z(\mathcal{T}) \mid \text{$\mathcal{T}$ is a $(n,n)$-tangle}\nonumber\\
 \hphantom{\End_{\overline{U}}^{\mathrm{zRT}}(V^{\otimes n})=\{F_V\circ z(\mathcal{T}) \mid{}} \text{with upward oriented boundary points}\}.\label{eq:zRT}
\end{gather}
We say that a representation $V$ is \emph{$z$-weakly ambidextrous} if for every $f\in \End_{\overline{U}}^{\mathrm{zRT}}\bigl(V^{\otimes 2}\bigr)$ there is an equality $\tr_R(f)=\tr_L(f)$.

\begin{Lemma}\label{lem:P1P2Vt}
 For all $\t\in\P$, the following hold:
 \begin{enumerate}[label=\textup{(\arabic*)}]\itemsep=0pt
 \item \label{lem:P2Vt} $\Vt$ is $z$-weakly ambidextrous.
 \item \label{lem:P1Vt} If \smash{$f\in \End^{\mathrm{zRT}}_{\overline{U}}(\Vt)$}, then $f=\langle f\rangle \cdot \mathrm{id}_{\Vt}$.
 \item \label{lem:RTVt} The map $\mathcal{L}\mapsto F'_{\Vt}\circ z(\mathcal{L})$ is an invariant of unframed links.
 \end{enumerate}
 \end{Lemma}
The proof is a consequence of the case on $\overline{U}^H$-modules $\VH$ from Lemma~\ref{lem:P1P2} and the identity $F_{\Vt}\circ z(\mathcal{T})=
\mathsf{F}_H\circ F_{\VH}\circ z(\mathcal{T})$. As in Lemma~\ref{lem:unframedinvariants}, we denote the invariant in Lemma~\ref{lem:P1P2Vt}\,\ref{lem:RTVt} by $\Delta_{\Vt}(\mathcal{L})$.

\begin{Remark}\label{rem:invt} Alternatively, $\Delta_{\Vt}$ could be defined from the expression in Proposition~\ref{prop:braid} by composing $\psi_n^H$ with $\mathsf{F}_H$ or by setting $\zeta^{\mvec{\lambda}}=\t$.
\end{Remark}

\begin{Definition}\label{defn:invtVt}
	 Let $\Delta_{\mathfrak{sl}_3}$ be the invariant of unframed links given by the map $\mathcal{L}\mapsto \Delta_{\Vt}(\mathcal{L})$ considered as a Laurent polynomial in the variables $t_1,t_2$.
\end{Definition}

\begin{Remark}
	Although we consider only singly-colored links here, $F'$ is more generally defined in \cite{GPT} as an invariant of multi-colored framed links. With the appropriate normalizations, $\Delta_{\mathfrak{sl}_3}$~extends to an invariant of multi-colored links.
\end{Remark}
	\subsection{The Alexander--Conway polynomial from representations of \texorpdfstring{$\boldsymbol{\overline{U}_\zeta(\mathfrak{sl}_3)}$}{Uz(sl\_3)}}\label{ss:Alexander}
 Recall the four-dimensional $\overline{U}_\zeta(\mathfrak{sl}_3)$ representations $W_\alpha(\t)$ associated to each $\alpha\in\proots$, introduced in Sections~\ref{sec:Wi} and~\ref{sec:W12}.
 We consider the invariants of unframed links colored by each $W_\alpha(\t)$ for $\t\in\R_\alpha$ and show that they agree with the Alexander--Conway polynomial. It is important to note that although the invariant is the Alexander--Conway polynomial, the $R$-matrix itself does not satisfy the Alexander--Conway skein relation. Instead, the skein relation only holds after taking a modified trace.

	Let $\redbraid^\alpha_{\t}$ denote the action of $\redbraid_{\t}$ on $W_\alpha(\t)\otimes W_\alpha(\t)$ as a quotient of $\Vt\otimes \Vt$. Note that the matrix elements of $\redbraid^\alpha_{\t}$ are expressible in terms of $\t$. By Theorem~\ref{thm:directsum}, $W_\alpha(\t)^{\otimes 2}$ is semisimple and multiplicity free which implies $\redbraid^{\alpha}_{\t}$ is central in \smash{$\End_{\overline{U}}\bigl(W_\alpha(\t)^{\otimes 2}\bigr)$}. Therefore, by \cite[Lemma~1]{GPT} $W_\alpha(\t)$ is an ambidextrous representation. Since $W_\alpha(\t)$ is irreducible for $\t\in\R_\alpha$, Lemma~\ref{lem:unframedinvariants} implies that $F'_{W_\alpha(\t)}\circ z$ is an invariant of unframed links colored by $W_\alpha(\t)$ which evaluates to 1 on the unknot. We denote this invariant by $\Delta_{W_\alpha(\t)}$. Let
	\begin{align*}
\delta_{W_\alpha(\t)}=\redbraid_{\t}^\alpha-(\redbraid_{\t}^\alpha)^{-1}- \bigl(t^2-t^{-2}\bigr)\mathrm{id}_{W_\alpha(\t)^{\otimes 2}},
	\end{align*}
	which we identify with the Alexander--Conway skein relation given in Figure~\ref{fig:skein}.
	
	\begin{figure}[!ht]\centering
$
		 \begin{tikzpicture}[scale=.4,baseline=9]
		 \draw[very thick,->] (2,0) to (0,2);
\draw[over,very thick,->]
(0,0) to (2,2);
		\end{tikzpicture}
 -
 \begin{tikzpicture}[scale=.4,baseline=9]
		 \draw[very thick,->] (0,0) to (2,2);
\draw[over, very thick,->] (2,0) to (0,2);
		\end{tikzpicture}
 = \bigl(t^2-t^{-2}\bigr)
 \begin{tikzpicture}[scale=.4,baseline=9]
\draw[very thick,<-] (0,2) to[out=-60,in=90] (.5,1)
to[out=-90,in=60] (0,0);
\draw[very thick,->] (1.5,0) to[out=120,in=-90] (1,1) to[out=90,in=-120] (1.5,2);
\end{tikzpicture}
$
		\caption{Alexander--Conway skein relation in the variable $\bigl(t^{\frac{1}{2}}\bigr)^4$.}\label{fig:skein}
	\end{figure}

	\begin{Lemma}\label{lem:dV}
 		Let $\alpha\in\proots$ be a positive root. The action of $\delta_{W_\alpha(\t)}$ is zero on the four-dimensional direct summands of $W_\alpha(\t)^{\otimes 2}$.
	\end{Lemma}
\begin{proof}
	We first consider $W_i(\t)$ where $\alpha=\alpha_i$ with $i\in\{1,2\}$ and $\t\in \R_i$. Let $j$ be the complement of $i$ in $\{1,2\}$ and set $t=t_j$. Although $V(\t)^{\otimes 2}$ does not completely decompose as a~sum of irreducibles, Corollary \ref{cor:Rdecomp} can still be applied to compute the action of the braiding on vectors in~$W_i(\t)$.

 Consider $\redbraid_{\s}$ for some typical $\s\in\P$. Note that the highest weight vectors
 \[
 \Delta\bigl(E^{(111)}F^{(111)}\bigr)\bigl(v_0^{\s}\otimes v_0^{\s}\bigr),
 \quad \Delta\bigl(E^{(111)}F^{(111)}\bigr)\bigl(v_0^\s\otimes F_jF_iF_jv_0^{\s}\bigr)
 \]
 are proportional to the respective vectors
 \[
 v_0^{\s}\otimes v_0^{\s},
 \qquad
 \Delta(E_jE_iE_jF_jF_iF_j)\bigl(v_0^{\s}\otimes F_jF_iF_jv_0^{\s}\bigr)
 =
 \Delta(E_jE_iE_j)\bigl(F_jF_iF_jv_0^{\s}\otimes F_jF_iF_jv_0^{\s}\bigr).
 \]
 Therefore, they are also highest weight vectors. By Corollary \ref{cor:Rdecomp}, $\redbraid_{\s}$ acts on this pair of vectors by $s_i^2s_j^2$ and $-s_j^{-2}$, respectively. Lemma~\ref{lem:FRTVt} now implies that the actions of $\mathsf{F}\circ \redbraid_{\s}$ and $\mathsf{F}\circ \redbraid_{\t}$ on these vectors are related by the substitution $\s\mapsto\t$. It follows that $\redbraid_{\t}$ acts on
 \[
 v_0^{\t}\otimes v_0^{\t}
 \qquad \text{and}\qquad
 \Delta(E_jE_iE_j)\bigl(F_jF_iF_jv_0^{\t}\otimes F_jF_iF_jv_0^{\t}\bigr)
 \]
 by the respective scalars $t^2$ and $-t^{-2}$.

 Denote the quotient map $V(\t)\twoheadrightarrow W_i(\t)$ by $p$. Then $p\otimes p$ defines a surjection $\Vt^{\otimes 2}\twoheadrightarrow W_i(\t)^{\otimes 2}$. Naturality of the braiding implies the following diagram commutes
 \begin{equation*}
		\begin{tikzcd}
		\Vt \otimes \Vt \arrow{r}{p\otimes p} \arrow[swap]{d}{\redbraid_{\t}} & W_i(\t) \otimes W_i(\t) \arrow{d}{\redbraid_{\t}^\alpha} \\%
		\Vt \otimes \Vt
 \arrow{r}{ p\otimes p}& W_i(\t) \otimes W_i(\t),
		\end{tikzcd}	
	\end{equation*}
 and therefore $\redbraid_{\t}^\alpha$ acts on the vectors
 \[
 w_0^{\t}\otimes w_0^{\t}
 \quad \mbox{and}
 \quad
 \Delta(E_jE_iE_j)\bigl(F_jF_iF_jw_0^{\t}\otimes F_jF_iF_jw_0^{\t}\bigr)
 \]
 by those same scalars $t^2$ and $-t^{-2}$. Notice also that these vectors are the highest weight \mbox{vectors}~indicated in Theorem~\ref{thm:directsum} for the four-dimensional direct summands $W_i\bigl(\mvec{t^2}\bigr)$ and \linebreak ${ W_i\bigl(\mvec{t^2}\cdot \mvec{\sigma_{(010)}} \mvec{\sigma_{j}}\bigr)}$ of $W_i(\t) \otimes W_i(\t)$. Thus, $\redbraid_{\t}^\alpha-(\redbraid_{\t}^\alpha)^{-1}$ acts by $t^2-t^{-2}$ on both of these summands and, therefore, $\delta_{W_\alpha(\t)}$ acts on them by zero.

 For $W_{12}(\t)$, we set $\alpha=\alpha_1+\alpha_2$ and $\t=\bigl(\zeta t, \pm t^{-1}\bigr)\in\R_{12}$. A similar argument as above shows that $\redbraid_{\t}$ acts on the highest weight vectors
 \[
 \Delta(E_1)\bigl(F_1v_0^{\t}\otimes F_1 v_0^{\t}\bigr)
 \qquad \mbox{and}\qquad
 \Delta(E_2)\bigl(F_2v_0^{\t}\otimes F_2 v_0^{\t}\bigr)
 \]
 by the respective scalars $-t_2^2=-t^{-2}$ and $-t_1^2=t^2$. Under the quotient map $\Vt^{\otimes 2}\twoheadrightarrow W_{12}(\t)^{\otimes 2}$, these highest weight vectors correspond to the four-dimensional direct summands $W_{12}\bigl(\mvec{t}^2\cdot \mvec{\sigma_i}\bigr)$ and $W_{12}\bigl(\mvec{t}^2\cdot \mvec{\sigma_j}\bigr)$ of Theorem~\ref{thm:directsum}. Naturality of the braiding again implies $\redbraid_{\t}^\alpha-(\redbraid_{\t}^\alpha)^{-1}$ acts by $t^2-t^{-2}$ on these summands and that $\delta_{W_{\alpha}(\t)}$ acts on them by zero.
\end{proof}

\begin{Remark}\label{rem:convention}
	In the case $\alpha=\alpha_1+\alpha_2$, we specifically avoided $\t=\bigl(t, \pm \zeta t^{-1}\bigr)$ which would yield $\delta_{W_\alpha(\t)}=- 2\bigl(t^2-t^{-2}\bigr)$ on the four-dimensional summands of $W_\alpha(\t)^{\otimes 2}$. Replacing $\redbraid_\t^\alpha$ with its inverse, or~$t$ by $\zeta t$ in~$\delta_{W_\alpha(\t)}$ resolves this discrepancy. Since we recover the Alexander polynomial in the variable~$t^4$, which does not distinguish mirror images, using either convention is consistent with Theorem~\ref{thm:smallAlexander}.
\end{Remark}

\begin{Remark}
 We will not prove it here, but one can show that $\delta_{W_\alpha(\t)}$ acts by $-\bigl(t^2-t^{-2}\bigr)$ on the eight-dimensional summand of $W_\alpha(\t)^{\otimes 2}$.
\end{Remark}

 Recall that $W_\alpha(\t)^{\otimes2}$ is semisimple and decomposes as a multiplicity free direct sum by Theorem~\ref{thm:directsum}. Therefore, any $f\in\End_{\overline{U}}\bigl(W_\alpha(\t)^{\otimes2}\bigr)$ is expressible as a sum of scalars acting on each summand, which we write as
 \begin{align}\label{eq:EndDecomp}
 f = f_+ \cdot p_+ + f_- \cdot p_- + f_V \cdot p_V,
 \end{align}
 where $p_+,p_-,p_V\in \End_{\overline{U}}\bigl(W_\alpha(\t)^{\otimes2}\bigr)$ is the projection onto the corresponding summand according to the decomposition above.
\begin{Lemma}\label{lem:intertwineraction}
	Let $\t\in\R_\alpha$. For any $f\in\End_{\overline{U}}\bigl(W_\alpha(\t)^{\otimes2}\bigr)$,
	\begin{align*}
		\tr_R(f)=
		\begin{cases}
			\dfrac{f_+-f_-}{t_j^2+t_j^{-2}}\cdot \mathrm{id}_{W_\alpha(\t)}, & \mbox{where $\alpha=\alpha_i$ and $\{i,j\}=\{1,2\}$}\\[1em]
			\dfrac{f_+-f_-}{t_1^2-t_2^{-2}}\cdot \mathrm{id}_{W_\alpha(\t)}& \text{for } \alpha=\alpha_1+\alpha_2
		\end{cases}.
	\end{align*}
\end{Lemma}

The proof of Lemma~\ref{lem:intertwineraction} is given in Appendix \ref{sec:intertwineraction}.

	\begin{Theorem}\label{thm:smallAlexander}
 Suppose $W_\alpha(\t)$ is irreducible. Then the link invariant $\Delta_{W_\alpha(\t)}$ is equal to the Alexander--Conway polynomial evaluated at $t^4$.
\end{Theorem}
	\begin{proof}
		Fix $\alpha\in\proots$ and $\t\in \R_\alpha$. If $\alpha=\alpha_1+\alpha_2$, we may assume $\t=(\zeta t, \pm t^{-1})$ for some generic $t$, as explained in Remark \ref{rem:convention}. The Alexander--Conway relation is encoded by $\delta_{W_\alpha(\t)}$. In the notation of \eqref{eq:EndDecomp}, Lemma~\ref{lem:dV} shows $(\delta_{W_\alpha(\t)})_+=(\delta_{W_\alpha(\t)})_-=0$. Fix any intertwiner $f\in \End_{\overline{U}}\bigl(W_\alpha(\t)^{\otimes 2}\bigr)$. The first equality below follows from Lemma~\ref{lem:intertwineraction},\[\tr_R(\delta_{W_\alpha(\t)}f)=\frac{1}{r_\alpha}
 (
 (\delta_{W_\alpha(\t)}f)_+-
 (\delta_{W_\alpha(\t)}f)_-
 )
 =
 \frac{1}{r_\alpha}
 (
 (\delta_{W_\alpha(\t)})_+f_+-
 (\delta_{W_\alpha(\t)})_-f_-
 )=0,
 \]
 where $r_\alpha\in\mathbb{C}^\times$.
 In particular, if $f\in \End_{\overline{U}}^{\mathrm{zRT}}\bigl(W_\alpha(\t)^{\otimes 2}\bigr)$, considered as the image of a $(2,2)$-tangle as defined in~\eqref{eq:zRT}, we have shown that $f$ composed with the skein relation has partial quantum trace zero. Therefore, $\Delta_{W_{\alpha}(\mvec{t})}$ satisfies the Alexander--Conway skein relation.
	\end{proof}

	\section{Properties of \texorpdfstring{$\boldsymbol{\Delta_{\mathfrak{sl}_3}}$}{Delta sl\_3}}\label{sec:properties}
 \subsection{Evaluation to the Alexander polynomial}

	Here we prove $\Delta_{\mathfrak{sl}_3}$, as in Definition~\ref{defn:invtVt}, evaluates to the Alexander polynomial. We also discuss the skein relation for~$\Delta_{\mathfrak{sl}_3}$ and apply it to compute the invariant for $(2,2n+1)$ torus knots. Basic symmetry properties of the invariant are also given.

\begin{Lemma}\label{lem:invtseq}
		Suppose that $\t\in\R_\alpha$, for some $\alpha\in \proots$, such that there exist irreducible representations $V_1,V_2\in\Uwmod$ so that there is an exact sequence
		\begin{align*}
			0\rightarrow V_1\to \Vt\rightarrow V_2\rightarrow0
		\end{align*}
 as in either Proposition~{\rm \ref{prop:Xexact}} or Proposition~{\rm \ref{prop:Hexact}}. For every knot $\mathcal{K}$, there is an equality
 \[
 \Delta_{\Vt}(\mathcal{K})=\Delta_{V_1}(\mathcal{K})=\Delta_{V_2}(\mathcal{K}).
 \]
	\end{Lemma}
	\begin{proof}
 Note that $V_1$ and $V_2$ are irreducible representations $W_\alpha(\t)$ and $W_\alpha(\mvec{t'})$ for some ${\t,\t'\hspace{-0.2pt}\in \R_\alpha}$ and $\alpha\in \proots$. As discussed in the introduction to Section~\ref{ss:Alexander}, these representations produce well-defined link invariants. By Lemma~\ref{lem:P1P2Vt}\,\ref{lem:RTVt}, $\Delta_{\Vt}$ is well defined for all $\t\in \P$.

 Fix a knot $\mathcal{K}$ and let $\mathcal{K}^{\rm cut}$ be a $(1,1)$-tangle with closure $\mathcal{K}$. Now, $\Delta_{\Vt}(\mathcal{K})$ is computed as the scalar multiple of the identity $F_{\Vt}\circ z(\mathcal{K}^{\rm cut})= \Delta_{\Vt}(\mathcal{K})\cdot \mathrm{id}_{\Vt}$.
 Naturality of the braiding and pivotal structure discussed in Section~\ref{sec:unrolled} imply that the inclusion $i\colon V_1\hookrightarrow\Vt$ satisfies the intertwiner relation \begin{align*}
			\Delta_{\Vt}(\mathcal{K})\cdot i=F_{\Vt}(z(\mathcal{K}^{\rm cut}))\circ i=i\circ F_{V_1}(z(\mathcal{K}^{\rm cut}))=	\Delta_{V_1}(\mathcal{K})\cdot i.
		\end{align*}
	Therefore, $\Delta_{V_1}(\mathcal{K})=\Delta_{\Vt}(\mathcal{K})$. Similarly, the surjection $\Vt\twoheadrightarrow V_2$ intertwines the scalar action.
	\end{proof}
	\begin{Theorem}\label{thm:plugin}
	 Let $\mathcal{K}$ be any knot. Then
		\begin{equation*}
		\Delta_{\mathfrak{sl}_3}(\mathcal{K})(t,\pm 1)=\Delta_{\mathfrak{sl}_3}(\mathcal{K})(\pm 1,t)=\Delta_{\mathfrak{sl}_3}(\mathcal{K})\bigl(t,\pm \sqrt{-1}/t\bigr)=\Apoly(\mathcal{K})\bigl(t^4\bigr).
		\end{equation*}
		Moreover, if for all knots $\mathcal{K}$ the equality $\Delta_{\mathfrak{sl}_3}(\mathcal{K})(t,s)=\Apoly(\mathcal{K})\bigl(t^4\bigr)$ holds, then either $s^2=1$ or $s^2=-t^{-2}$.
	\end{Theorem}
	\begin{proof}
 Fix a knot $\mathcal{K}$. Consider the first evaluation $\Delta_{\mathfrak{sl}_3}(\mathcal{K})(t,\pm 1)=\Delta_{\Vt}(\mathcal{K})$ for $\t\in \R_{\alpha_2}$. In this case, there is an exact sequence of modules as given in Proposition~\ref{prop:Xexact}. Now Lemma~\ref{lem:invtseq} applies and we have the equality $\Delta_{\Vt}(\mathcal{K})=\Delta_{W_{\alpha_2}(\t)}(\mathcal{K})$, which is then equal to $\Apoly(\mathcal{K})\bigl(t^4\bigr)$ by Theorem~\ref{thm:smallAlexander}. The remaining evaluations correspond to $\t\in \R_{\alpha_1}$ and $\t\in \R_{\alpha_{12}}$, respectively, and are proven using similar arguments.

 The second claim follows from checking which evaluations of $\Delta_{\mathfrak{sl}_3}$ simultaneously yield the Alexander polynomial on the knots~$\mathsf{3_1}$ and~$\mathsf{4_1}$.
	\end{proof}

Lemma~\ref{lem:invtseq} only applies to knots. If a link were colored by reducible representations $\Vt$, only the color of the open strand could be replaced by $V_1$ or $V_2$ under the naturality transformation. All other components of the diagram remain colored by $\Vt$.

 \begin{Example}
 	We give an example of how Theorem~\ref{thm:plugin} and therefore Lemma~\ref{lem:invtseq} does not apply to links. We begin by stating the nontrivial fact that the multi-colored invariant of links is well defined by \cite{GPT}, which follows from the ambidexterity of $\VH$ and also requires the modified dimension function $\mathsf{mdim}$. For a pair of irreducible representations $V$ and $W$, the modified dimension function satisfies the property
 \[
			\mdim(V)~
			\begin{tikzpicture}[baseline=15, xscale=.75, yscale=.75]
				\draw[thick,->] (0,1) to (0,2); 
				\draw[thick,over] (-1/2,1) to [out=90, in=180] (0,3/2) to [out=0, in=90] (1/2,1); 
				\draw[thick,->] (1/2,1) to [out=-90, in=0] (0,1/2) to [out=180, in=-90] (-1/2,1);
				\draw[thick,over] (0,0) to (0,1); 
				\node[right] at (0,0) {$V$};
				\node[right] at (1/2,1) {$W$};
			\end{tikzpicture}
			=
			\mdim(W)~
			\begin{tikzpicture}[baseline=15, xscale=.75, yscale=.75]
				\draw[thick] (-1/2,1) to [out=90, in=180] (0,3/2) to [out=0, in=90] (1/2,1); 
				\draw[thick,over,->] (0,1) to (0,2); 
				\draw[thick] (0,0) to (0,1); 
				\draw[thick,over,<-] (1/2,1) to [out=-90, in=0] (0,1/2) to [out=180, in=-90] (-1/2,1); 
				\node[right] at (0,0) {$W$};
				\node[right] at (1/2,1) {$V$};
			\end{tikzpicture}.
			\] Although the modified dimension function is unique up to a global scalar, it can be chosen to be the reciprocal of the Hopf link. This renormalization by the Hopf link also allows us to consider the specializations of Theorem~\ref{thm:plugin} without the polynomials vanishing. A routine calculation following \cite[Proposition~45]{GPqgps} shows that
 \begin{align*}
\Delta_{\mathfrak{sl}_2}\bigl(\mathsf{2_1^2}\bigr)(t)=\bigl(t-t^{-1}\bigr)\qquad \mbox{and} \qquad \Delta_{\mathfrak{sl}_3}\bigl(\mathsf{2_1^2}\bigr)=\bigl(t_1-t_1^{-1}\bigr)\bigl(t_2-t_2^{-1}\bigr) \bigl(t_1t_2+t_1^{-1}t_2^{-1}\bigr).
	\end{align*}
 Now consider these invariants of the singly-colored $(4,2)$ torus link $T_{4,2}$ with these normalizations. Observe
 \begin{align*}
\frac{\Apoly({T_{4,2}})\bigl(t_1^4\bigr)}{\Apoly\bigl(\mathsf{2_1^2}\bigr)\bigl(t_1^4\bigr)}=t_1^4+t_1^{-4},\qquad \dfrac{\Delta_{\mathfrak{sl}_3}({T_{4,2}})(t_1,t_2)}{\Delta_{\mathfrak{sl}_3}\bigl(\mathsf{2_1^2}\bigr)(t_1,t_2)} =t_1^{4}t_2^{4}+t_1^{4
	}+t_2^{4}+t_2^{-4}+t_1^{-4}+t_1^{-4}t_2^{-4},
	\end{align*}
 and that these polynomials are not related by either of the evaluations in Theorem~\ref{thm:plugin}, i.e., $t_2\mapsto 1$ and $t_2\mapsto \zeta t_1^{-1}$.
\end{Example}

 \subsection{Symmetry transformation on variables}\label{sec:sym}
 The statements of Theorem~\ref{Thm:list}\,\ref{thm:sym} are a consequence of the identities noted in \cite[equation~(103)]{GK} for the knot invariant $\Lambda_{-1}$. In \cite[Theorem~1.2]{GHKST}, $\Lambda_{-1}$ was extended to a link invariant and proven to be equivalent to $\Delta_{\mathfrak{sl}_3}$. This identification implies $\Delta_{\mathfrak{sl}_3}$ is palindromic and valued in $\mathbb{Z}\bigl[t_1^{\pm2},t_2^{\pm2}\bigr]$. We include a self-contained proof of the symmetry under exchange of variables $t_1\leftrightarrow t_2$.

 \begin{Remark}
 There is an additional symmetry, first noted in \cite[equation~(103)]{GK} for $\Lambda_{-1}$, which implies that $\Delta_{\mathfrak{sl}_3}(t_1,t_2)=\Delta_{\mathfrak{sl}_3}\bigl(\frac{\zeta}{ t_1t_2} , t_2\bigr)$. The author thanks an anonymous referee for pointing out this identity.
 \end{Remark}

 \begin{Remark}
 The palindrome and symmetry properties, and the evaluations of $\Delta_{\mathfrak{sl}_3}$ to the Alexander polynomial in $t^4$ taken together can be shown to imply $\Delta_{\mathfrak{sl}_3}$ is valued in $\mathbb{Z}\bigl[t_1^{\pm2},t_2^{\pm2}\bigr]$.
 \end{Remark}

 Let $\tau$ be an automorphism of the Dynkin diagram of $\mathfrak{sl}_3$. Define $\widehat{\tau}$ to be an algebra automorphism of $\overline{U}^H$ so that $\widehat{\tau}(X_i)=X_{\tau(i)}$ for $X\in\{E,F,K,H\}$. The Hopf algebra structure for $\overline{U}^H$ is intertwined by $\widehat{\tau}$.
	\begin{Lemma}\label{lem:symm}
		 The automorphism $\widehat{\tau}$ determines an automorphism $\widetilde{\tau}$ of $\UHwmod$ as a ribbon category.
	\end{Lemma}
	\begin{proof}
		We check that ${\Rt}$ is invariant under $\widehat{\tau}\otimes \widehat{\tau}$. Recall
		\begin{align*}
			{\Rt}={}& (1\otimes 1+ \bigl(\zeta-\zeta^{-1}\bigr)E_{1}\otimes F_{1})(1\otimes 1+ \bigl(\zeta-\zeta^{-1}\bigr)E_{12}\otimes F_{12})\bigl(1\otimes 1+ \bigl(\zeta-\zeta^{-1}\bigr)E_{2}\otimes F_{2}\bigr)\!
			\\={} &
			1\otimes 1
			+2\zeta(E_{1}\otimes F_{1}+E_{2}\otimes F_{2})+(2\zeta E_{12}\otimes F_{12}-4E_{1}E_2\otimes F_{1}F_2)\\
			&-4(E_{1}E_{12}\otimes F_{1}F_{12}+E_{12}E_{2}\otimes F_{12}F_2)-8\zeta(E_1E_{12}E_2\otimes F_1F_{12}F_2)\\
			={}&
			1\otimes 1
			+2\zeta(E_{1}\otimes F_{1}+E_{2}\otimes F_{2})+2\zeta(E_1E_2\otimes F_2F_1+E_2E_1\otimes F_1F_2)\\
& -2(E_1E_2\otimes F_1F_2+E_2E_1\otimes F_2F_1)
			-4\zeta ( E_{1}E_2E_1\otimes F_{1}F_2F_1+ E_2E_{1}E_{2}\otimes F_2F_{1}F_2)\\
&+8(E_1E_2E_1E_2\otimes F_1F_{2}F_1F_2)
		\end{align*}
		and that $(E_1E_2)^2=(E_2E_1)^2$ and $(F_1F_2)^2=(F_2F_1)^2$.
		From this, we see $\widehat{\tau}\otimes \widehat{\tau}({\Rt})={\Rt}$.

Let $\widetilde{\tau}$ be the endofunctor on $\UHwmod$ defined by $\widetilde{\tau} ((V,\rho) )=(V,\rho\circ\widehat{\tau})$ on representations and is the identity on morphisms as linear maps. That is, if $\mathsf{F}\colon \UHwmod\to\text{Vect}$ is the fiber functor which forgets the action of $\overline{U}^H$, then $\mathsf{F}\circ \widetilde{\tau}=\mathsf{F}$. Since $\widehat{\tau}$ is a Hopf algebra morphism, $\widetilde{\tau}$ is canonically a~strict $\otimes$-functor and $\widetilde{\tau}(V^*)=\widetilde{\tau}(V)^*$ up to canonical isomorphism. Therefore, \smash{$\widetilde{\tau}\bigl(\lev_V\bigr)=\lev_{\widetilde{\tau}(V)}$} up to canonical isomorphism and similarly for the other duality maps.
				
We prove $\widetilde{\tau}\bigl(\braid_{V,W}\bigr)=\braid_{\widetilde{\tau}(V),\widetilde{\tau}(W)}$ for any weight representations $(V,\rho)$ and $(W,\rho')$, noting that $\rho$ and $\rho'$ are suppressed in our notation for the braiding. Since $\mathsf{F}$ is injective on morphisms, it is enough to show that \smash{$\mathsf{F}\bigl(\widetilde{\tau}(\braid_{V,W})\bigr) =\mathsf{F}\bigl(\braid_{\widetilde{\tau}(V),\widetilde{\tau}(W)}\bigr)$}, which is the same as showing \smash{$\mathsf{F}\bigl(\braid_{V,W}\bigr)=\mathsf{F}\bigl(\braid_{\widetilde{\tau}(V),\widetilde{\tau}(W)}\bigr)$}. For this proof and its corollary, we~distinguish the braiding~$\braid_{V,W}$ as an abstract morphism in $\UHwmod$ from the linear map realizing it. To be more precise, the realization given in~\eqref{eqn:R} is, in fact, $\mathsf{F}\bigl(\braid_{V,W}\bigr)$. Since $\widehat{\tau}\otimes \widehat{\tau}({\Rt})={\Rt}$, we have
		\begin{align}
 \mathsf{F}(\braid_{\widetilde{\tau}(V),\widetilde{\tau}(W)})&=P_{\widetilde{\tau}(V),\widetilde{\tau}(W)}\circ \E_{\rho\circ\widehat{\tau},\rho'\circ\widehat{\tau}}\circ(\rho\circ\widehat{\tau}\otimes\rho'\circ\widehat{\tau} )({\Rt})\nonumber\\
 & =P_{V,W}\circ \E_{\rho\circ\widehat{\tau},\rho'\circ\widehat{\tau}}\circ(\rho\otimes\rho')({\Rt}).\label{eq:Fc}
		\end{align}
		Suppose that $\rho(H_i)v=\lambda_iv$ and $\rho'(H_i)w=\mu_iw$, then $\rho\circ\widehat{\tau}(H_i)v=\rho(H_{\tau(i)})v=\lambda_{\tau(i)}v$ and similarly $\rho'\circ\widetilde{\tau}(H_i)w=\mu_{\tau(i)}w$. Therefore, \begin{align*}
			\E_{\rho\circ\widehat{\tau},\rho'\circ\widehat{\tau}}(v\otimes w)
&=\zeta^{{\sum_{ij}(A^{-1})_{ij}\lambda_{\tau(i)}\mu_{\tau(j)}}}(v\otimes w)\\
&=\zeta^{{\sum_{ij}(A^{-1})_{\tau^{-1}(i)\tau^{-1}(j)}\lambda_{i}\mu_{j}}}(v\otimes w)=\E_{\rho,\rho'}(v\otimes w)
		\end{align*}
		by invariance of the Cartan matrix under $\tau$. Continuing from~\eqref{eq:Fc}, \begin{align*}
			\mathsf{F}\bigl(\braid_{\widetilde{\tau}(V),\widetilde{\tau}(W)}\bigr)=P\circ \E_{\rho\circ\widehat{\tau},\rho'\circ\widehat{\tau}}\circ(\rho\otimes \rho')({\Rt})=\mathsf{F}\bigl(\braid_{V,W}\bigr).
		\end{align*} Thus, $\widetilde{\tau}\bigl(\braid_{V,W}\bigr)=\braid_{\widetilde{\tau}(V),\widetilde{\tau}(W)}$.
		
In Proposition~\ref{lem:rib}, we expressed the ribbon structure of $\UHwmod$ in terms of the braiding and pivotal action by $\tr_R\bigl(\braid_{V,V}\bigr)=\theta_{V}$. Therefore, $\widetilde{\tau}(\theta_{{V}})=\theta_{\widetilde{\tau}(V)}$ and $\widetilde{\tau}$ is an automorphism of $\UHwmod$ as a ribbon category.
	\end{proof}
\begin{Corollary}\label{cor:symm}
		Let $\tau$ be an automorphism of the Dynkin diagram of $\mathfrak{sl}_3$ and $\mathcal{L}$ any link. Then there is a symmetry of the polynomial,
	\begin{align*}
	\Delta_{\mathfrak{sl}_3}(\mathcal{L})(t_1, t_2)=\Delta_{\mathfrak{sl}_3}(\mathcal{L})(t_{\tau(1)}, t_{\tau(2)}).
	\end{align*}
\end{Corollary}
\begin{proof}
	As above, $\tau$ induces the automorphism $\widetilde{\tau}$ on $\UHwmod$. In a slight abuse of notation, we will also use $\widetilde{\tau}$ to denote the automorphism on $\Uwmod$. Let $\mvec{\tau t}$ denote $(t_{\tau(1)},t_{\tau(2)})$. If $v_0$ is the highest weight vector in $\Vt$, then $\rho\circ \widehat{\tau}(K_i)v_0=\rho(K_{\tau(i)})v_0=t_{\tau(i)}v_0$. Thus, $\widetilde{\tau}(\Vt)=V^H(\mvec{\tau t})$.
	
	Let $\mathcal{L}$ be a framed link with $(1,1)$-tangle representative $\mathcal{L}^{{\rm cut}}$ and $\Vt$ an irreducible representation. Since $F_{\Vt}(\mathcal{L}^{{\rm cut}})$ is given by a composition of normalized braidings, evaluations, and coevaluations, Lemma~\ref{lem:symm} implies $\widetilde{\tau}\circ F_{\Vt}(\mathcal{L}^{{\rm cut}})=F_{\widetilde{\tau}'(\Vt)}(\mathcal{L}^{{\rm cut}})=F_{V^H(\mvec{\tau t})}(\mathcal{L}^{{\rm cut}})$. Applying the forgetful functor $\mathsf{F}$, we have the equality of linear maps $\mathsf{F}\circ \widetilde{\tau}\circ F_{\Vt}(\mathcal{L}^{{\rm cut}})=\mathsf{F}\circ F_{V^H(\mvec{\tau t})}(\mathcal{L}^{{\rm cut}})$. This now implies the equality
 $\Delta_{\mathfrak{sl}_3}(\mathcal{L})(t_1,t_2)=\Delta_{\mathfrak{sl}_3}(\mathcal{L})(t_{\tau(1)},t_{\tau(2)})$.
\end{proof}

	\begin{Lemma}\label{lem:invert}
 Let $\mathcal{L}$ be an oriented link and $-\mathcal{L}$ the same link with all orientations reversed. Then
		\[\Delta_{\mathfrak{sl}_3}(-\mathcal{L})(t_1,t_2)=\Delta_{\mathfrak{sl}_3}(\mathcal{L})\bigl(-t_1^{-1},-t_2^{-1}\bigr).\]
	\end{Lemma}
	\begin{proof}
		From Remark \ref{rem:dual}, $\Vt^*\cong V\bigl(\mvec{-t^{-1}}\bigr)$ and by Theorem~\ref{thm:ambi} both $\Vt$ and its dual are ambidextrous for typical $\t$. Since the morphism assigned to the open Hopf link colored by both~$\Vt$ and~$\Vt^*$ is nonzero, we may apply \cite[Proposition~19]{GPT}. Thus, reversing the orientation of a component of~$\mathcal{L}$ is equivalent to coloring it by~$\Vt^*$. Therefore, $\Delta_{\mathfrak{sl}_3}(-\mathcal{L})$ is computed from coloring all components of~$\mathcal{L}$ by~$V\bigl(\mvec{-t^{-1}}\bigr)$.
	\end{proof}

\begin{Remark}
 For every link $\mathcal{L}$, the inversion symmetry
 \begin{align*}
 \Delta_{\mathfrak{sl}_3}(\mathcal{L})(t_1,t_2)=\Delta_{\mathfrak{sl}_3}(\mathcal{L})\bigl(-t_1^{-1}, -t_2^{-1}\bigr)
 \end{align*}
 together with Lemma~\ref{lem:invert} implies $\Delta_{\mathfrak{sl}_3}$ does not detect link inversion.
\end{Remark}	

	\subsection{Skein relation} The skein relation and values of $\Delta_{\mathfrak{sl}_3}$ on two strand torus knots are both derived from the characteristic (minimal) polynomial of $\redbraid_{\t}$. The former is obtained from \eqref{eqn:charpoly}, and the latter is stated in Theorem~\ref{thm:torusknots}.
	\begin{Proposition}\label{prop:skein}
		There is a nine-term skein relation for $\Delta_{\mathfrak{sl}_3}$.
	\end{Proposition}
	\begin{proof}
		Let $r$ be the $8\times8$ matrix which appears in Corollary \ref{cor:Rdecomp}. By the Cayley--Hamilton theorem, the characteristic polynomial of $r$ determines a relation among powers of itself. Therefore, the characteristic polynomial of $\redbraid_{\t}$ is the characteristic polynomial of $r$ raised to the power $\dim \Vt$. Thus, $\redbraid_{\t}$ is a solution to the equation given by $r$. On $\Vt^{\otimes 2}$, this relation takes the form
		\begin{align}\label{eqn:charpoly}
		\bigl(\redbraid_{\t}^2+\mathrm{id}\bigr)\bigl(t_1^2\mathrm{id}+\redbraid_{\t}\bigr) \bigl(t_1^2\redbraid_{\t}+\mathrm{id}\bigr)\bigl(t_2^2\mathrm{id}+\redbraid_{\t}\bigr)\bigl(t_2^2\redbraid_{\t} +\mathrm{id}\bigr)\bigl(t_1^2t_2^2\mathrm{id}-\redbraid_{\t}\bigr)\bigl(t_1^2t_2^2\redbraid_{\t}-\mathrm{id}\bigr)=0.
		\end{align}
		After expansion and normalization, this implies the palindromic relation
		\begin{align*}
		c_0\mathrm{id}_{\Vt^{\otimes 2}}+\sum_{i=1}^{4}c_{i}\bigl(\redbraid_{\t}^{i}+(\redbraid_{\t})^{-i}\bigr)=0,
		\end{align*} where
 \begin{gather*}
 c_0=-2 \cdot\frac{t_1^8 t_2^6 + t_1^6 t_2^8 - t_1^6 t_2^6 + t_1^6 t_2^4 - t_1^6 t_2^2 + t_1^4 t_2^6 - 3 t_1^4 t_2^4 + t_1^4 t_2^2 - t_1^2 t_2^6 + t_1^2 t_2^4 - t_1^2 t_2^2 + t_1^2 + t_2^2}{t_1^4 t_2^4},
 \\
 c_1 = -\frac{t_1^8 t_2^8\hspace{-0.5pt} + t_1^8 t_2^4\hspace{-0.5pt} + 3 t_1^6 t_2^6\hspace{-0.5pt} - 3 t_1^6 t_2^4\hspace{-0.5pt} + t_1^4 t_2^8\hspace{-0.5pt} - 3 t_1^4 t_2^6\hspace{-0.5pt} + 2 t_1^4 t_2^4\hspace{-0.5pt} - 3 t_1^4 t_2^2\hspace{-0.5pt} + t_1^4\hspace{-0.5pt} - 3 t_1^2 t_2^4\hspace{-0.5pt} + 3 t_1^2 t_2^2\hspace{-0.5pt} + t_2^4\hspace{-0.5pt} + 1}{t_1^4 t_2^4},
 \\
 c_2 = -\frac{\bigl(t_1^4 t_2^2 + t_1^2 t_2^4 - t_1^2 t_2^2 - 1\bigr) \bigl(t_1^4 t_2^4 + t_1^2 t_2^2 - t_1^2 - t_2^2\bigr)}{t_1^4 t_2^4},
 \\
 c_3 = -\frac{t_1^4 t_2^4 - t_1^4 t_2^2 - t_1^2 t_2^4 - t_1^2 - t_2^2 + 1}{t_1^2 t_2^2},
 \\
 c_4 = 1,
 \end{gather*}
 as determined by \eqref{eqn:charpoly}. Replacing each factor of $\redbraid_{\t}$ with a diagrammatic strand crossing and $\mathrm{id}_{\Vt^{\otimes 2}}$ by two vertical strands, we obtain the skein relation.
	\end{proof}
	
	Similar to how we used the characteristic polynomial of the braiding to determine the skein relation, other characteristic polynomials yield relations among families of torus knots. Let $q$ be a prime number, and $r$ any positive integer less than $q$. Then for each $0<n<q$, we have that $qn+r$ and $q$ are coprime. Define
	\begin{align*}
	\beta_q=\left(\prod_{i=0}^{q-2}\mathrm{id}^{\otimes i}\otimes \redbraid_{\t}\otimes \mathrm{id}^{\otimes q-i-2}\right),
	\end{align*}
	which acts on $\Vt^{\otimes q}$. Then the characteristic polynomial of $\beta_q^q$ is some equation of the form
	\begin{equation}\label{eqn:char}
	\sum_{i=0}^{8^q} a_i\beta_q^{qi}=0.
	\end{equation}
	Multiplying this equation by $\beta_q^r$ implies that the invariants of the torus knots of types $(r,q),(q+r,q),\dots,((8^q-1)q+r,q)$ determine the invariant for the $(8^qq+r,q)$ torus knot. With this information and after multiplying equation \eqref{eqn:char} by $\beta_q^{r+1}$, we can deduce the invariant for the $((8^q+1)q+r,q)$ torus knot and so on. This implies a recursion relation for all torus knots $T_{nq+r,q}$, which can then be converted to an explicit function of $n$. The resulting expression for the $q=2,r=1$ case is stated as a theorem below.

	\begin{Theorem}\label{thm:torusknots}
 The value of $\Delta_{\mathfrak{sl}_3}$ on a $(2n+1,2)$ torus knot is given by
 \begin{align*}
		\frac{
			\bigl(t_1-t_1^{-1}\bigr)\bigl(t_1^{4n+2}+t_1^{-(4n+2)}\bigr)} {\bigl(t_2+t_2^{-1}\bigr)\bigl(t_1^2+t_1^{-2}\bigr)\bigl(t_1t_2-t_1^{-1}t_2^{-1}\bigr)}&+
		\frac{	\bigl(t_2-t_2^{-1}\bigr)\bigl(t_2^{4n+2}+t_2^{-(4n+2)}\bigr)}{\bigl(t_1+t_1^{-1}\bigr) \bigl(t_2^2+t_2^{-2}\bigr)\bigl(t_1t_2-t_1^{-1}t_2^{-1}\bigr)}\\
&+
		\frac{\bigl(t_1t_2+t_1^{-1}t_2^{-1}\bigr) \bigl(t_1^{4n+2}t_2^{4n+2}+t_1^{-(4n+2)}t_2^{-(4n+2)}\bigr)}{\bigl(t_1^2t_2^2+t_1^{-2}t_2^{-2}\bigr) \bigl(t_1+t_1^{-1}\bigr)\bigl(t_2+t_2^{-1}\bigr)}.
		\end{align*}
\end{Theorem}
		Observe that the expression for these torus knots can be separated into three terms: one pair of terms exchange the roles of~$t_1$ and~$t_2$, the other is symmetric in~$t_1$ and~$t_2$.

	\section{Values of \texorpdfstring{$\boldsymbol{\Delta_{\mathfrak{sl}_3}}$}{Delta sl\_3}}\label{sec:values}
	
	In this section, we give the value of the unrolled restricted quantum $\mathfrak{sl}_3$ invariant for all prime knots with at most seven crossings, as well as some other examples.
 Among these examples are knots that compare $\Delta_{\mathfrak{sl}_3}$ to other well-known invariants. The HOMFLY polynomial does not distinguish the knot $\mathsf{11_{n34}}$ from $\mathsf{11_{n42}}$ nor does it distinguish $\mathsf{5_1}$ and $\mathsf{10_{132}}$ but $\Delta_{\mathfrak{sl}_3}$ does. The Jones polynomial differentiates $\mathsf{6_1}$ and $\mathsf{9_{46}}$, but $\Delta_{\mathfrak{sl}_3}$ does not. The Jones polynomial and the $\mathfrak{sl}_3$ invariant both distinguish $\mathsf{8_9}$ from $\mathsf{10_{155}}$; however, the Alexander polynomial does not. These data for the HOMFLY, Jones, and Alexander polynomials are from \cite{TKA}.

 We also refer to \cite{TKA} for braid presentations of prime knots. These invariants were computed locally with Python (\textsc{SymPy 1.14.0}), and previously using \textsc{Maple 2018.0} with the unity high performance computing cluster at The Ohio State University. Both sets of code used to produce these invariants are available on the author's GitHub repository \cite{code}.
	
 By the symmetry results of Section~\ref{sec:properties}, it is enough to specify the coefficient of $t_1^{2a}t_2^{2b}$ in $\Delta_{\mathfrak{sl}_3}(\mathcal{L})$ for each $(a,b)$ in the cone
	\begin{align*}
	C=\bigl\{(a,b)\in \mathbb{Z}^2\mid a\geq 0 \text{ and } |b|\leq a\bigr\}.
	\end{align*}
	The coefficients of various knots can be found in Figures \ref{fig:sl3leq7} and \ref{fig:sl3geq8} below.
	We have boxed the leftmost value on each cone, it has coordinates $(0,0)$ and is the constant term in the polynomial invariant for the indicated knot. From the values given, we can reconstruct $\Delta_{\mathfrak{sl}_3}$ since the coefficient in position $(a,b)$ is equal to those in positions $(b,a)$, $(-a,-b)$, and $(-b,-a)$. For example, the polynomial
 \begin{align*}
	\Delta_{\mathfrak{sl}_3}(\mathsf{3_1})(t_1,t_2) ={} & \bigl(t_1^4t_2^4+t_1^{-4}t_2^{-4}\bigr) -\bigl(t_1^4t_2^2+t_1^2t_2^4+t_1^{-4}t_2^{-2}+t_1^{-2}t_2^{-4}\bigr)+\bigl(t_1^4+t_2^4+t_1^{-4}+t_2^{-4}\bigr)\\
	&+2\bigl(t_1^2t_2^2+t_1^{-2}t_2^{-2}\bigr)-2\bigl(t_1^2+t_2^2+t_1^{-2}+t_2^{-2}\bigr)+ \bigl(t_1^2t_2^{-2}+t_1^{-2}t_2^2\bigr)+1
	\end{align*}
 is determined from the entries given in Figure~\ref{fig:sl3leq7}.
	
\begin{figure}[!ht]\centering
		\savebox{\largestimage}{\includegraphics{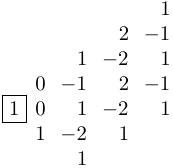}}%
		\begin{subfigure}[t]{.15\linewidth}
			\centering
			\raisebox{\dimexpr.33\ht\largestimage-.33\height}
			{\includegraphics[]{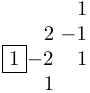}}
			\caption{$\mathsf{3_1}$}
			\label{fig:3_1}
		\end{subfigure}
		\begin{subfigure}[t]{.15\linewidth}
			\centering
		 \raisebox{\dimexpr.33\ht\largestimage-.33\height}{				\includegraphics[]{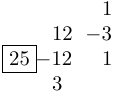}}
			\caption{$\mathsf{4_1}$}
			\label{fig:4_1}
		\end{subfigure}
		\begin{subfigure}[t]{.25\linewidth}
			\centering
				\usebox{\largestimage}
			\caption{$\mathsf{5_1}$}
			\label{fig:5_1}
		\end{subfigure}
		\begin{subfigure}[t]{.15\linewidth}
			\centering
		\raisebox{\dimexpr.33\ht\largestimage-.33\height}{				\includegraphics[]{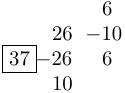}}
		\caption{$\mathsf{5_2}$}
		\label{fig:5_2}
	\end{subfigure}
\begin{subfigure}[t]{.15\linewidth}
	\centering
	\raisebox{\dimexpr.33\ht\largestimage-.33\height}{				\includegraphics[]{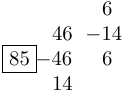}}
	\caption{$\mathsf{6_1}$}
	\label{fig:6_1}
\end{subfigure}
	\savebox{\largestimage}{\includegraphics{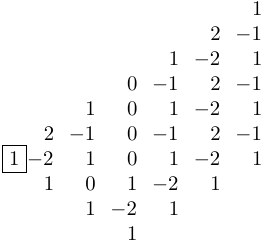}}%
\begin{subfigure}[t]{.3\linewidth}
	\centering
	\raisebox{\dimexpr.33\ht\largestimage-.33\height}{				\includegraphics[]{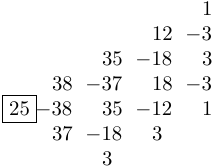}}
	\caption{$\mathsf{6_2}$}
	\label{fig:6_2}
\end{subfigure}
\begin{subfigure}[t]{.3\linewidth}
	\centering
	\raisebox{\dimexpr.33\ht\largestimage-.33\height}{				\includegraphics[]{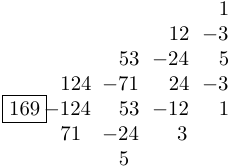}}
	\caption{$\mathsf{6_3}$}
	\label{fig:6_3}
\end{subfigure}
\begin{subfigure}[t]{.35\linewidth}
	\centering
	\raisebox{\dimexpr.33\ht\largestimage-.33\height}{				\includegraphics[]{71.pdf}}
	\caption{$\mathsf{7_1}$}
	\label{fig:71}
\end{subfigure}\\
	\savebox{\largestimage}{\includegraphics{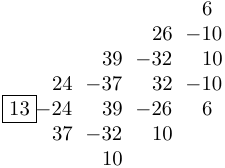}}%
\begin{subfigure}[t]{.15\linewidth}
	\centering
	\raisebox{\dimexpr.33\ht\largestimage-.33\height}{				\includegraphics[]{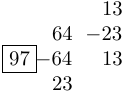}}
	\caption{$\mathsf{7_2}$}
	\label{fig:72}
\end{subfigure}
\begin{subfigure}[t]{.3\linewidth}
	\centering
	\raisebox{\dimexpr.33\ht\largestimage-.33\height}{				\includegraphics[]{73.pdf}}
	\caption{$\mathsf{7_3}$}
	\label{fig:73}
\end{subfigure}
\begin{subfigure}[t]{.16\linewidth}
	\centering
	\raisebox{\dimexpr.33\ht\largestimage-.33\height}{				\includegraphics[]{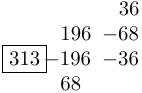}}
	\caption{$\mathsf{7_4}$}
	\label{fig:74}
\end{subfigure}
\begin{subfigure}[t]{.3\linewidth}
	\centering
	\raisebox{\dimexpr.33\ht\largestimage-.33\height}{				\includegraphics[]{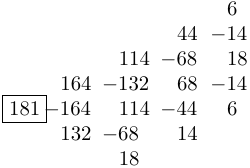}}
	\caption{$\mathsf{7_5}$}
	\label{fig:75}
\end{subfigure}\\
\begin{subfigure}[t]{.3\linewidth}
	\centering
	\raisebox{\dimexpr.33\ht\largestimage-.33\height}{				\includegraphics[]{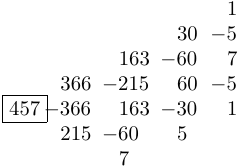}}
	\caption{$\mathsf{7_6}$}
	\label{fig:76}
\end{subfigure}
\begin{subfigure}[t]{.3\linewidth}
	\centering
	\raisebox{\dimexpr.33\ht\largestimage-.33\height}{				\includegraphics[]{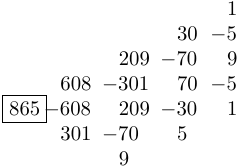}}
	\caption{$\mathsf{7_7}$}
	\label{fig:77}
\end{subfigure}
		\caption{The value of $\Delta_{\mathfrak{sl}_3}$ for all prime knots with at most seven crossings.
		}\label{fig:sl3leq7}
	\end{figure}
	
\begin{figure}[!ht]
	\savebox{\largestimage}{\includegraphics{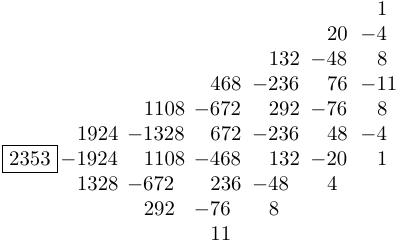}}%
	\begin{subfigure}[t]{.45\linewidth}
		\centering
		\raisebox{\dimexpr.33\ht\largestimage-.33\height}
		{\includegraphics[]{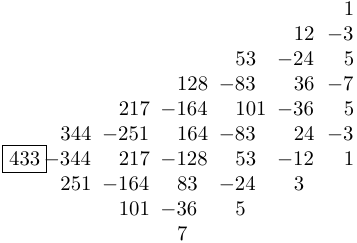}}
		\caption{$\mathsf{8_9}$}
		\label{fig:8_9}
	\end{subfigure}
\begin{subfigure}[t]{.45\linewidth}
	\centering
	\raisebox{\dimexpr.33\ht\largestimage-.33\height}
	{\includegraphics[]{817.pdf}}
	\caption{$\mathsf{8_{17}}$}
	\label{fig:8_17}
\end{subfigure}\\
\savebox{\largestimage}{\includegraphics{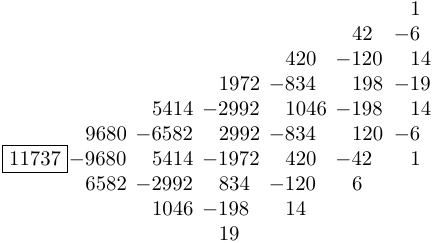}}%
\begin{subfigure}[t]{.45\linewidth}
\centering
\raisebox{\dimexpr.33\ht\largestimage-.33\height}
{\includegraphics[]{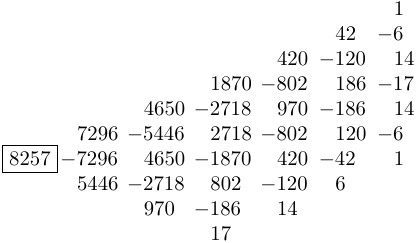}}
\caption{$\mathsf{9_{32}}$}
\label{fig:9_32}
\end{subfigure}
\begin{subfigure}[t]{.5\linewidth}
\centering
\raisebox{\dimexpr.33\ht\largestimage-.33\height}
{\includegraphics[]{933.pdf}}
\caption{$\mathsf{9_{33}}$}
\label{fig:9_33}
\end{subfigure}\\
\savebox{\largestimage}{\includegraphics{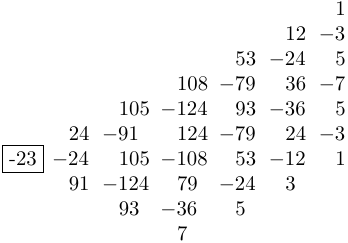}}%
\begin{subfigure}[t]{.15\linewidth}
	\centering
	\raisebox{\dimexpr.33\ht\largestimage-.33\height}
	{\includegraphics[]{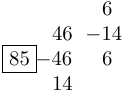}}
	\caption{$\mathsf{9_{46}}$}
	\label{fig:9_46}
\end{subfigure}
\begin{subfigure}[t]{.35\linewidth}
	\centering
	\raisebox{\dimexpr.33\ht\largestimage-.33\height}
	{\includegraphics[]{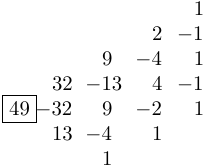}}
	\caption{$\mathsf{10_{132}}$}
	\label{fig:10_132}
\end{subfigure}
\begin{subfigure}[t]{.37\linewidth}
	\centering
	\raisebox{\dimexpr.33\ht\largestimage-.33\height}
	{\includegraphics[]{10155.pdf}}
	\caption{$\mathsf{10_{155}}$}
	\label{fig:10_155}
\end{subfigure}
\\
\savebox{\largestimage}{\includegraphics{1134.pdf}}%
\begin{subfigure}[t]{.45\linewidth}
	\centering
	{\includegraphics[]{1134.pdf}}
	\caption{$\mathsf{11_{n34}}$}
	\label{fig:11_34}
\end{subfigure}
\begin{subfigure}[t]{.3\linewidth}
	\centering
	\raisebox{\dimexpr.5\ht\largestimage-.5\height\relax}
	{\includegraphics[]{1142.pdf}}
	\caption{$\mathsf{11_{n42}}$}
	\label{fig:11_42}
\end{subfigure}
\begin{subfigure}[t]{.2\linewidth}
	\centering
	\raisebox{\dimexpr.5\ht\largestimage-.75\height\relax}
	{\includegraphics[]{Wh31.pdf}}
	\caption{$\text{Wh}^0(\mathsf{3_{1}})$}
\end{subfigure}
\caption{The value of $\Delta_{\mathfrak{sl}_3}$ for some prime knots with more than seven crossings.}\label{fig:sl3geq8}
	\end{figure}

In Figure~\ref{fig:sl3geq8}, observe that the mutant pair $\mathsf{11_{n34}}$ and $\mathsf{11_{n42}}$ have different $\mathfrak{sl}_3$ polynomials. In \cite[Lemma~1 and p.~199]{GPT}, three sufficient criteria for ambidexterity of a module $V$ are given. One of their criteria, given in Lemma~1, is that the braiding on $V\otimes V$ is central in $\End_{\overline{U}^H} (V\otimes V)$. If the braiding were central, then the corresponding link invariant could not detect mutation \cite[Theorem~5]{MortonCromwell}.

\appendix

 \section{Proof of Proposition~\ref{prop:quasiR}}\label{sec:quasiR}
 We show that ${\Rt}$ is a quasi-$R$-matrix.

\begin{proof}
	We first prove that ${\Rt}\Delta(x)=\Psi_\zeta(\Delta^{\rm op}(x)){\Rt}$ for all $x\in \overline{U}^H$.
	We then show that $(\Psi_{\zeta})_{23}(\Rt_{13})\Rt_{23}=(\Delta\otimes 1)({\Rt})$ and $(\Psi_{\zeta})_{12}(\Rt_{13})\Rt_{12}=(1\otimes \Delta)({\Rt})$ both hold. We give an explicit computation proving the former, while the latter is a similar computation.

	To prove ${\Rt}\Delta(x)=\Psi_\zeta(\Delta^{\rm op}(x)){\Rt}$, we first note that $K_i$ and $H_i$ have symmetric coproducts which are preserved by $\Psi_\zeta$ and commute with ${\Rt}$. For these generators the relation holds trivially.

 It is sufficient to now consider only $x=E_1$, as the computation is similar on other root generators. It will then follow that the relation holds for all $x\in\overline{U}^H$. We verify that
	\begin{align*}
		{\Rt}(E_1\otimes K_1+1\otimes E_1)={\Rt}\Delta(E_1)=\Psi_\zeta(\Delta^{\rm op}(E_1)){\Rt}=(E_1\otimes K_1^{-1}+1\otimes E_1){\Rt}.
	\end{align*}
	Computing each term of ${\Rt}\Delta(E_1)$ directly yields
	\begin{align*}
		{\Rt}(E_1\otimes K_1)={} & (1+2\zeta E_1\otimes F_1)(1+2\zeta E_{12}\otimes F_{12})\\
& \times (E_1\otimes 1
 +2\zeta (\zeta E_1E_2+\zeta E_{12})\otimes F_2)(1\otimes K_1)\\
		={} &
		(E_1\otimes 1)(1+2\zeta (-\zeta E_{12})\otimes F_{12})( 1-2E_2\otimes F_2)(1\otimes K_1)
 \\&+
		(1+2\zeta E_1\otimes F_1)(-2E_{12}\otimes F_2)(1\otimes K_1)\\
		={} &
		(E_1\otimes 1)(1+2 E_{12}\otimes F_{12})( 1-2E_2\otimes F_2)(1\otimes K_1)
 \\&+(1\otimes K_1)(1-2\zeta E_1\otimes F_1)(2\zeta E_{12}\otimes F_2)\\
		={} &
		(E_1\otimes K_1)(1+2\zeta E_{12}\otimes F_{12})( 1+2\zeta E_2\otimes F_2)
 \\&+
		(1\otimes K_1)(1-2\zeta E_1\otimes F_1)(2\zeta E_{12}\otimes F_2)
	\end{align*}
	and
	\begin{align*}
		{\Rt}(1\otimes E_1)={}& (1+2\zeta E_1\otimes F_1)(1+2\zeta E_{12}\otimes (E_1F_{12}-\zeta F_2K_1))(1+2\zeta E_2\otimes F_2)\\
		={}&
		(1\otimes E_1){\Rt}+(-2\zeta E_1\otimes \floor{K_1})(1+2\zeta E_{12}\otimes F_{12})(1+2\zeta E_2\otimes F_2)\\
& +
		(1+2\zeta E_1\otimes F_1)(2E_{12}\otimes F_2K_1)\\
		={}&
		(1\otimes E_1){\Rt}+(- E_1\otimes (K_1-K_1^{-1}))(1+2\zeta E_{12}\otimes F_{12})(1+2\zeta E_2\otimes F_2)\\
& +
		(1\otimes K_1)(1-2\zeta E_1\otimes F_1)(-2\zeta E_{12}\otimes F_2).
	\end{align*} Thus,
	\begin{align*}
		{\Rt}(E_1\otimes K_1+1\otimes E_1)& =
		(1\otimes E_1){\Rt}+\bigl( E_1\otimes K_1^{-1}\bigr)(1+2\zeta E_{12}\otimes F_{12})(1+2\zeta E_2\otimes F_2)\\
		& =
		\bigl(E_1\otimes K_1^{-1}+1\otimes E_1\bigr){\Rt}.
	\end{align*}
	
	To prove the next condition, we observe
	\begin{align*}
		(\Psi_{\zeta})_{23}(\Rt_{13})\Rt_{23}
		=
		\prod_{\beta\in\proots}(1+ 2\zeta E_{\beta}\otimes K_{\beta}\otimes F_{\beta})\prod_{\beta\in\proots}(1+ 2\zeta\otimes E_{\beta}\otimes F_{\beta}).
	\end{align*}
	For simple roots $\alpha$, \begin{align*}
		(1+2\zeta E_\alpha\otimes K_\alpha\otimes F_\alpha)(1+2\zeta \otimes E_\alpha\otimes F_\alpha)=(\Delta\otimes 1)(1+2\zeta E_\alpha\otimes F_\alpha)
	\end{align*}
	and for $\alpha=\alpha_{12}$,
	\begin{align*}
		& (1+2\zeta E_\alpha\otimes K_\alpha\otimes F_\alpha)(1+2\zeta 1\otimes E_\alpha\otimes F_\alpha)\\
&\qquad =(\Delta\otimes 1)(1+2\zeta E_\alpha\otimes F_\alpha)+4\zeta E_2\otimes E_1K_2\otimes F_{12}.
	\end{align*}
	We commute the terms appearing in $(\Psi_{\zeta})_{23}(\Rt_{13})\Rt_{23}$ so that the above product expressions for the coproduct appear and simplify to $(\Delta\otimes 1)({\Rt})$. The following equalities are readily verified:
	\begin{align*}
		&[1+2\zeta E_2\otimes K_2\otimes F_2, 1+2\zeta \otimes E_1\otimes F_1]=-4\zeta E_2\otimes E_1K_2\otimes F_{12},\\
		&[1+2\zeta E_{12}\otimes K_1K_2\otimes F_{12}, 1+2\zeta \otimes E_1\otimes F_1]=0,\\
		&[1+2\zeta E_2\otimes K_2\otimes F_2, 1+2\zeta \otimes E_{12}\otimes F_{12}]=0.
	\end{align*}
	Thus,
	\begin{align*}
	\Psi_{\zeta,23}(\Rt_{13})\Rt_{23}
		={}&
		\prod_{\beta\in\proots}\bigl(1^{\otimes 3}+ 2\zeta E_{\beta}\otimes K_{\beta}\otimes F_{\beta}\bigr)\prod_{\beta\in\proots}\bigl(1^{\otimes 3}+ 2\zeta\otimes E_{\beta}\otimes F_{\beta}\bigr)\\
		={}& \prod_{\beta\in\{1,12\}}(1+ 2\zeta E_{\beta}\otimes K_{\beta}\otimes F_{\beta})\\
& \times ( (1+2\zeta \otimes E_1\otimes F_1)(1+2\zeta E_2\otimes K_2\otimes F_2)-4\zeta E_2\otimes E_1K_2\otimes F_{12} )\\
& \times
\prod_{\alpha\in\{12,2\}}(1+ 2\zeta\otimes E_{\beta}\otimes F_{\beta})\\
		={}& (\Delta\otimes \mathrm{id})(1+2\zeta E_1\otimes F_1)\\
& \times
((\Delta\otimes \mathrm{id})(1+2\zeta E_{12}\otimes F_{12})+4\zeta E_2\otimes E_1K_2\otimes F_{12}
)\\
& \times (\Delta\otimes \mathrm{id})(1+2\zeta E_2\otimes F_2)-4\zeta (1+ 2\zeta E_{1}\otimes K_{1}\otimes F_{1})\\
& \times(E_2\otimes E_1K_2\otimes F_{12}) (1+ 2\zeta\otimes E_{2}\otimes F_{2})\\
		={}& (\Delta\otimes \mathrm{id})(1+2\zeta E_1\otimes F_1)(\Delta\otimes \mathrm{id})(1+2\zeta E_{12}\otimes F_{12})\\
& \times
(\Delta\otimes \mathrm{id})(1+2\zeta E_2\otimes F_2)\\
={}& (\Delta\otimes \mathrm{id})({\Rt}).
	\end{align*}
	This gives the desired equality.
\end{proof}

\section{Proof of Lemma~\ref{lem:intertwineraction}}\label{sec:intertwineraction}
The proof follows the structure of \cite[Lemma A.17]{Ohtsuki}.
\begin{proof}
	Recall that any $f\in\End_{\overline{U}}\bigl(W_\alpha(\t)^{\otimes2}\bigr)$ is expressible as a sum of scalars acting on each summand of the tensor product decomposition. We use the notation of \eqref{eq:EndDecomp}. For such $f$ we also have that \smash{$\tr_R(f)=\trace_2((\mathrm{id}_{W_\alpha(\t)}\otimes h_{W_\alpha(\t)})\cdot f)\in \End_{\overline{U}^H}(W_\alpha(\t))$}. For each $\alpha$, assume $\t\in \R_\alpha$ so that $W_\alpha(\t)$ is an irreducible representation, which implies that~$\tr_R(f)$ is a scalar multiple of the identity. Therefore, it is sufficient to compute its action on a highest weight vector of~$W_\alpha(\t)$.
	
It is straightforward to verify the following equalities hold in $W_i(\t)\otimes W_i(\t)$:
	\begin{align*}
		&w_0^{i,\t}\otimes F_jw_0^{i,\t}
=\dfrac{\floor{t_j}}{\floor{t_j^2}}\Delta(F_j)\bigl(w_0^{i,\t}\otimes w_0^{i,\t}\bigr)-
		\dfrac{1}{\floor{t_j^2}}\Delta(E_j)\bigl(F_jw_0^{i,\t}\otimes F_jw_0^{i,\t}\bigr),\\
	&	w_0^{i,\t}\otimes F_iF_jw_0^{i,\t}
=t_i\left(\dfrac{\floor{t_j}}{\floor{t_j^2}}\Delta(F_iF_j)\bigl(w_0^{i,\t}\otimes w_0^{i,\t}\bigr)-
		\dfrac{1}{\floor{t_j^2}}\Delta(F_iE_j)\bigl(F_jw_0^{i,\t}\otimes F_jw_0^{i,\t}\bigr)\right),\\
		&w_0^{i,\t}\otimes F_jF_iF_jw_0^{i,\t}=
		\dfrac{t_i}{2\floor{\zeta t_j^2}}\Delta(F_jF_iF_j)\bigl(w_0^{i,\t}\otimes w_0^{i,\t}\bigr)\\
&\hphantom{w_0^{i,\t}\otimes F_jF_iF_jw_0^{i,\t}=}{}
		-\dfrac{t_i}{2\floor{t_j}\floor{t_j^2}}\Delta(F_iF_jE_j)\bigl(F_jw_0^{i,\t}\otimes F_jw_0^{i,\t}\bigr)
		\\
&\hphantom{w_0^{i,\t}\otimes F_jF_iF_jw_0^{i,\t}=}{}
+\dfrac{2}{\floor{t_j^4}}\Delta(E_jE_iE_j)\bigl(F_jF_iF_jw_0^{i,\t}\otimes F_jF_iF_jw_0^{i,\t}\bigr).
	\end{align*}
Thus,
\begin{align*}
	&f\bigl(	w_0^{i,\t}\otimes w_0^{i,\t}\bigr)=f_+\bigl(w_0^{i,\t}\otimes w_0^{i,\t}\bigr),\\
	&f\bigl(	w_0^{i,\t}\otimes F_jw_0^{i,\t}\bigr)
=\left(\dfrac{\floor{t_j}}{t_j\floor{t_j^2}}f_++\dfrac{t_j\floor{t_j}}{\floor{t_j^2}}f_V\right) \bigl(w_0^{i,\t}\otimes F_jw_0^{i,\t}\bigr) +\cdots,\\
	&f\bigl(	w_0^{i,\t}\otimes F_iF_jw_0^{i,\t}\bigr)
=\left(\dfrac{\floor{t_j}}{t_j\floor{t_j^2}}f_++
	\dfrac{t_j\floor{t_j}}{\floor{t_j^2}}f_V\right)\bigl(	w_0^{i,\t}\otimes F_iF_jw_0^{i,\t}\bigr)+\cdots,\\
	&f\bigl(	w_0^{i,\t}\otimes F_jF_iF_jw_0^{i,\t}\bigr)
=\left(\dfrac{1}{2t_j^2\floor{\zeta t_j^2}}f_+
	+\dfrac{2\floor{\zeta t_j}\floor{t_j}t_j^2}{\floor{t_j^4}}f_-\right) \bigl(w_0^{i,\t}\otimes F_jF_iF_jw_0^{i,\t}\bigr)+\cdots
\end{align*}
with ``$\cdots$'' indicating terms that are outside the span of the given vector, i.e., off-diagonal matrix entries. We can see that $\tr_R(f)$ acts as multiplication by \begin{gather*}
	t_j^{-2}\left(\begin{matrix}
	 f_+-\left(\dfrac{\floor{t_j}}{t_j\floor{t_j^2}}f_++\dfrac{t_j\floor{t_j}}{\floor{t_j^2}}f_V\right)
 \vspace{1mm}\\
 +
 \left(\dfrac{\floor{t_j}}{t_j\floor{t_j^2}}f_++
\dfrac{t_j\floor{t_j}}{\floor{t_j^2}}f_V\right)-\left(\dfrac{1}{2t_j^2\floor{\zeta t_j^2}}f_+
+\dfrac{2\floor{\zeta t_j}\floor{t_j}t_j^2}{\floor{t^4_j}}f_-\right)
	\end{matrix}\right),
\end{gather*}
which simplifies to the desired scalar.

We now consider $\alpha=\alpha_1+\alpha_2$ and $W_\alpha(\t)$, where $t_1t_2=\zeta \sigma$ and $\sigma ^2=1$. The following equalities are easily verified:
\begin{align*}
	&w_0^{\alpha,\t}\otimes F_1w_0^{\alpha,\t}
=\dfrac{\floor{t_1}}{\floor{t_1^2}}\Delta(F_1)\bigl(w_0^{\alpha,\t}\otimes w_0^{\alpha,\t}\bigr)-
	\dfrac{1}{\floor{t_1^2}}\Delta(E_1)\bigl(F_1w_0^{\alpha,\t}\otimes F_1w_0^{\alpha,\t}\bigr),\\
&	w_0^{\alpha,\t}\otimes F_2w_0^{\alpha,\t}
=\dfrac{\floor{t_2}}{\floor{t_2^2}}\Delta(F_2)\bigl(w_0^{\alpha,\t}\otimes w_0^{\alpha,\t}\bigr)-
	\dfrac{1}{\floor{t_2^2}}\Delta(E_2)\bigl(F_2w_0^{\alpha,\t}\otimes F_2w_0^{\alpha,\t}\bigr)\\
	& w_0^{\alpha,\t}\otimes F_1F_2w_0^{\alpha,\t}
=
	\dfrac{\sigma \floor{\zeta t_1}}{2\floor{t_1}}\Delta(F_1F_2)\bigl(w_0^{\alpha,\t}\otimes w_0^{\alpha,\t}\bigr)
	+\dfrac{ \sigma }{2}\Delta(F_2F_1)\bigl(w_0^{\alpha,\t}\otimes w_0^{\alpha,\t}\bigr)
	\\
&\hphantom{w_0^{\alpha,\t}\otimes F_1F_2w_0^{\alpha,\t}=}{}
+\dfrac{\sigma }{2\floor{t_1}\floor{\zeta t_1^2}}\Delta(F_2E_1)\bigl(F_1w_0^{\alpha,\t}\otimes F_1w_0^{\alpha,\t}\bigr)
 \\
	&
\hphantom{w_0^{\alpha,\t}\otimes F_1F_2w_0^{\alpha,\t}=}{}
-\dfrac{1}{2\floor{t_1}\floor{\zeta t_1^2}}\Delta(F_1E_2)\bigl(F_2w_0^{\alpha,\t}\otimes F_2w_0^{\alpha,\t}\bigr).
\end{align*}
Applying $f$, we have
\begin{align*}
	&f\bigl(w_0^{\alpha,\t}\otimes w_0^{\alpha,\t}\bigr)=f_V\bigl(w_0^{\alpha,\t}\otimes w_0^{\alpha,\t}\bigr),\\
	&f\bigl(w_0^{\alpha,\t}\otimes F_1w_0^{\alpha,\t}\bigr)
=\left(\dfrac{\floor{t_1}}{\floor{t_1^2}t_1}f_V+t_1
	\dfrac{\floor{t_1}}{\floor{t_1^2}}f_+\right)\bigl(w_0^{\alpha,\t}\otimes F_1w_0^{\alpha,\t}\bigr)+\cdots,\\
	&f\bigl(w_0^{\alpha,\t}\otimes F_2w_0^{\alpha,\t}\bigr)
=\left(\dfrac{\floor{t_2}}{\floor{t_2^2}t_2}f_V+
	\dfrac{t_2\floor{t_2}}{\floor{t_2^2}}f_-\right)\bigl(w_0^{\alpha,\t}\otimes F_2w_0^{\alpha,\t}\bigr)+\cdots,\\
	&f\bigl(w_0^{\alpha,\t}\otimes F_1F_2w_0^{\alpha,\t}\bigr)
=\begin{pmatrix}
		\dfrac{\sigma \floor{\zeta t_1}}{2\floor{t_1}t_1t_2}f_V
		+\dfrac{ \sigma \floor{t_1}}{2t_1t_2\floor{\zeta t_1}}f_V
		\vspace{1mm}\\
-\dfrac{\sigma \floor{t_1}t_1\floor{t_1}}{2\floor{t_1}\floor{\zeta t_1^2}t_2\floor{\zeta t_1}}f_+
		+\dfrac{\floor{t_2}t_2}{2\floor{t_1}\floor{\zeta t_1^2}t_1}f_-
	\end{pmatrix}\\
& \hphantom{f\bigl(w_0^{\alpha,\t}\otimes F_1F_2w_0^{\alpha,\t}\bigr)=}{}\times
	\bigl(w_0^{\alpha,\t}\otimes F_1F_2w_0^{\alpha,\t}\bigr)+\cdots.
\end{align*}
Since $h_{W_\alpha(\t)}=-1$, the scalar action of $\tr_R(f)$ is multiplication by
\begin{align*}
	-\begin{pmatrix}
		f_V\left(1-\dfrac{\floor{t_1}}{\floor{t_1^2}t_1}-\dfrac{\floor{t_2}}{\floor{t_2^2}t_2}+\dfrac{\floor{\zeta t_1}}{2\floor{t_1}\zeta}
		+\dfrac{ \floor{t_1}}{2\zeta\floor{\zeta t_1}}\right)\vspace{1mm}\\
		+f_+\left(-t_1
		\dfrac{\floor{t_1}}{\floor{t_1^2}}-\dfrac{\sigma \floor{t_1}t_1}{2\floor{\zeta t_1^2}t_2\floor{\zeta t_1}}\right)
		+f_-\left(
		-\dfrac{t_2\floor{t_2}}{\floor{t_2^2}}+\dfrac{\floor{t_2}t_2}{2\floor{t_1}\floor{\zeta t_1^2}t_1}\right)
	\end{pmatrix}=\dfrac{f_+-f_-}{2\floor{\zeta t_1^2}}.
\end{align*}
This may be written as \smash{$\frac{f_+-f_-}{t_1^2-t_2^2}$}.
\end{proof}

\subsection*{Acknowledgements}

The author is very grateful to Thomas Kerler for numerous insightful discussions. The author also thanks Sergei Chmutov, Sachin Gautam, Nathan Geer, Ben-Michael Kohli, Simon Lentner, Peter Samuelson, Vladimir Turaev, Emmanuel Wagner, and many anonymous referees for their helpful comments and suggestions which have greatly improved the quality of this manuscript. This work was partially supported through the NSF-RTG grants \#DMS-1547357 and \#DMS-2135960. The computation of invariants in Sections~\ref{sec:properties} and~\ref{sec:values} was done in \textsc{Python} (\textsc{SymPy~1.14.0}) and previously in \textsc{Maple 2018.0}. Both sets of code can be downloaded from the author's \textsc{GitHub} repository~\cite{code}. The author thanks the maintainers of the domainmatrix \textsc{SymPy} module for providing tools for efficient computation and The Ohio State University for access to the unity high performance computing cluster.

\pdfbookmark[1]{References}{ref}
\LastPageEnding

\end{document}